\documentclass{olplainarticle}

\usepackage{amssymb}
\usepackage{amsmath}
\usepackage{placeins}
\usepackage{graphicx}
\usepackage[colorlinks=true, linkcolor=blue, citecolor=blue, urlcolor=blue]{hyperref}
\usepackage{gensymb}
\usepackage{caption,subcaption}
\usepackage{tabularray}
\usepackage[acronym, nonumberlist]{glossaries}
\usepackage{xcolor}
\usepackage{graphicx}
\usepackage[T1]{fontenc}   
\usepackage[utf8]{inputenc} 
\usepackage{titlecaps}

\usepackage{bm}
\usepackage{tabularx}

\usepackage{comment}

\usepackage[linesnumbered,ruled,vlined]{algorithm2e}
\usepackage{multirow}
\usepackage{makecell}
\usepackage{color,soul}
\usepackage{todonotes}
\usepackage{siunitx}

\usepackage{pdflscape}
\usepackage{booktabs} 

\title{One-Way Thermo-Mechanical Coupled System Identification Using Displacement and Temperature Measurements}

\author[1,3]{Talhah Ansari}
\author[1]{Suneth Warnakulasuriya}
\author[1,4]{Ihar Antonau}
\author[2]{Harbir Antil}
\author[2,3]{Rainald Löhner}
\author[1]{Roland Wüchner}

\affil[1]{Chair of Structural Analysis, Technical University of Munich, Germany}
\affil[2]{College of Science, George Mason University, Fairfax, Virginia, United States}
\affil[3]{Institute for Advanced Study, Technical University of Munich, Germany}
\affil[4]{Cluster of Excellence SE2A -- Sustainable and Energy-Efficient Aviation, Technical University of Braunschweig}

\keywords{System identification, high-fidelity, damage detection, thermal field reconstruction, adjoint-method, Structural Health Monitoring}

\begin{abstract}
Structural system identification in the presence of thermal loads is challenging, as unmeasured or poorly modeled thermal effects can mask or mimic damage, leading to unreliable conclusions. This work presents an optimization-driven, adjoint-based high-fidelity system identification framework for localizing structural weakness and recovering the temperature field in one-way thermo-mechanical coupled structures. The methodology builds upon a standard optimization formulation that minimizes weighted discrepancies between simulated responses and measured data from a sparse displacement and temperature sensor network. To account for thermal effects, two strategies are proposed: a \textit{monolithic} approach, which simultaneously identifies Young’s modulus and temperature distributions, and a \textit{partitioned} approach, which iteratively couples two inexact sub‑problems through a Gauss-Seidel type fixed-point scheme. 
The proposed approaches are evaluated using two numerical examples—a Plate With a Hole and a Footbridge model -- under linearly varying and localized thermal fields, and for different sensor layouts. Both approaches successfully recover the Young’s modulus and temperature distributions, even when sensor placement does not fully capture the underlying thermal trends. Compared with a constant‑temperature assumption and interpolation of the temperature field from sensor data, the proposed approach achieves the most accurate damage localization and temperature reconstruction. The largest gains occur when localized thermal features are poorly sampled by sensors, where interpolation and constant--temperature assumptions underperform. Furthermore, results show that the location of the temperature sensors is as influential as the number of sensors: well--placed sensors substantially improve identification, while additional sensors that miss critical thermal features provide limited benefit.
\end{abstract}

\begin{document}

\flushbottom
\maketitle
\thispagestyle{empty}
\noindent
\printglossary[type=\acronymtype]

\section*{Nomenclature}

{\renewcommand\arraystretch{1.0}
\noindent\begin{longtable*}{@{}l @{\quad=\quad} l@{}}
$J$  & Cost function  \\
$\mathcal{L}$ &    Lagrangian function \\
$\boldsymbol{\phi}$& system identification parameters \\
$\mathbf{E}$ & Young's modulus distribution\\
$\boldsymbol{\Delta \mathbf{T}}$   & temperature difference compared to the ambient \\
$\boldsymbol{\psi}$   & computed model response vector\\
$\boldsymbol{\psi}^{\text{meas}}$ & measured data at the sensor location \\
$\boldsymbol{\psi}^{\text{comp}}$ & model computed data at the sensor location \\
$\mathbf{u}^{\text{meas}}$ & measured displacement at the sensor location \\
$\mathbf{u}^{\text{comp}}$ & model computed displacement at the sensor location \\
$\Delta \text{T}^{\text{meas}}$ & measured temperature difference at the sensor location \\
$\Delta \text{T}^{\text{comp}}$ & model computed temperature difference at the sensor location \\
$\mathbf{u}, \tilde{\mathbf{u}}$   & displacement and corresponding adjoint variables \\
$\mathbf{R}$  & FE residual \\
$\underline{\mathbf{K}}$  & stiffness matrix \\
$\mathbf{f}_{\text{int}}$ & internal force vector \\
$\mathbf{f}_{\text{ext}}$ & external mechanical force vector\\
$\mathbf{f}_{\boldsymbol{\Delta \mathbf{T}}}$ & thermal load vector\\
$x$ & scalar value x \\
$\mathbf{x}$ & vector x \\
$\underline{\mathbf{x}}$ & matrix x \\
\end{longtable*}}

\section{Introduction}
\label{s:intro}

Throughout their service life, civil and mechanical structures undergo gradual changes in material and system properties due to aging mechanisms such as damage, corrosion, and fatigue. Capturing these changes in a timely and reliable manner is central to \acrfull{SHM} and to the construction of high‑fidelity digital twins that evolve alongside the physical asset \citep{HAntil_2024a}. The continued proliferation of sensors, improvements in sensing technology, and advances in numerical simulation enable increasingly accurate digital representations of structural behavior. In many practical settings, however, responses are driven by the interaction of multiple physical fields (e.g., thermo‑mechanical coupling), necessitating digital twins that consistently account for these interacting domains. A key step in such pipelines is \acrfull{SI}—the estimation of the current state of the structure, including localization of weakness/damage \citep{airaudo2023adjoint} and the update of environmental fields that influence the response.

Thermal loads are among the dominant environmental drivers under normal operation. Bridges, for instance, experience markedly non‑uniform temperatures through the depth of the deck and soffit due to differential solar exposure \citep{kulprapha2012structural}; long‑span systems can also develop longitudinal temperature gradients \citep{ma2023statistical}. For long‑term monitoring, temperature‑induced response changes can be comparable to—or even exceed—the effects of live loads or moderate damage, complicating the interpretation of measured displacements and strains \citep{glashier2024temperature, kromanis2014predicting, zhou2014summary, bayraktar2022long}. Consequently, \acrshort{SI} methods that ignore or oversimplify thermal loads risk attributing temperature‑driven changes to damage, leading to false positives, biased parameter estimates, or failure to localize genuine damage. In the context of digital twins, separating temperature‑induced variability from true material degradation is therefore crucial for actionable diagnostics and decision making.

While prior work has addressed temperature field identification as a standalone problem  \citep{ansari2025adjoint}, real‑world applications often demand simultaneous reconstruction of environmental fields and structural properties, because these quantities are entangled in the measured response. The combined inverse problem is intrinsically ill‑conditioned: sparse sensors limit observability; temperature and stiffness fields can compensate for one another in the mechanical response; and distinct parameter combinations may yield similar error levels, admitting misleading local minima. As a result, algorithms require stabilization strategies and careful coupling to avoid misidentification and divergence.

This work develops an optimization‑driven, adjoint‑based high-fidelity \acrshort{SI} framework for one‑way thermo‑ mechanical coupled systems that jointly localizes damage (via Young’s modulus reduction) and reconstructs the temperature field using displacement and temperature measurements from sparse sensors. We propose two complementary strategies within a standard error-minimization formulation: (i) a monolithic approach that updates Young’s modulus and temperature distributions simultaneously in a single optimization, and (ii) a partitioned approach that staggered‑updates the two fields via a Gauss–Seidel type fixed‑point coupling with intentionally loose sub‑problem convergence to promote stable, incremental progress. 

This paper is organized as follows. Section~\ref{s:literature review} provides a focused, yet non‑exhaustive, overview of thermal field estimation and the treatment of thermal effects in \acrshort{SHM}. Section~\ref{s:methodology} introduces the proposed \acrshort{SI} framework, including sensor normalization, the optimization formulation, regularization, and the basis for result comparison. Section~\ref{s:numerical} demonstrates the methodology through two numerical examples—a Plate With a Hole and a Footbridge \acrfull{FE} model—examined under different thermal fields and sensor configurations. Section~\ref{s:conclusion} summarizes the main findings and outlines directions for future research.

\section{Literature Review}
\label{s:literature review}
Because temperature profoundly influences structural response, extensive efforts have been devoted to modeling and reconstructing thermal fields in structures. In general, thermal fields in civil structures can be reconstructed through three main routes \citep{lin20213d}. The most direct is measurement‑based sensing, including point thermometers, infrared thermography, and distributed temperature sensing using optical fibers, which provides dense measurements along fiber paths. A second route is analytical heat‑conduction solutions, suitable mainly for simple geometries \citep{santillan2015new}. The third route is to numerically simulate the temperature profile, commonly using \acrshort{FE} modeling  \citep{bui2019evaluation, chen2019temperature, yang2012fem}.

A widely used strategy is to spatially interpolate sparse temperature measurements using inverse distance weighting, kriging, natural neighbor interpolation, or spline methods. Kriging has been widely applied to dams and massive concrete structures, where 2D/3D thermal fields are reconstructed from limited sensor measurements or optical fibers \citep{zheng2019simulation, zhou2019temperature}. Interpolation quality depends heavily on sensor layout, motivating sensor positioning studies such as \cite{peng2020positioning}. More recently, learning‑based methods have been used to integrate climatic variables, e.g., convolutional neural network‑based reconstruction of dam surface temperatures under solar radiation and weather conditions \citep{pan2022novel}.
Field observations show that structures exposed to uneven solar radiation develop non‑linear cross‑sectional temperature profiles, both vertically and longitudinally \citep{kulprapha2012structural, fan2022efficient}. These spatially complex patterns highlight the limitations of simplified thermal assumptions. During construction of large concrete structures, temperature control is equally important: hydration heat, boundary conditions, and cooling strategies influence temperature evolution and can induce cracking \citep{lin20213d, liu2015precise}. Thermal regimes also vary significantly across different boundary zones, such as exposed versus submerged dam faces \citep{bui2019evaluation}.

Vibration‑based \acrshort{SHM} has been extensively investigated for temperature effects on damage detection \citep{bhuyan2019vibration, tefera2023challenges}. Natural frequencies are known to be highly temperature‑sensitive, while higher‑order mode shapes, though harder to extract, carry more damage information \citep{kostic2017vibration}.
\cite{jin2015structural} observed a negative correlation between the first seven eigenfrequencies and temperature, with more sensitivity observed for temperatures below the freezing point.
The main challenge is that real thermal fields are non‑uniform and time‑varying, and environmental changes (e.g., freezing, humidity, wind) can modify boundary conditions and material properties, often producing response variations greater than those caused by moderate damage \citep{glashier2024temperature, kromanis2014predicting}.

Dimensionality‑reduction approaches, such as Principal Component Analysis (PCA) and its variants, treat environmental effects as latent variables during damage detection \citep{nguyen2014damage, reynders2014output}. Hybrid PCA–Gaussian Mixture Model strategies further improve robustness under varying temperature states \citep{wang2020damage}.
Regression‑based and machine‑learning methods are widely used when historical monitoring data are available. Linear regression, artificial neural network (ANN), support vector regression, and autoregressive and exogenous input (ARX) time‑series models have been studied to determine the effect of temperature on the system properties \citep{peeters2001vibration, jin2016damage, jin2015structural}. Studies such as \cite{ni2009generalization} and \cite{kromanis2014predicting} established that data‑driven models can reliably reproduce thermal response trends, while more advanced frameworks, such as extended Kalman filter-based ANNs \citep{jin2016damage} and time-series/autoregressive-based ANN hybrids \citep{kostic2017vibration,huang2020damage}, further enhance compensation capability.

Several studies handle temperature via baseline‑oriented strategies, where damage is inferred as deviations from these baselines. \cite{yarnold2015temperature} advocated a temperature‑based \acrshort{SHM} paradigm that builds response--temperature baselines and flags damage as deviations, allowing for possible nonlinearity in the response–temperature relation. \cite{soo2020new} proposed a single baseline temperature profile to separate global thermal trends from local stiffness changes, estimating total elemental stiffness change and subtracting the global thermal component so that peaks in the residual indicate damage.

Beyond baselines, some works incorporate temperature within the identification variables.  \cite{meruane2012structural} formulate an inverse problem that jointly estimates stiffness and simplified temperature parameters (bridge end temperatures in a linearly varying field) using correlated mode shapes and eigenfrequencies, highlighting the coupling between thermal and material parameters. \cite{huang2019damage} embed temperature directly into material properties (stiffness-temperature relation) and use a hybrid particle swarm–cuckoo search to jointly identify damage and temperature effects based on a composite natural frequency, mode shape, and modal strain energy cost function, noting that uniform temperature assumptions reduce accuracy and recommending gradient‑type profiles (e.g., American Association of State Highway and Transportation Officials (AASHTO)‑style) as closer to real working states.

Despite substantial progress, many treatments of thermal effects still rely on simplified representations (uniform or code‑prescribed gradients), interpolation from sparse sensors, or statistical compensation in feature space, which may be insufficient when spatial heterogeneity and limited sensor coverage obscure critical thermal features. To address this more generally, we propose an optimization‑driven adjoint framework for one‑way thermo‑mechanical coupling that jointly identifies the Young’s modulus and temperature distributions from sparse displacement and temperature measurements. 
The system identification is posed as minimizing a weighted error functional that aggregates discrepancies between measured and simulated responses. We consider two formulations: a monolithic approach that updates both fields simultaneously in one global optimization, and a partitioned approach that splits the problem into two inexact sub‑optimizations (one for Young’s modulus, one for temperature) coupled via a Gauss–Seidel‑type fixed‑point iteration with loose sub‑problem convergence—i.e., early termination to allow small, incremental corrections that stabilize the coupled updates. To mitigate ill‑conditioning (few measurements vs. many design variables and the tendency of temperature and stiffness to compensate each other), we employ \acrfull{VM} filtering to smooth variables and gradients, suppressing spurious oscillations while retaining large‑scale structure. The underlying \acrshort{SI} formulation is consistent with the generic \acrshort{SI} used for damage detection under static and transient loads, as well as model calibration, such as constitutive parameters, boundary conditions, and load identification \citep{Ansari2025AdjointMembranes,ansari2025modelqual,ansari2026load,ansari2026transient},  thereby providing a unified and extensible approach.

\section{Methodology}
\label{s:methodology}

This section describes the proposed \acrshort{SI} methodology for identifying the Young's modulus and thermal field distributions of a structure from displacement and temperature measurements at limited discrete locations for one-way thermo-mechanical coupled systems.  
Firstly, the canonical system identification formulation is described. 
Section \ref{s:dd_with_temp_ident} then expands on the proposed methodology of obtaining both Young's modulus and temperature distributions by specializing the canonical \acrshort{SI} form. 
Sections \ref{s:sensor_normalization}, \ref{s:vertex_morphing}, and \ref{s:optimization} describe the cost function non-dimensionalization, \acrshort{VM} regularization, and the optimization procedure, respectively.
Lastly, section \ref{sec:comparison} discusses how the results  evaluation criteria and comparison metrics.

Let $\boldsymbol{\phi}(\mathbf{x})$ be the spatial distribution of the design variable (also referred to here as the control variable) throughout the structure. For a discretized system, determining the $\boldsymbol{\phi}(\mathbf{x})$ field from measurements may be posed as an optimization problem: Given $n$ test (load) cases and $l$ sensor types with each recording measurements at $m$ measuring locations $\mathbf{x}_j$, $j = 1,\dots, m$ of their respective responses (such as displacement, strains, temperature, etc) $\boldsymbol{\psi}_{kij}^{\text{meas}}$, $k=1,l$, $i = 1, n$, $j = 1, m$; determine the spatial distribution of $\boldsymbol{\phi}(\mathbf{x})$ that minimizes the errors between the sensor measurements and the model response. For a discretized system, the canonical form of the cost function for system identification can be written as:
\begin{eqnarray}
\min_{\boldsymbol{\phi}} \quad J(\boldsymbol{\psi}(\boldsymbol{\phi})) =\frac{1}{2} \sum_{k=1}^{l} \Omega_{k} \sum_{i=1}^{n} \zeta_{ki} \sum_{j=1}^{m} \omega_{kij} \left\lVert \boldsymbol{\psi}_{kij}^{\text{comp}}(\boldsymbol{\phi}) - \boldsymbol{\psi}_{kij}^{\text{meas}} \right\rVert_{2}^{2} \quad, 
\label{eq:cost_function}
\end{eqnarray}
where, 
$k=1,\dots,l$ refers to the type of sensor, such as displacement, strain, temperature sensors; $i=1,\dots,n$ refers to the $n$ given load cases; $j=1,\dots,m$ refers to the $m$ measuring locations; $\boldsymbol{\psi}_{kij}^{\text{meas}}$ and $\boldsymbol{\psi}_{kij}^{\text{comp}}$ are the measured and model computed quantities of interest at the sensor locations; $\Omega_{k}$ are individual sensor type cost function weights; $\zeta_{ki}$ are the individual load case weights; and $\omega_{kij}$ are the individual sensor weights.

\subsection{One-Way Thermo-Mechanical Coupled System Identification}
\label{s:dd_with_temp_ident}

This section proposes a combined multiple-field (Young's modulus and temperature) system identification method for one-way thermo-mechanical coupled structural systems using deformation and temperature measurements at limited discrete locations. 
The main motivation for this multi-field \acrshort{SI} is that structures may be exposed to thermal effects which may not be exactly known throughout the structure. Therefore, some temperature measurements from the structure may be recorded and incorporated into the \acrshort{SI} formulation to improve estimates of the Young's modulus distribution (i.e., damage localization) and the temperature distribution (i.e., thermal load).

Let $\text{E}(\mathbf{x})$ and $\Delta \text{T}(\mathbf{x})$  be the spatial distribution of the Young's modulus and temperature throughout the structure. In this work, the temperature difference compared to the ambient is used instead of absolute temperatures, hence the $\Delta$ symbol in the variable $\Delta \text{T}$. For the discretized system considered, $\mathbf{E}$ is the vector containing the elemental Young's modulus and $\boldsymbol{\Delta \textbf{T}}$ is the vector containing the nodal temperature differences compared to the ambient. For brevity, this vector $\boldsymbol{\Delta \textbf{T}}$ may be interchangeably referred to as the temperature field throughout this work.
Under thermal effects consideration, the finite element residual ($\mathbf{R}$) can be written as:
\begin{align}\label{eq:residual_equations_temp}
    \mathbf{R}(\mathbf{u}(\mathbf{E},\boldsymbol{\Delta \mathbf{T}}), \mathbf{E},\boldsymbol{\Delta \mathbf{T}}) =  \mathbf{f}_{\text{int}} 
    - \mathbf{f}_{\text{ext}} - \mathbf{f}_{\boldsymbol{\Delta \mathbf{T}}}  = \mathbf{0} \quad ,
\end{align}
where $\mathbf{f}_{\text{int}} = \underline{\mathbf{K}} \cdot \mathbf{u} $ for linear elastic with small deformation cases considered herein.
The Lagrangian, considering the residual, can be expressed as:
\begin{equation}
    \mathcal{L}(\mathbf{u}(\mathbf{E},\boldsymbol{\Delta \mathbf{T}}),\mathbf{E}, \boldsymbol{\Delta \mathbf{T}}, \mathbf{\tilde{u}}) = J(\mathbf{u}(\mathbf{E},\boldsymbol{\Delta \mathbf{T}}),\boldsymbol{\Delta \mathbf{T}}) + \mathbf{\tilde{u}}^\top \cdot \mathbf{R}(\mathbf{u}(\mathbf{E},\boldsymbol{\Delta \mathbf{T}}),\mathbf{E},\boldsymbol{\Delta \mathbf{T}})
    \quad,
\label{eq:lagrangian_temp}
\end{equation}
where $\mathbf{\tilde{u}}$ are the Lagrange multipliers (adjoint variables).
The following sub-sections discuss two different approaches to achieve this two-field (Young's modulus and temperature distributions) \acrshort{SI}. The first is the `monolithic approach', where Young's modulus distribution (also referred to here as `weakness' or `damage') and the temperature distributions are identified simultaneously in a monolithic \acrshort{SI}.
The second is a partitioned scheme in which the Young's modulus and temperature distributions are identified via a Gauss-Seidel iterative method. 

For both approaches, a simplified formulation considering one load case, i.e., $n=1$ and $\zeta=1$, shall be presented to improve readability of the equations and focus on demonstrating the key aspects of the methodology. However, there are two groups/types of sensors: displacement sensors and temperature sensors. In general case, the method's formulation can be extended to multiple loads/sensor types.

\subsubsection{Monolithic Approach}
\label{s:dd_with_temp_ident_mono}

This section describes the \textit{monolithic} approach for weakness detection and thermal field identification using displacement and temperature sensors. The control variables are the Young's modulus distribution ($\mathbf{E}$) and the temperature distribution ($\boldsymbol{\Delta \mathbf{T}}$).
Given displacement sensor measurements $\mathbf{u}_{j}^{\text{meas}}$ at $m$ measuring locations $\mathbf{x}_j$, $j = 1,\dots, m$ and temperature sensor measurements $\Delta\text{T}_{p}^{\text{meas}}$ at $o$ measuring locations $\mathbf{x}_p$, $p = 1,\dots, o$, determine the spatial distribution of $\mathbf{E}$ and $\boldsymbol{\Delta \text{T}}$ that minimizes the errors between the sensor measurements and the model response. The canonical form in Eq.~\eqref{eq:cost_function} can be modified to obtain the composite cost function formulated as:
\begin{align}
\min_{\mathbf{E}, \boldsymbol{\Delta \mathbf{T}}} \quad
J(\mathbf{u}(\mathbf{E},\boldsymbol{\Delta \mathbf{T}}),
  \boldsymbol{\Delta \mathbf{T}})
\notag 
& = 
\underbrace{\frac{1}{2} \Omega_1 
\sum_{j=1}^{m} \omega_{j}
\left\lVert
\mathbf{u}_{j}^{\text{comp}}(\mathbf{E},\boldsymbol{\Delta \mathbf{T}})
- \mathbf{u}_{j}^{\text{meas}}
\right\rVert_2^2}_{J_{D}} 
\notag \\[0.0ex]
& \quad +
\underbrace{\frac{1}{2} \Omega_2 
\sum_{p=1}^{o} \omega_{p}
\left\lVert
\Delta \text{T}_{p}^{\text{comp}}
- \Delta \text{T}_{p}^{\text{meas}}
\right\rVert_2^2}_{J_{T}} .
\label{eq:cost_function_monolithic}
\end{align}
where, $j=1,\dots,m$ refers to the $m$ measuring locations of $\mathbf{u}_{j}^{\text{meas}}$; $p=1,\dots,o$ refers to the $p$ measuring locations;  $\Delta \text{T}_{p}^{\text{meas}}$ and $\Delta \text{T}_{p}^{\text{comp}}$ are the measured and model computed temperature differences at the sensor locations; $\omega_{j}$, $\omega_{p}$ are the individual displacement and temperature sensor weights; and $\Omega_{1}$, $\Omega_{2}$ are the displacement and temperature response weights.

The variation of the Lagrangian (Eq.~\eqref{eq:residual_equations_temp}) yields: the `forward' or `primal' problem ($\frac{\partial \mathcal{L}}{\partial \mathbf{\tilde{u}}}=\mathbf{R}=\mathbf{0}$) and the `adjoint' problem ($\frac{\partial \mathcal{L}}{\partial \mathbf{u}}=\mathbf{0}$). The adjoint system solving for the adjoint variables ($\tilde{\mathbf{u}}$) is expressed as:
\begin{equation}
    \left[ \frac{\partial \mathbf{R}}{\partial \mathbf{u}} \right]^{\top} \cdot \mathbf{\tilde{u}} = -\frac{\partial J}{\partial \mathbf{u}} \implies \underline{\mathbf{K}}^{\top} \cdot \mathbf{\tilde{u}} = -\frac{\partial J}{\partial \mathbf{u}} \quad.
\label{eq:adjoint_problem}
\end{equation}
For most structural systems $\underline{\mathbf{K}}^{\top}=\underline{\mathbf{K}}$, and thus can be reused from the primal solve.
It is noted here that the temperature differences ($\Delta\text{T}$) considered in this work are considered to be small, such that the Young's modulus of the material is independent of temperature. This presumption would not naturally hold for large temperature changes, where Young's modulus depends on temperature and the Young's modulus-temperature relation must be accounted for.

Considering this simplified assumption, using the residual and Lagrangian formulations from Eqs.~\eqref{eq:residual_equations_temp},~\eqref{eq:lagrangian_temp}, the gradient with respect to the Young's modulus distribution can be written as:
\begin{align}
\frac{d\mathcal{L}}{d\mathbf{E}}  
    &=  \frac{\partial \mathcal{L}}{\partial \mathbf{E}} 
    +  \frac{\partial \mathcal{L}}{\partial \mathbf{u}} \cdot  \frac{d \mathbf{u}}{d \mathbf{E}} 
    +  \frac{\partial \mathcal{L}}{\partial \boldsymbol{\Delta \mathbf{T}}} \cdot  \frac{d \boldsymbol{\Delta \mathbf{T}}}{d  \mathbf{E}}  
    \nonumber 
    \\
    \begin{split}
    &=  \frac{\partial J}{\partial \mathbf{E}} + \mathbf{\tilde{u}}^\top  \cdot \frac{\partial \mathbf{R}}{\partial \mathbf{E}} +  \underbrace{ \left( \frac{\partial J}{\partial \mathbf{u}} + \mathbf{\tilde{u}}^\top \cdot \frac{\partial \mathbf{R}}{\partial \mathbf{u}} \right) }_{\text{Adjoint system}}   \cdot  \frac{d \mathbf{u}}{d \mathbf{E}} 
    \\
    &\qquad \qquad \qquad \qquad+ \left( \frac{\partial J}{\partial\boldsymbol{\Delta \mathbf{T}}} + \mathbf{\tilde{u}}^\top \cdot \frac{\partial \mathbf{R}}{\partial \boldsymbol{\Delta \mathbf{T}}} \right)  \cdot  \frac{d \boldsymbol{\Delta \mathbf{T}}}{d \mathbf{E}}   
    \end{split}
    \nonumber 
    \\
    &= \frac{\partial J}{\partial \mathbf{E}} +  \mathbf{\tilde{u}}^\top  \cdot \frac{\partial \mathbf{R}}{\partial \mathbf{E}}  \quad.
\label{eq:gradient_monolithic_E}
\end{align}
Along similar lines, the gradient with respect to the temperature field can be written as:
\begin{align}
\frac{d\mathcal{L}}{d\boldsymbol{\Delta \mathbf{T}}} 
    &=  \frac{\partial \mathcal{L}}{\partial \boldsymbol{\Delta \mathbf{T}}} 
    +  \frac{\partial \mathcal{L}}{\partial \mathbf{u}} \cdot  \frac{d \mathbf{u}}{d \boldsymbol{\Delta \mathbf{T}}} 
    +  \frac{\partial \mathcal{L}}{\partial\mathbf{E}} \cdot  \frac{d \mathbf{E}}{d \boldsymbol{\Delta \mathbf{T}}}  
    \nonumber 
    \\
    \begin{split}
    &=  \frac{\partial J}{\partial \boldsymbol{\Delta \mathbf{T}}} + \mathbf{\tilde{u}}^\top  \cdot \frac{\partial \mathbf{R}}{\partial \boldsymbol{\Delta \mathbf{T}}} +  \underbrace{ \left( \frac{\partial J}{\partial \mathbf{u}} + \mathbf{\tilde{u}}^\top \cdot \frac{\partial \mathbf{R}}{\partial \mathbf{u}} \right) }_{\text{Adjoint system}}   \cdot  \frac{d \mathbf{u}}{d \boldsymbol{\Delta \mathbf{T}}} 
    \\
    &\qquad \qquad \qquad \qquad \qquad+ \left( \frac{\partial J}{\partial\mathbf{E}} + \mathbf{\tilde{u}}^\top \cdot \frac{\partial \mathbf{R}}{\partial \mathbf{E}} \right)  \cdot  \frac{d \mathbf{E}}{d \boldsymbol{\Delta \mathbf{T}}}   
    \end{split}
    \nonumber 
    \\
    &= \frac{\partial J}{\partial \boldsymbol{\Delta \mathbf{T}}} +  \mathbf{\tilde{u}}^\top  \cdot \frac{\partial \mathbf{R}}{\partial \boldsymbol{\Delta \mathbf{T}}}  \quad.
\label{eq:gradient_monolithic_T}
\end{align}
%
In Eq.~\eqref{eq:gradient_monolithic_E}: $\frac{d\mathcal{L}}{d\mathbf{E}}=\frac{dJ}{d\mathbf{E}}$, and in Eq.~\eqref{eq:gradient_monolithic_T}: $\frac{d\mathcal{L}}{d\boldsymbol{\Delta\mathbf{T}}}=\frac{dJ}{d\boldsymbol{\Delta\mathbf{T}}}$, because at equilibrium the residual $\mathbf{R}=\mathbf{0}$ and so does its change (Ref. \cite{antil2018brief}).
As explained earlier, the cases considered in this study work under the assumption that the Young's modulus and temperature are independent fields, thus the terms $\frac{d\boldsymbol{\Delta \mathbf{T}}}{d\mathbf{E}}$ and $\frac{d\mathbf{E}}{d\boldsymbol{\Delta \mathbf{T}}}$ are dropped from Eqs.~\eqref{eq:gradient_monolithic_E},~\eqref{eq:gradient_monolithic_T}. 

The adjoint variables computed in Eq.~\eqref{eq:adjoint_problem} are inserted in Eqs.~\eqref{eq:gradient_monolithic_E},~\eqref{eq:gradient_monolithic_T} to compute the gradients.
The partial sensitivity term $\frac{\partial J}{\partial  \mathbf{E}} = \mathbf{0}$ because the cost function (Eq.~\eqref{eq:cost_function_monolithic}) does not have a direct dependence on $\mathbf{E}$. 
However, in the gradient w.r.t. $\boldsymbol{\Delta \mathbf{T}}$, the partial sensitivity term $\frac{\partial J}{\partial \boldsymbol{\Delta \mathbf{T}}} \neq \mathbf{0}$ because the cost function (Eq.~\eqref{eq:cost_function_monolithic}) has a direct dependence on $\boldsymbol{\Delta \mathbf{T}}$.
In terms of computational resources, the adjoint method computes gradients with $\mathcal{O}(1)$ effort, compared to $\mathcal{O}(n_\text{design variables})$ with direct sensitivity.

\FloatBarrier
\subsubsection{Partitioned Approach}
\label{s:dd_with_temp_ident_part}

This section describes the \textit{partitioned} approach for weakness detection and thermal field identification using displacement and temperature sensors. 
The main idea is to optimize $\mathbf{E}$ and $\boldsymbol{\Delta \mathbf{T}}$ in two separate optimizations.
A simple two‑step sequential scheme alternates once between the two fields: first identify $\boldsymbol{\Delta \mathbf{T}}$ with $\mathbf{E}$ held fixed, then identify $\mathbf{E}$ with $\boldsymbol{\Delta \mathbf{T}}$ held fixed. In ill‑conditioned settings, this highly path‑dependent convergence may hinge on the initial guess and the update order, and can yield non‑physical or spurious solutions. For example, when $\mathbf{E}$ is fixed, the optimizer can compensate displacement errors in damaged regions by introducing localized heating/cooling in $\boldsymbol{\Delta \mathbf{T}}$ (unless temperature sensors directly constrain those zones). The subsequent $\mathbf{E}$ update then inherits this bias, making accurate weakness identification difficult.

To robustly couple the two fields, we propose a Gauss–Seidel-type fixed-point iteration outer loop encompassing the two sub‑optimizations. Each sub‑problem is solved inexactly with loose convergence (early termination), so that updates remain incremental while the other field is still imperfect. This prevents either sub‑problem from over‑fitting inconsistencies in the counterpart field and promotes stable, simultaneous correction. Small successive refinements across coupling iterations guide the combined system toward a consistent solution, whereas driving each sub‑problem to full convergence in isolation tends to over‑compensate and can diverge.

Figure~\ref{fig:Figure_1} summarizes the workflow: sub‑problem ‘A’ updates $\boldsymbol{\Delta \mathbf{T}}$ with $\mathbf{E}$ fixed; sub‑problem ‘B’ updates $\mathbf{E}$ with $\boldsymbol{\Delta \mathbf{T}}$ fixed. Eqs.~\eqref{eq:partitioned_T_1}–\eqref{eq:partitioned_T_3} define optimization `A’ (objective $J=J_D+J_T$ since $\boldsymbol{\Delta \mathbf{T}}$ affects both displacement via thermal expansion and temperature sensor errors), and Eqs.~\eqref{eq:partitioned_E_1}–\eqref{eq:partitioned_E_3} define optimization `B’ (objective $J=J_D$ under the assumption that $\mathbf{E}$ does not affect $J_T$). After each sub‑optimization is solved, an optional relaxation with parameter $\beta$ is applied (Eqs.~\eqref{eq:relaxation_1}–\eqref{eq:relaxation_2}), and the outer loop is checked against the overall convergence criterion in Eq.~\eqref{eq:gs_convergence}. 

Since the underlying PDE‑constrained problem is unchanged, the gradients coincide with the monolithic case (Eqs.~\eqref{eq:gradient_monolithic_E}, \eqref{eq:gradient_monolithic_T}). While both monolithic and partitioned schemes have $\mathcal{O}(1)$ cost per \textit{optimization iteration} for gradient computation, the partitioned approach incurs higher computational cost due to multiple sub‑optimization solves per coupling iteration.

\begin{figure}[t]
\centering
\fbox{%
\begin{minipage}{0.95\linewidth}
\textbf{Subproblem A: Optimization w.r.t. field $\boldsymbol{\Delta \mathbf{T}}$ (fixed $\mathbf{E}^{(r)}$)}  \vspace{-0.4em}
\begin{align}
\widehat{\boldsymbol{\Delta \mathbf{T}}}^{(r+1)} &\approx \arg\min_{\boldsymbol{\Delta \mathbf{T}}} \;
J(\mathbf{E}^{(r)}, \boldsymbol{\Delta \mathbf{T}})
\label{eq:partitioned_T_1} \\[0.8ex]
J(\mathbf{E}^{(r)}, \boldsymbol{\Delta \mathbf{T}}) &=
\underbrace{\frac{1}{2} \Omega_1  \sum_{j=1}^{m} \omega_{j}
\left\lVert  \mathbf{u}_{j}^{\text{comp}}(\mathbf{E}^{(r)},\boldsymbol{\Delta \mathbf{T}}) - \mathbf{u}_{j}^{\text{meas}}
\right\rVert_2^2}_{J_{D}}
\notag \\[0.5ex]
&\quad+ \underbrace{\frac{1}{2} \Omega_2  \sum_{p=1}^{o} \omega_{p}
\left\lVert \Delta \text{T}_{p}^{\text{comp}}
- \Delta \text{T}_{p}^{\text{meas}} \right\rVert_2^2}_{J_{T}} ,
\label{eq:partitioned_T_2} \\[0.8ex]
\frac{d\mathcal{L}}{d\boldsymbol{\Delta \mathbf{T}}} = \frac{dJ}{d\boldsymbol{\Delta \mathbf{T}}} &=
\frac{\partial J}{\partial \boldsymbol{\Delta \mathbf{T}}} + \tilde{\mathbf{u}}^{\top} \frac{\partial \mathbf{R}}{\partial \boldsymbol{\Delta \mathbf{T}}} \quad . \label{eq:partitioned_T_3}
\end{align}
\emph{Inner solve terminated early using a loose convergence criterion.}
\end{minipage}}

\vspace{0.9em}
\begin{minipage}{0.95\linewidth}
\qquad \qquad \qquad \qquad \qquad \qquad \qquad \qquad\quad $\Downarrow$ \quad
\textit{Gauss--Seidel relaxation}
\quad $\Uparrow$
\vspace{-0.3em}
\begin{align}
     \mathbf{E}^{(r+1)} &= \beta \, \widehat{\mathbf{E}}^{(r+1)}
+ (1-\beta)\, \mathbf{E}^{(r)}, \quad 0 < \beta < 2 \label{eq:relaxation_1} \\
 \boldsymbol{\Delta \mathbf{T}}^{(r+1)} &= \beta \, \widehat{\boldsymbol{\Delta \mathbf{T}}}^{(r+1)} + (1-\beta)\, \boldsymbol{\Delta \mathbf{T}}^{(r)},
\quad 0 < \beta < 2  \label{eq:relaxation_2}
\end{align}
\vspace{-1.1em}
\end{minipage}
\vspace{1em}
\fbox{%
\begin{minipage}{0.95\linewidth}
\textbf{Subproblem B: Optimization w.r.t. field $\mathbf{E}$ (fixed $\boldsymbol{\Delta \mathbf{T}}^{(r+1)}$)} \vspace{-0.4em}
\begin{align}
\widehat{\mathbf{E}}^{(r+1)} &\approx \arg\min_{\mathbf{E}} \;
J_{D}(\mathbf{E}, \boldsymbol{\Delta \mathbf{T}}^{(r+1)})
\label{eq:partitioned_E_1} \\[0.8ex]
J_D(\mathbf{E}, \boldsymbol{\Delta \mathbf{T}}^{(r+1)}) &=
\frac{1}{2} \Omega_1 \sum_{j=1}^{m} \omega_{j} \left\lVert 
 \mathbf{u}_{j}^{\text{comp}}(\mathbf{E},\boldsymbol{\Delta \mathbf{T}}^{(r+1)}) - \mathbf{u}_{j}^{\text{meas}} \right\rVert_2^2 ,
\label{eq:partitioned_E_2} \\[0.8ex]
\frac{d\mathcal{L}}{d\mathbf{E}} = \frac{dJ}{d \mathbf{E}} &= \frac{\partial J_D}{\partial \mathbf{E}} + \tilde{\mathbf{u}}^{\top} \frac{\partial \mathbf{R}}{\partial \mathbf{E}} \quad . \label{eq:partitioned_E_3}
\end{align}
\emph{Inner solve terminated early using a loose convergence criterion.} 
\end{minipage}
}
\fbox{%
\begin{minipage}{0.95\linewidth}
\textbf{Outer loop convergence}
\vspace{-0.8em}
\begin{equation}
    J = J_D + J_T \le J_{\text{target}} \quad \text{or} \quad r \ge r_{\max} \label{eq:gs_convergence}
\end{equation}

\end{minipage}}
\caption{Gauss--Seidel--type partitioned fixed point optimization with inexact subproblem solves, where $r$ refers to the coupling iteration and $\beta$ refers to the relaxation parameter.}
\label{fig:Figure_1}
\end{figure}

\FloatBarrier

\subsection{Sensor Normalization}
\label{s:sensor_normalization}

From the cost functions in Eqs.~\eqref{eq:cost_function},~\eqref{eq:cost_function_monolithic}, it can be noted that different types of sensors have different dimensions (m, \degree C, dimensionless, etc) and orders of magnitude of sensor measured values and sensor errors. For example, displacements might be in the order of $10^{-6}-10^{-1}$ m, while temperatures might be in the order of $10^{0}-10^{2}$ \degree C. To avoid the composite cost function from favoring one set of sensors and thereby one design variable, it is preferred to normalize and non-dimensionalize the cost function.
\cite{airaudo2023adjoint} presented several methods for sensor normalization, including local, average, maximum, and local/maximum weighting. As this work is not focused on sensor normalization, the `maximum measured value' normalization is selected and consistently used throughout to ensure consistent comparisons among the results. It is written as:
\begin{align}
    \mathrm{u}_{\text{max}}^{\text{meas}} &= \text{max} (\left\lVert\mathbf{u}_{j}\right\rVert_{\infty}, j=1,\dots,m) \quad &; \quad \omega_{j} =& \frac{1}{(\mathrm{u}_{\text{max}}^{\text{meas}})^{2}} \quad , \label{eq:sensor_weight_a} \\
    \Delta \text{T}_{\text{max}}^{\text{meas}} &= \text{max} (\left\lVert\Delta \text{T}_{p}\right\rVert_{\infty}, p=1,\dots,o) \quad &; \quad \omega_{p} =& \frac{1}{(\Delta \text{T}_{\text{max}}^{\text{meas}})^{2}} \quad ,
    \label{eq:sensor_weight_b}
\end{align}
where $j=1,\dots,m$ refers to the $m$ displacement sensors, and $p=1,\dots,o$ refers to the $o$ temperature sensors. 
In general, the sensor normalizing factor may be carried for each summation in the canonical form of the cost function (Eq.~\eqref{eq:cost_function}), i.e., per sensor type, per load case, and per sensor.  
In this work, since only one mechanical load case ($n=1$) is considered, and maximum weighting is used, the above equations are obtained.

\FloatBarrier
\subsection{Vertex Morphing}
\label{s:vertex_morphing}

This section presents a brief description of the vertex morphing (\acrshort{VM}) smoothing method, which is commonly used in node-based shape optimization for smooth shape updates \citep{hojjat2014vertex}. The central idea is to identify the physical design field that is controlled by a non-physical control field via a mapping function.

However, before \acrshort{VM} is applied, the Sigmoid function is used to project the control variables (here, the elemental Young's modulus and nodal temperature differences) to a $[0,1]$ range, referred to as \textit{physical phi field} ($\boldsymbol{\Phi}$). 
\acrshort{VM} is then applied to this field. 
Eq.\eqref{eq:vm_forward_filter} defines the operation known as \textit{forward filtering}, in discretized form. 
\begin{equation}
    \boldsymbol{\Phi} = \underline{\mathbf{A}} \cdot \boldsymbol{\tilde{\Phi}} \qquad or \qquad \Phi_p = \text{A}_{pq} \cdot \tilde{\Phi}_q \qquad ,
\label{eq:vm_forward_filter}
\end{equation}
where $\boldsymbol{\Phi}$ is the vector of physical phi field, $\boldsymbol{\tilde{\Phi}}$ is the vector of corresponding control field,  and $\underline{\mathbf{A}}$ is the kernel function matrix mapping the control field $\boldsymbol{\tilde{\Phi}}$ to the physical field $\boldsymbol{\Phi}$. 
$\text{A}_{pq}$ represents the filtering interaction between node/element $p$ and node/element $q$ based on their positions and Euclidean distance. 
Although \acrshort{VM} supports different discretizations for the physical and control fields, here, the same discretization is used in both spaces.  

Eq.~\eqref{eq:vm_backward_filter} defines the operation known as \textit{backward filtering}, in discretized form, which maps the physical field gradients ($dJ/d\boldsymbol{\Phi}$) to the control field gradients ($dJ /d\boldsymbol{\tilde{\Phi}}$).
\begin{subequations}\label{eq:vm_backward_filter} 
    \begin{gather}
         \frac{d J}{d\tilde{\Phi}_p} = \frac{d J}{d \Phi_q} \cdot \frac{d\Phi_q}{d\tilde{\Phi}_p} \quad, 
         \label{eq:vm_backward_filter_a}\\
         \implies \frac{d J}{d\tilde{\Phi}_p} = \text{A}_{qp} \cdot \frac{d J}{ d \Phi_q} = \text{A}_{qp} \cdot \text{b}_{q} \quad \text{or} \quad  \frac{d J}{d\boldsymbol{\tilde{\Phi}}} = \underline{\mathbf{A}}^\top \cdot \mathbf{b} = \underline{\mathbf{A}}^\top \cdot \frac{d J}{d\boldsymbol{\Phi}} \quad,
         \label{eq:vm_backward_filter_b}
    \end{gather} 
\end{subequations}
where $\text{b}_q = \frac{d J}{d\Phi_q}$ is the gradient of cost function $J$ for the $q$-th physical field design variable. 

Since these high-fidelity \acrshort{SI} problems are highly ill-conditioned, i.e., a unique solution is generally not possible because different combinations of design variables can produce the same response, some form of regularization is often required. Thus, \acrshort{VM} acts as a smoothing function and can be understood as a convolution operator within a specified radius, thereby yielding smoother variable updates and smoother fields. 
The readers are referred to \cite{antonau2025comparison,antonau2022latest, ghantasala2021realization, hojjat2014vertex} for a detailed explanation of \acrshort{VM}. Simple point averaging, Laplacian smoothing, Pseudo-Laplacian smoothing, etc, are other gradient smoothing methods available in literature \citep{airaudo2023adjoint, lohner2024high}.

\subsection{Optimization}
\label{s:optimization}

\begin{algorithm}[H]
\setcounter{AlgoLine}{0}
\caption{Adjoint-based optimization in a single-loop (monolithic) system identification}
\label{alg:adjoint_optimization_mono}

\SetKwInput{KwData}{Start}

\KwData{Initial control variables $(\boldsymbol{\phi}^{(0)},\boldsymbol{\theta}^{(0)})$, iteration $r\gets0$}

\While{not converged ($r\le r_{\max}$ and $J>J_{\text{target}}$)}{

Forward analysis of all $n$ load cases:
$\mathbf{u}^{(r)}=[\mathbf{u}_1^{(r)},\dots,\mathbf{u}_n^{(r)}]$\;

Cost function calculation: $J^{(r)}$\;

Adjoint analysis of all $n$ load cases:
$\tilde{\mathbf{u}}^{(r)}=[\tilde{\mathbf{u}}_1^{(r)},\dots,\tilde{\mathbf{u}}_n^{(r)}]$\;

Gradient computation for each load case:
$\boldsymbol{\frac{d J}{d \boldsymbol{\phi}}}^{(r)} =  [ \frac{d J}{d \boldsymbol{\phi}_1}^{(r)} , \dots, \frac{d J}{d \boldsymbol{\phi}_n}^{(r)} ] $ , $\boldsymbol{\frac{d J}{d \boldsymbol{\theta}}}^{(r)} =  [ \frac{d J}{d \boldsymbol{\theta}_1}^{(r)} , \dots, \frac{d J}{d \boldsymbol{\theta}_n}^{(r)} ] $ \;

Aggregate and smoothen (if required):
$\left(\frac{d J}{d \boldsymbol{\phi}}\right)_{\text{smooth}}^{(r)}, \left(\frac{d J}{d \boldsymbol{\theta}}\right)_{\text{smooth}}^{(r)}$  \;

Update the control variables:
$ \boldsymbol{\phi}^{(r+1)}  = \boldsymbol{\phi}^{(r)}  -  \gamma \left(\frac{d J}{d \boldsymbol{\phi}}\right)_{\text{smooth}}^{(r)}$ , $ \boldsymbol{\theta}^{(r+1)}  = \boldsymbol{\theta}^{(r)}  -  \gamma \left(\frac{d J}{d \boldsymbol{\theta}}\right)_{\text{smooth}}^{(r)}$ \;

$r \gets r + 1$ \;

}

\end{algorithm}

\begin{algorithm}[H]
\setcounter{AlgoLine}{0}
\caption{Adjoint-based optimization in a Gauss-Seidel-based partitioned system identification}
\label{alg:adjoint_optimization_partitioned}
\SetKwInput{KwData}{Start}
\KwData{Initial control variables $(\boldsymbol{\phi}^{(0)}, \boldsymbol{\theta}^{(0)})$, iteration $r=r_{\text{A}} + r_{\text{B}} \gets 0$}
\While{\text{\acrshort{GS} loop not converged} ($r \le r_{\text{max}}$ and $J=J_{\text{A}}+J_{\text{B}} > J_{\text{target}}$)}{
\SetKwInput{KwData}{Start}
\KwData{sub-optimization `A': Initial control variables $\boldsymbol{\phi}^{(r_{\text{A}})}$, fixed $\boldsymbol{\theta}^{(r)}$, iteration counter $r_{\text{A}} \gets 0$}
\While{\text{not converged} ($r_{\text{A}} \le r_{\text{A,max}}$ and $J_{\text{A}} > J_{\text{A,target}}$)}{
     \text{Forward analysis of all $n$ load cases:} $ \mathbf{u}^{(r_{\text{A}})} = [\mathbf{u}_1^{(r_{\text{A}})}, \dots, \mathbf{u}_n^{(r_{\text{A}})}] $ \;
     
     \text{Cost function calculation:} $J_{\text{A}}^{(r_{\text{A}})}$ \;
     
    \text{Adjoint analysis of all $n$ load cases:} $\tilde{\mathbf{u}}^{(r_{\text{A}})} = \left[ \tilde{\mathbf{u}}_1^{(r_{\text{A}})}, \dots, \tilde{\mathbf{u}}_n^{(r_{\text{A}})} \right]$ \;
    
    \text{Gradient computation for each load case:} $\boldsymbol{\frac{d J}{d \boldsymbol{\phi}}}^{(r_{\text{A}})} =  [ \frac{d J}{d \boldsymbol{\phi}_1}^{(r_{\text{A}})} , \dots, \frac{d J}{d \boldsymbol{\phi}_n}^{(r_{\text{A}})} ] $ \;
    
    \text{Aggregate and smoothen (if required):} $\left(\frac{d J}{d \boldsymbol{\phi}}\right)_{\text{smooth}}^{(r_{\text{A}})}$ \;
    
    \text{Control variable update:}  $ \boldsymbol{\phi}^{(r_{\text{A}}+1)}  = \boldsymbol{\phi}^{(r_{\text{A}})}  -  \gamma \left(\frac{d J}{d \boldsymbol{\phi}}\right)_{\text{smooth}}^{(r_{\text{A}})}$ \;
    
    $r_A \gets r_A + 1$,  $r \gets r + 1$ , $\widehat{\boldsymbol{\phi}}^{(r)} \gets \boldsymbol{\phi}^{(r_{\text{A}})}$ \;
}
\SetKwInput{KwData}{Relaxation}
\KwData{$\boldsymbol{\phi}^{(r)} \gets \beta \cdot \widehat{\boldsymbol{\phi}}^{(r)} + (1-\beta)\cdot \boldsymbol{\phi}^{(r-1)} $}
\SetKwInput{KwData}{Start}
\KwData{sub-optimization `B': Initial control variables $\boldsymbol{\theta}^{(r)}$, fixed $\boldsymbol{\phi}^{(r)}$, iteration $r_{\text{B}} \gets 0$}
\While{\text{not converged} ($r_{\text{B}} \le r_{\text{B,max}}$ and $J_{\text{B}} > J_{\text{B,target}}$)}{
     \text{Forward analysis of all $n$ load cases:} $ \mathbf{u}^{(r_{\text{B}})} = [\mathbf{u}_1^{(r_{\text{B}})}, \dots, \mathbf{u}_n^{(r_{\text{B}})}] $ \;
     
     \text{Cost function calculation:} $J_{\text{B}}^{(r_{\text{B}})}$ \;
     
    \text{Adjoint analysis of all $n$ load cases:} $\tilde{\mathbf{u}}^{(r_{\text{B}})} = \left[ \tilde{\mathbf{u}}_1^{(r_{\text{B}})}, \dots, \tilde{\mathbf{u}}_n^{(r_{\text{B}})} \right]$ \;
    
    \text{Gradient computation for each load case:} $\boldsymbol{\frac{d J}{d \boldsymbol{\theta}}}^{(r_{\text{B}})} =  [ \frac{d J}{d \boldsymbol{\theta}_1}^{(r_{\text{B}})} , \dots, \frac{d J}{d \boldsymbol{\theta}_n}^{(r_{\text{B}})} ] $ \;
    
    \text{Aggregate and smoothen (if required):} $\left(\frac{d J}{d \boldsymbol{\theta}}\right)_{\text{smooth}}^{(r_{\text{B}})}$ \;
    
    \text{Control variable update:}  $ \boldsymbol{\theta}^{(r_{\text{B}}+1)}  = \boldsymbol{\theta}^{(r_{\text{B}})}  -  \gamma \left(\frac{d J}{d \boldsymbol{\theta}}\right)_{\text{smooth}}^{(r_{\text{B}})}$ \;
    
    $r_B \gets r_B + 1$ , $r \gets r + 1$, $\widehat{\boldsymbol{\theta}}^{(r)} \gets \boldsymbol{\theta}^{(r_{\text{B}})}$ \;
}
\SetKwInput{KwData}{Relaxation}
\KwData{$\boldsymbol{\theta}^{(r)} \gets \beta \cdot \widehat{\boldsymbol{\theta}}^{(r)} + (1-\beta)\cdot \boldsymbol{\theta}^{(r-1)} $}

}
\end{algorithm}

This section describes the optimization procedure used for system identification. The cost function and the adjoint-based gradients are computed as discussed in Section \ref{s:dd_with_temp_ident}. Algorithm \ref{alg:adjoint_optimization_mono} presents the optimization steps for the monolithic \acrshort{SI} approach, i.e., all control variables are updated within a single optimization loop.
Here, $\boldsymbol{\phi},\boldsymbol{\theta}$ are the two sets of control variables, $r$ is the iteration counter, $J$ is the cost function, and $\gamma$ is the step size. For all the cases, the search direction is scaled by the $L_2$ norm \citep{antonau2025comparison}.

Algorithm \ref{alg:adjoint_optimization_partitioned} presents the general optimization steps for the Gauss-Seidel-based partitioned \acrshort{SI} approach, i.e., each sub-optimization controls only one type of variable. Here, $\boldsymbol{\phi},\boldsymbol{\theta}$ are the two sets of control variables controlled by sub-optimizations `A' and `B', respectively. $r$ is the overall optimization iteration counter and keeps track of the sub-optimization iterations $r_{\text{A}}$ and $r_{\text{B}}$. $J_{\text{A}}$ and $J_{\text{B}}$ are the cost functions associated with sub-problems `A' and `B', respectively. $J$ is the cost function tracking overall convergence. $\beta$ is the relaxation parameter in the Gauss-Seidel loop. While the algorithm is presented in general form, the case-specific equations were given in Fig. \ref{fig:Figure_1}.

The step size $\gamma$ can be a constant or computed using methods such as line search to achieve optimal convergence. In this work, the steepest descent and Nesterov accelerated gradient \citep{nesterov1983method} optimization algorithms are used, along with the Barzilai-Borwein (\acrshort{BB}) and Quasi-Newton Barzilai-Borwein (\acrshort{QN-BB}) methods to compute the step size. 
The readers are referred to \cite{antonau2025comparison,antonau2023enhanced} for explanation of the algorithms and step size methods, which have not been included here to avoid repetition.

\FloatBarrier

\subsection{Comparison}
\label{sec:comparison}

This section describes the basis for the qualitative and quantitative comparison of the results.

\subsubsection{Comparison with other scenarios}
\label{sec:scenario_comparison}

For each example, four scenarios with varying thermal-effect considerations are analyzed.
In Scenario 1, the actual structure is subjected to thermal loading, but the system identification does not account for it. This scenario depicts the failure in detecting weaknesses when thermal effects are ignored during \acrshort{SI}. 
Scenario 2 accounts for thermal effects by assuming a constant temperature distribution across the structure. This is a simplistic way to account for thermal effects in the absence of temperature measurements.
In Scenario 3, a sparse network of temperature sensors is used to approximate the structure's thermal field via interpolation. In this work, a standard implementation of \acrfull{kNN} from the Python Scikit-learn library \citep{scikit_learn} is used. 
Scenario 4 concerns the proposed monolithic and partitioned \acrshort{SI} approaches, which incorporate the temperature sensor measurements directly into the \acrshort{SI} cost function.

\begin{table}[!b]
\caption{Different scenarios without and with accounting for thermal effects during \acrshort{SI}.}
\label{tab:Table_1}
\begin{tabular}{>{\raggedright\arraybackslash}p{0.35\textwidth} 
                >{\centering\arraybackslash}p{0.235\textwidth} 
                >{\centering\arraybackslash}p{0.16\textwidth} 
                >{\centering\arraybackslash}p{0.13\textwidth}}
\hline \hline
 \textbf{} & \multicolumn{3}{c}{\textbf{System Identification}} \\ \cline{2-4}
\textbf{Scenario}  & \textbf{Thermal effects accounted for?} & \textbf{If yes, how?} & \textbf{Identified quantity} \\ \hline
\textbf{1)} \acrshort{SI} \textit{without} considering temperature  & No & - & $\mathbf{E}$ \\ \hline
\textbf{2)} \acrshort{SI} considering \textit{constant} temperature field  & Yes & $\boldsymbol{\Delta \mathbf{T}}$ constant & $\mathbf{E}$ \\ \hline
\textbf{3)} \acrshort{SI} considering temperature \textit{interpolation}  & Yes & $\boldsymbol{\Delta \mathbf{T}}$ interpolated & $\mathbf{E}$ \\ \hline
\textbf{4)} \acrshort{SI} considering temperature \textit{identification}  & Yes & $\boldsymbol{\Delta \mathbf{T}}$ identified & $\mathbf{E}$, $\boldsymbol{\Delta \mathbf{T}}$ \\ \hline \hline
\end{tabular}
\end{table}

These scenarios are described in Table \ref{tab:Table_1}. As one moves down the table, the sophistication with which thermal effects are accounted for increases, and the \acrshort{SI} becomes more complex. For Scenario 4, the ill-conditioning of the \acrshort{SI} also increases because two distributions ($\mathbf{E}$, $\boldsymbol{\Delta \mathbf{T}}$) are being identified instead of one distribution ($\mathbf{E}$) in Scenarios 1--3.

\subsubsection{Identification Error Comparison}
\label{sec:comparison_error}

To evaluate the effectiveness of the proposed methodology, a qualitative visual examination of the Young's modulus and temperature distributions is first performed. Secondly, since this is an academic study, and the `target' distributions in the structures are prescribed, the relative discrete \acrfull{$L_2$} norm between the target and the reconstructed quantity can be computed. 
In this discrete setting, it is calculated as 
\begin{equation}
    \text{Relative discrete $L_2$ norm} = \epsilon_{L_2} = \frac{\sqrt{\sum_{i=1}^{N} ( Q_i^{\text{}} - Q_i^{\text{target}})^2 }}{\sqrt{\sum_{i=1}^{N} ( Q_i^{\text{target}})^2 }} \quad ,
    \label{eq:L2_norm}
\end{equation}
where $i=1,\dots,N$ is the loop through the nodes/elements with the total number of nodes/elements $N$; $Q_i^{\text{}}$ and $Q_i^{\text{target}}$ are the reconstructed (interpolated, identified, etc) and the target quantities of interest for the $i$-th node/element respectively. In this work, $Q$ corresponds either to the elemental Young’s modulus distribution or to the nodal temperature field, depending on the quantity being examined.

The percentage change in $\epsilon_{L_2}$ with respect to the start of the optimization is also computed to gauge the improvement (or deterioration) during \acrshort{SI}.
\begin{equation}
    \delta \epsilon_{L_2} = \frac{\epsilon_{L_2}^{} - \epsilon_{L_2}^{0}}{\epsilon_{L_2}^{0}} \cdot 100 \quad,
    \label{eq:L2_pc_change}
\end{equation}
where $\epsilon_{L_2}^{}$ refers to the relative discrete \acrshort{$L_2$} error for the case being examined and $\epsilon_{L_2}^{0}$ refers to the reference error, here, the relative discrete \acrshort{$L_2$} error at the start of the optimization. 

It is pointed out that the $\epsilon_{L_2}$ and $\delta \epsilon_{L_2}$ error metrics are different from the cost function sensor errors (Eqs.~\eqref{eq:cost_function}). In $\epsilon_{L_2}$ and $\delta \epsilon_{L_2}$, the entire distribution i.e., each node/element, is compared to the `target' distribution, whereas, in the cost function, the quantities are only known and compared at a few discrete sensor locations. In practice, since the actual `target' is unknown, only the cost function is available, and hence, the reason for system identification.


\section{Numerical Examples}
\label{s:numerical}

This section presents two numerical examples to demonstrate the applicability of the proposed combined Young's modulus and thermal field identification methodology. The first example is a Plate With a Hole modeled with shell elements and analyzed for different numbers of temperature sensors and configurations.  The second example of a Footbridge modeled with beam and shell elements illustrates the applicability of the approach for models with mixed finite element types. 
Both the examples are tested on a linearly varying and a localized (Gaussian-type) thermal load cases. 

The \acrshort{SI} and \acrshort{FE} analyses are performed using KratosMultiphysics \citep{dadvand2010object,vicente_mataix_ferrandiz_2025_15687676}, an open source multi-physics 
software supporting multiple element types and formulations, constitutive models, as well as optimization capabilities.

The numerical examples are setup as follows: To simulate the damaged structure subjected to thermal loads, a `target' Young's modulus distribution and a `target' temperature field is prescribed and the structure is simulated with the corresponding external mechanical load $f_{\text{ext}}$. The displacements and temperatures are recorded at the corresponding sensor locations (displacement at $\mathbf{x}_j$, $j=1,\dots,m$ and temperatures at $\mathbf{x}_p$, $p=1,\dots,o$). In this work, a single thermal field distribution ($\Delta \text{T}(\mathbf{x})$) shall be identified, implying that all mechanical load test cases shall be subjected to the same thermal field. This yields the measurement pairs $f_{\text{ext},i},i=1,\dots,n,\mathbf{u}_{ij}^{\text{meas}}, j=1,\dots,m, \Delta \text{T}_{p}^{\text{meas}}, p=1,\dots,o$ which is used as the `target measured data' during \acrshort{SI} to identify the damage i.e., Young's modulus distribution ($\mathbf{E}$) and the temperature distribution ($\boldsymbol{\Delta \mathbf{T}}$) in the structure. The `target measured data' is supposed to be recorded on the real structure whose state is being identified and hence also referred to as `actual structure'. 

For each example, the four scenarios described in Section \ref{sec:scenario_comparison} are analyzed, ranging from when thermal effects are not accounted for to when they are accounted for using different approaches. 

\subsection{Plate With Hole}
\label{s:plate_with_hole}

\begin{figure}[!t]
\centering
\subfloat[\centering]{\includegraphics[trim=0 450 0 450, clip, width=0.485\linewidth]{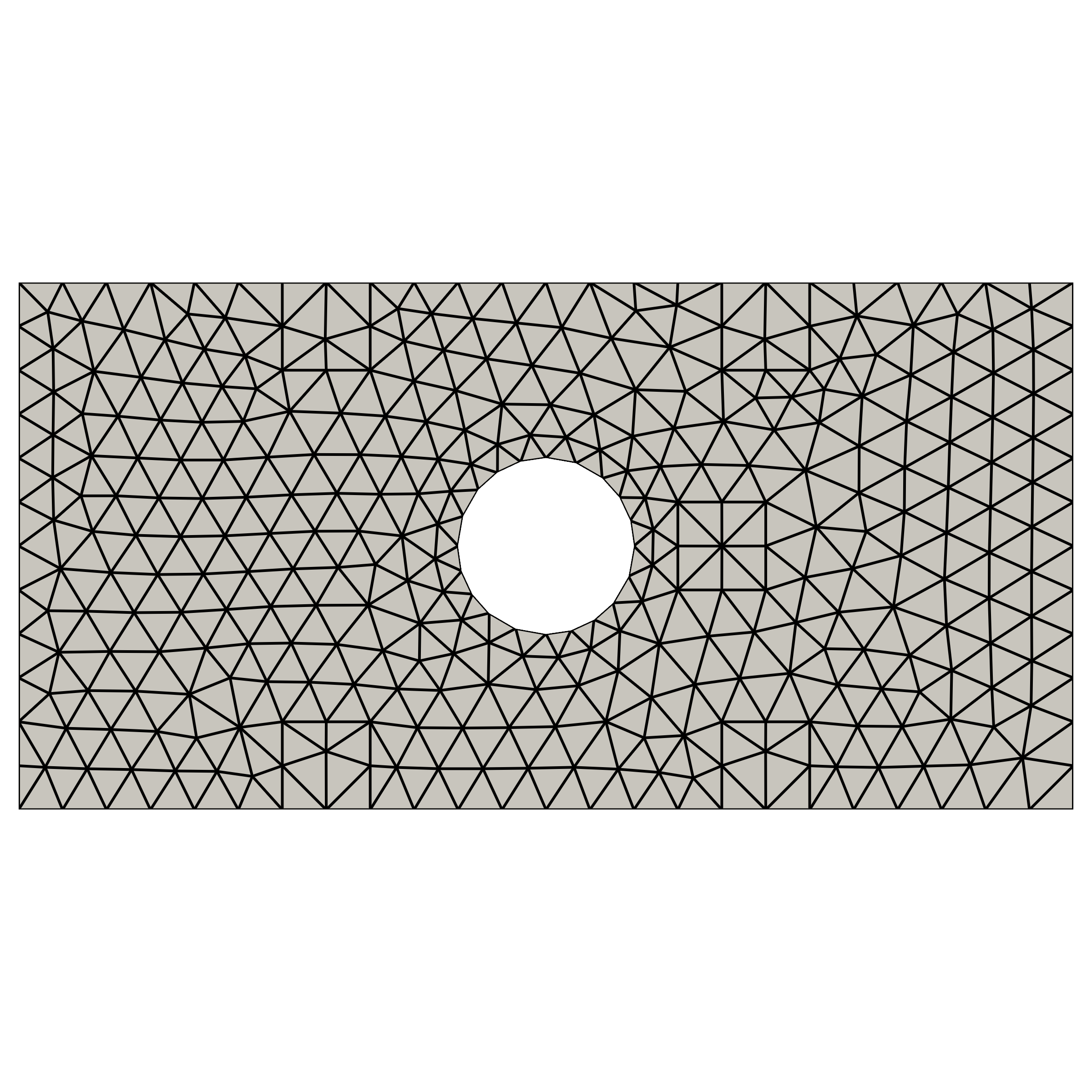}
\label{f:plate_mesh_a}}
\hfill
\subfloat[\centering]{\includegraphics[trim=0 470 0 470,, clip, width=0.485\linewidth]{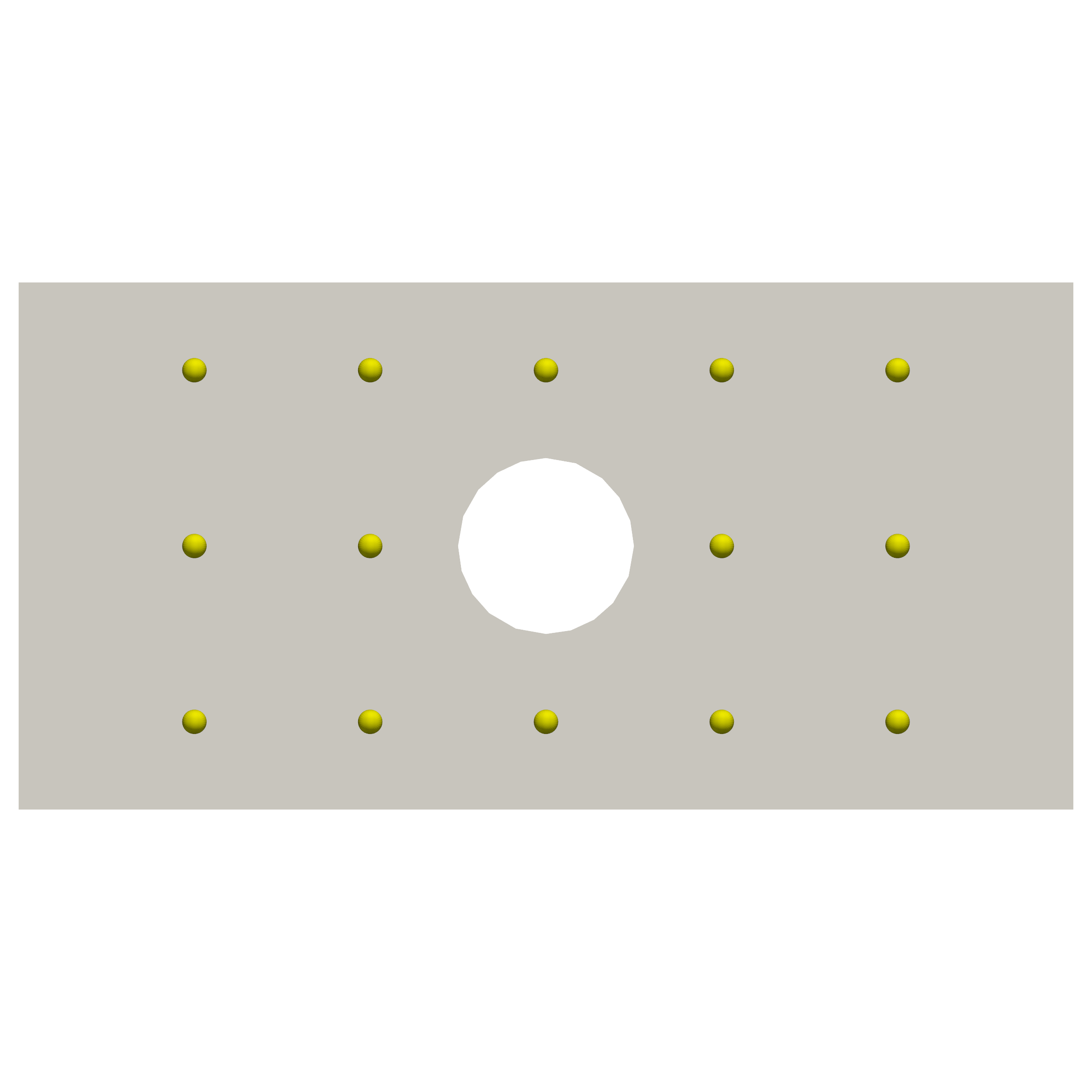}
\label{f:plate_mesh_b}}
\caption{\textbf{Plate with Hole.} \textbf{(a)} Mesh and \textbf{(b)} location of the 14 displacement sensors.}
\label{f:plate_mesh}
\end{figure} 

\begin{figure}[!b]
\centering
\includegraphics[trim=0 370 0 420, clip, width=0.5\linewidth]{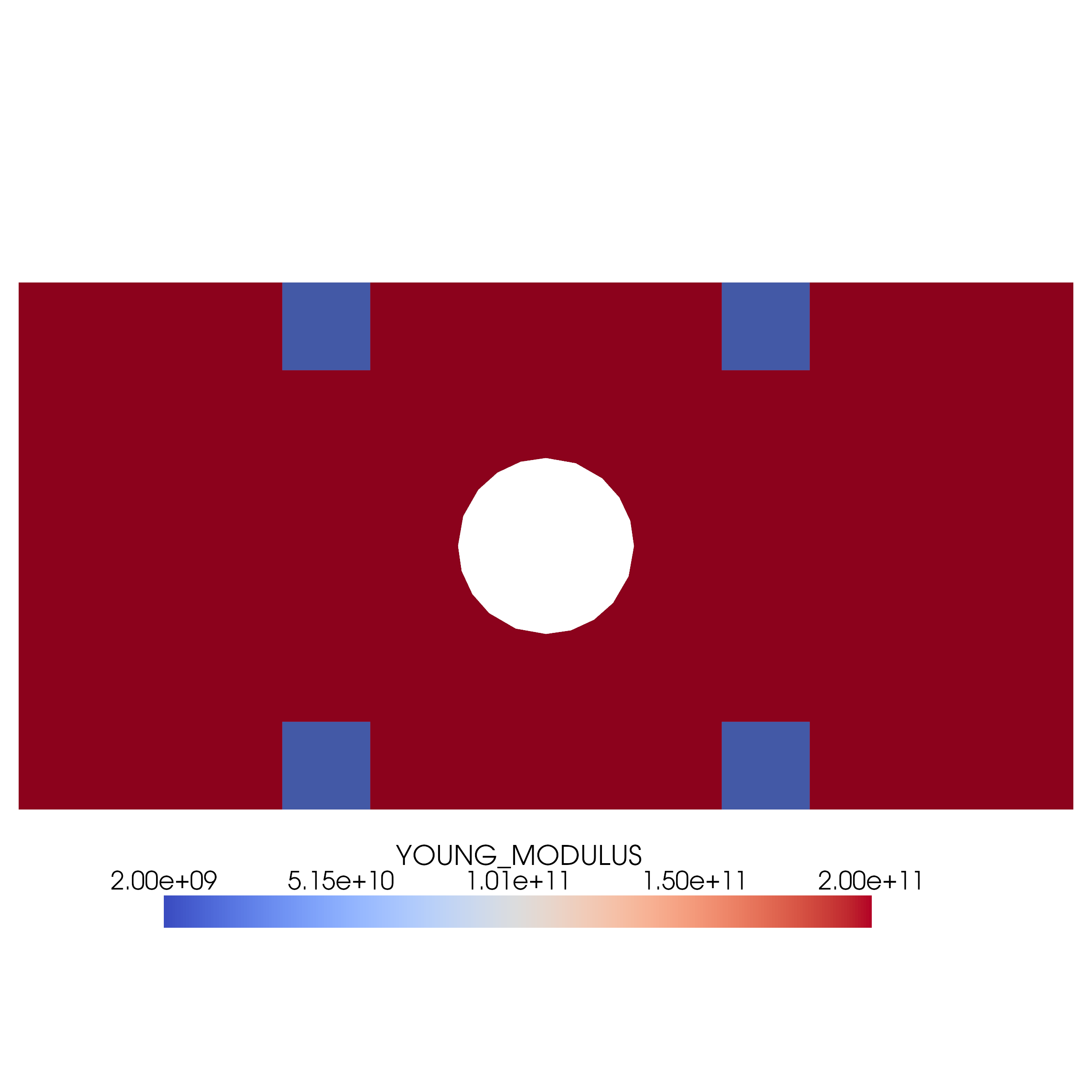}
\caption{\textbf{Plate with Hole.} Target Young's modulus distribution, i.e., localized damages.}
\label{f:plate_target_damage}
\end{figure}

The plate with a hole structure is shown in Fig. \ref{f:plate_mesh}. The plate has dimensions: $0 \le x \le 60$, $0 \le y \le 30$, and $0 \le z \le 0.1$. All quantities are in SI units, unless explicitly mentioned. A diameter $d=10$ hole is located in the plate with center at $x = 30$, $y = 15$. The material properties of the plate are: density $\rho = 7800$, Young's modulus (in pristine state) $\text{E} = 2\cdot 10^{11}$, Poisson's ratio $\nu = 0.3$, and coefficient of thermal expansion $\alpha = 1.0 \cdot 10^{-5}$ /\degree C. The plate is discretized using $646$ linear, triangular, plane stress shell elements as illustrated in Fig. \ref{f:plate_mesh_a}. The plate is fixed on the left edge ($\mathbf{u}(x=0)= \mathbf{0}$), and a static horizontal mechanical line load of $q=1\cdot10^{7}$ N/m is applied on the right edge ($x=60$).

\begin{figure}[!t]
\centering
\subfloat[\centering]{\includegraphics[trim=0 300 0 420, clip, width=0.49\linewidth]{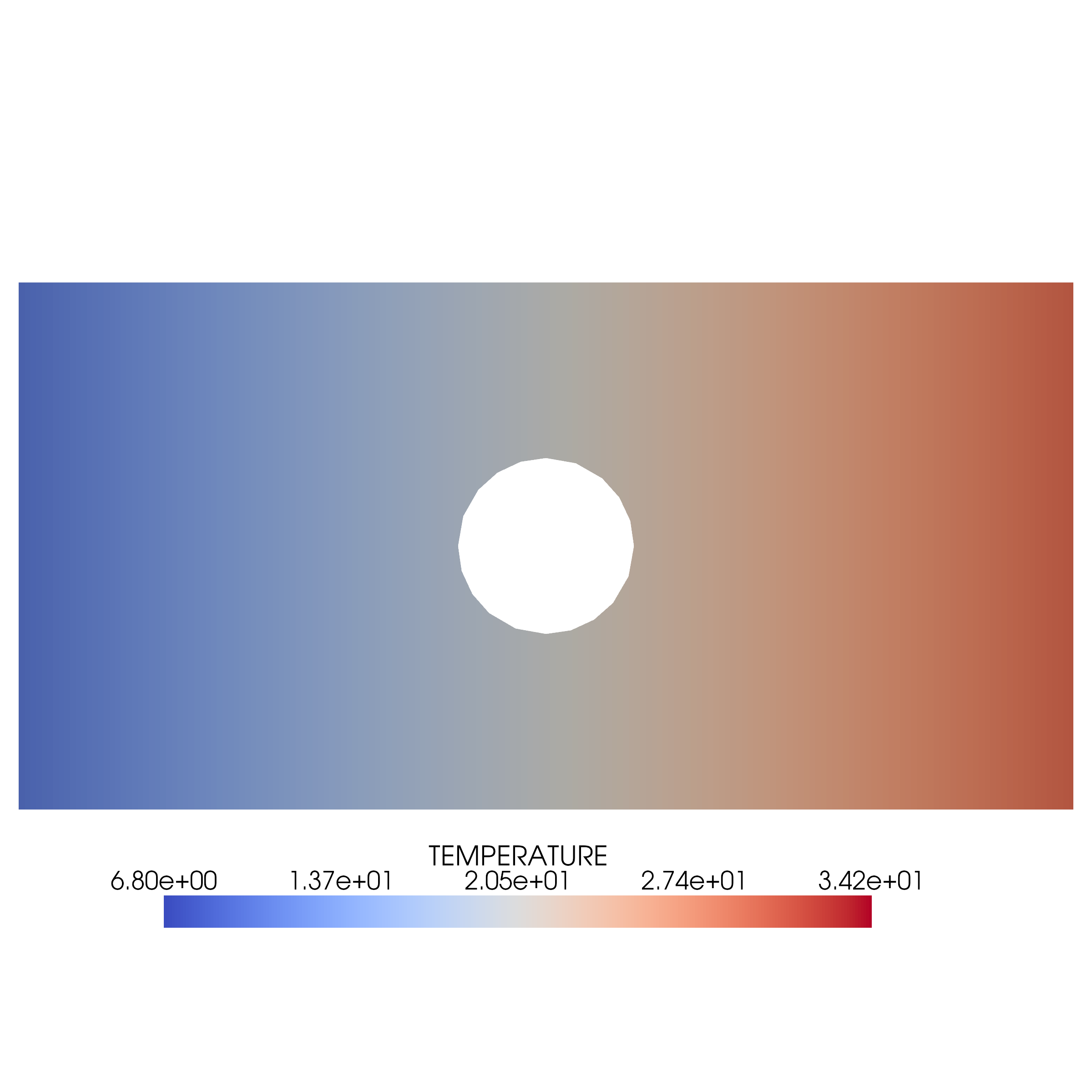}}
\hfill
\subfloat[\centering]{\includegraphics[trim=0 300 0 420, clip, width=0.49\linewidth]{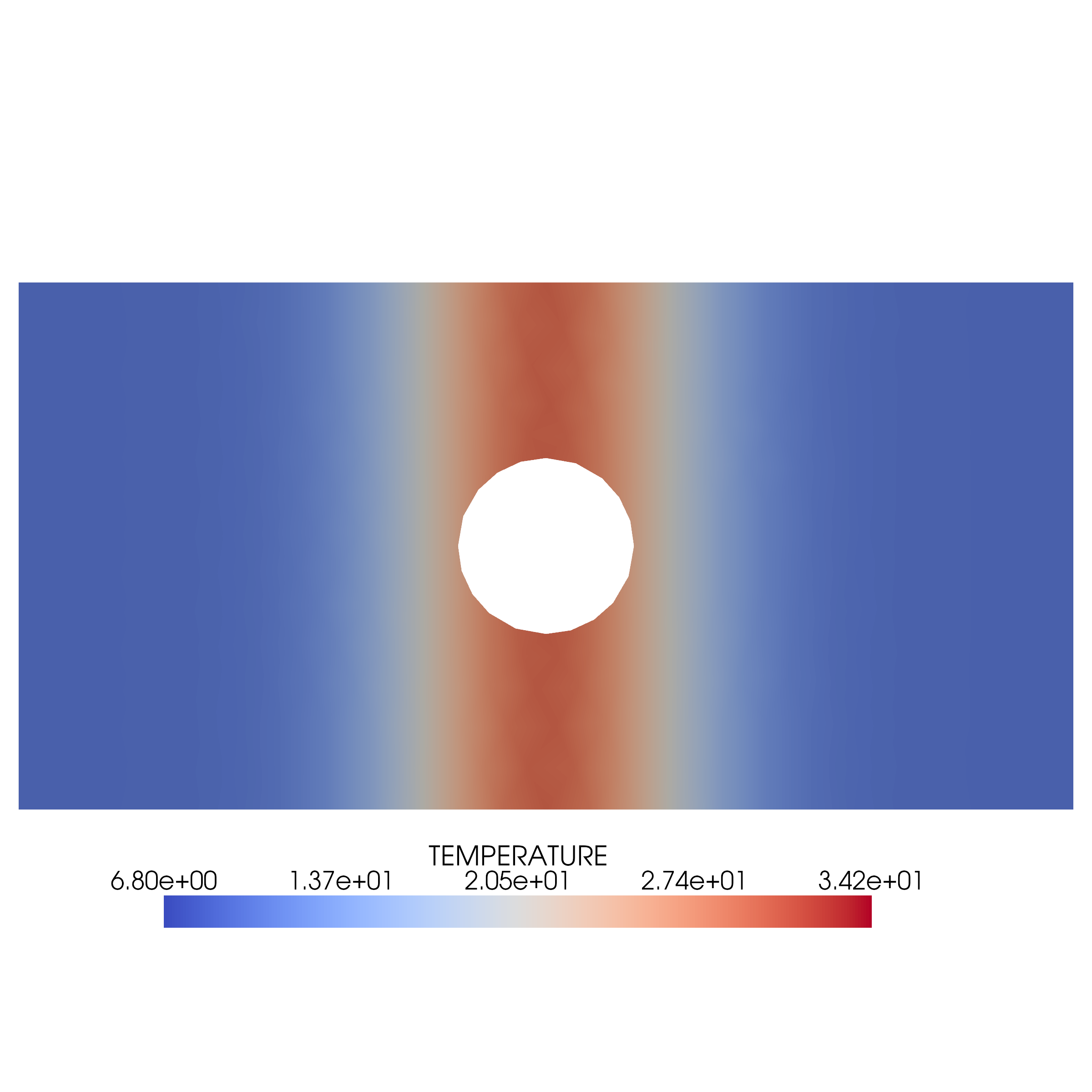}}
\caption{\textbf{Plate with Hole.} Target temperature distributions: (\textbf{a}) Linearly varying thermal field, and (\textbf{b}) Localized thermal field.}
\label{f:plate_target_temp}
\end{figure}

\begin{figure}[!b]
    \centering
    \begin{subfigure}[t]{\textwidth}
        \centering
        \begin{minipage}[t]{0.49\textwidth}
            \centering
            \includegraphics[trim=0 300 0 420, clip, width=\textwidth]{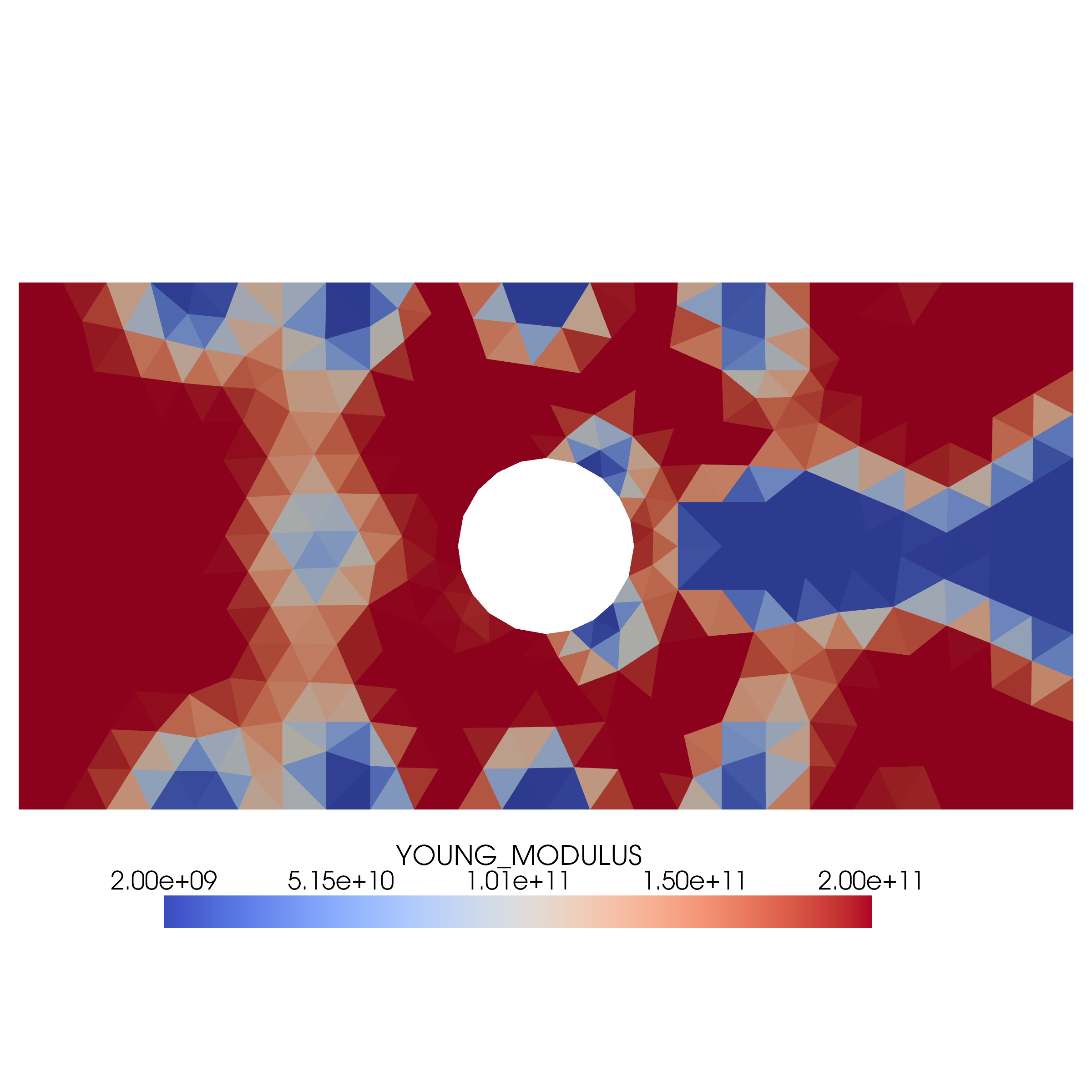}
        \end{minipage}
        \hfill
        \begin{minipage}[t]{0.49\textwidth}
            \centering
            \includegraphics[trim=0 0 0 0, clip, width=\textwidth]{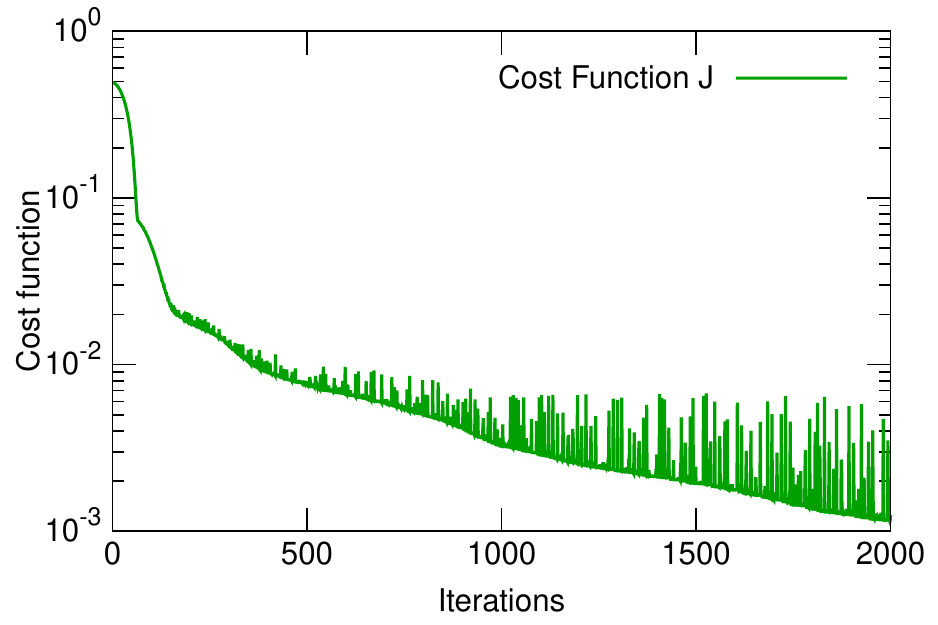}
        \end{minipage}
        \caption{}
        \label{f:plate_wo_temp_a}
    \end{subfigure}
    \begin{subfigure}[t]{\textwidth}
        \centering
        \begin{minipage}[t]{0.49\textwidth}
            \centering
            \includegraphics[trim=0 300 0 420, clip,width=\textwidth]{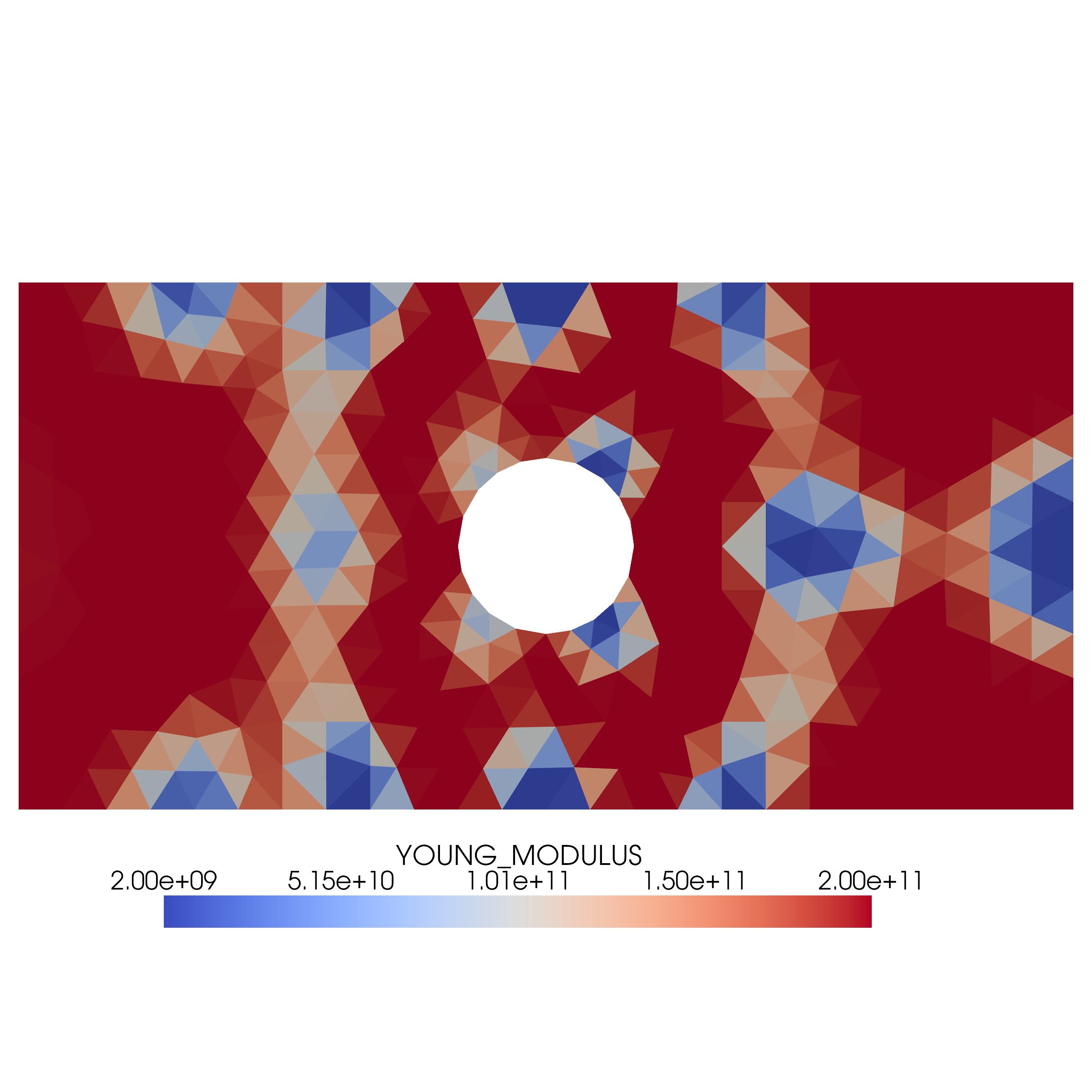}
        \end{minipage}
        \begin{minipage}[t]{0.49\textwidth}
            \centering
            \includegraphics[trim=0 0 0 0, clip,width=\textwidth]{Figure_5b.pdf}
        \end{minipage}  
        \caption{}
        \label{f:plate_wo_temp_b}
    \end{subfigure}
    \caption{\textbf{Plate with Hole.} Identified Young's modulus distributions when thermal load is \textit{not considered} during \acrshort{SI}, but the actual structure is subjected to a: (\textbf{a}) Linearly varying thermal field, (\textbf{b}) Localized thermal field.}
\label{f:plate_wo_temp}
\end{figure}

In the damaged state, the Young's modulus was prescribed to be reduced to $\text{E}=2\cdot10^{10}$ in the four locally damaged regions: upper left ($15 \leq x \leq 20$, $25 \leq y \leq 30$), upper right ($40 \leq x \leq 45$, $25 \leq y \leq 30$), lower left ($15 \leq x \leq 20$, $0 \leq y \leq 5$), and lower right ($40 \leq x \leq 45$, $0 \leq y \leq 5$). This `target' Young's modulus distribution is shown in Fig. \ref{f:plate_target_damage}.
The actual structure is also subjected to a thermal load. Here, two different types of thermal fields were investigated. The first was a linearly varying thermal field with $\Delta \text{T}(x=0)=30 $\degree C at the left edge and $\Delta \text{T}(x=60)=10 $\degree C at the right edge. The second was a localized Gaussian type thermal field defined as $\Delta \text{T}(\mathbf{x}) = \Delta \text{T}_{\text{min}} + (\Delta \text{T}_{\text{max}}- \Delta \text{T}_{\text{min}})e^{\frac{-(x-30)^{2}}{75}}$, $\Delta \text{T}_{\text{min}}=10 $\degree C, $\Delta \text{T}_{\text{max}}=30 $\degree C. These `target' temperature distributions are illustrated in Fig. \ref{f:plate_target_temp}. The color bar is shown in the range $[6.8, 34.2]$ \degree C for consistent visualization and better comparison with later results.

To measure the deformations, fourteen displacement sensors were distributed over the plate, as depicted in Fig. \ref{f:plate_mesh_b}. At the start of \acrshort{SI}, the plate was assumed to be in near-pristine condition with $\text{E}(\mathbf{x}) = 1.998\cdot 10^{11}$. 
The `maximum measured value' sensor weighting (Eq.\eqref{eq:sensor_weight_a}) was used for normalization. 
The steepest descent algorithm with \acrshort{BB} method and a maximum step size of $1\cdot 10^{-1}$ was used for the optimization.
The convergence criteria were set to a 6-magnitude reduction of the cost function, i.e., $J \leq (1\cdot 10^{-6} * J_0)$ or a maximum of 2000 iterations. The Young's modulus control variables were bounded in the range of $[2\cdot10^{9},2\cdot10^{11}]$.
To help with the ill-conditioned \acrshort{SI}, \acrshort{VM} with a linear kernel and radius $r = 5$ was employed to regularize gradients and updates. Based on the rule of thumb, the radius contained approximately $2-3$ element side lengths. 

To compare results, the damage locations and intensities, as well as the optimization and \acrshort{VM} settings, were kept the same across all cases in this example, except in Scenario 4, where the overall maximum number of iterations differed (discussed later).

\begin{figure}[!b]
    \centering
    \begin{subfigure}[t]{\textwidth}
        \centering
        \begin{minipage}[t]{0.49\textwidth}
            \centering
            \includegraphics[trim=0 300 0 420, clip, width=\textwidth]{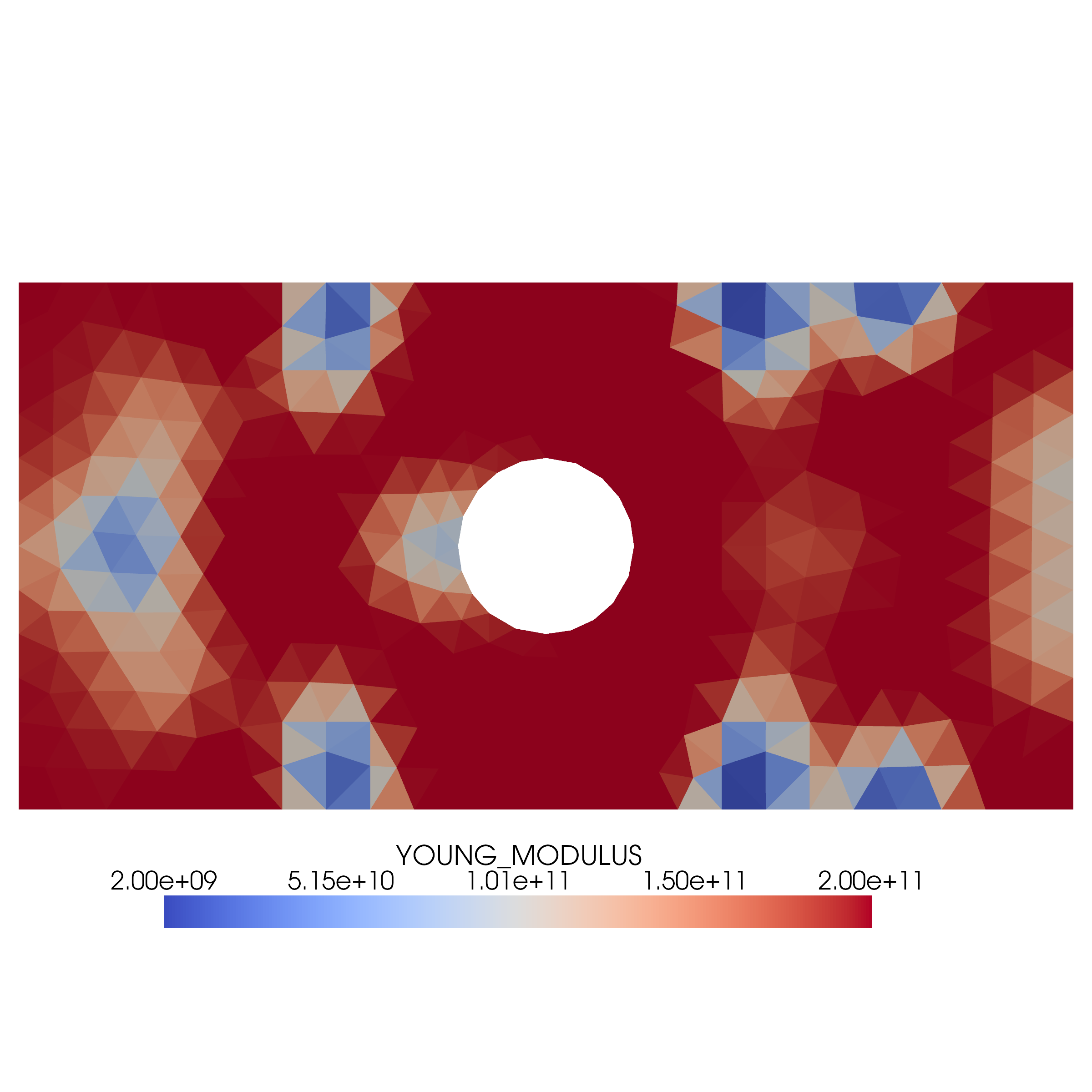}
        \end{minipage}
        \hfill
        \begin{minipage}[t]{0.49\textwidth}
            \centering
            \includegraphics[trim=0 0 0 0, clip, width=\textwidth]{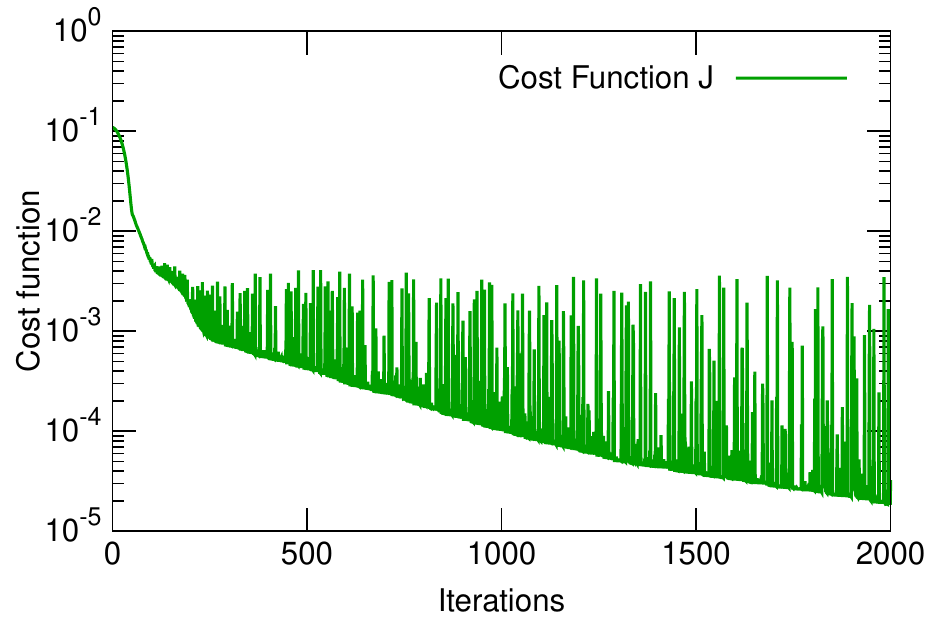}
        \end{minipage}
        \caption{}
        \label{f:plate_w_const_temp_a}
    \end{subfigure}
    \begin{subfigure}[t]{\textwidth}
        \centering
        \begin{minipage}[t]{0.49\textwidth}
            \centering
            \includegraphics[trim=0 300 0 420, clip,width=\textwidth]{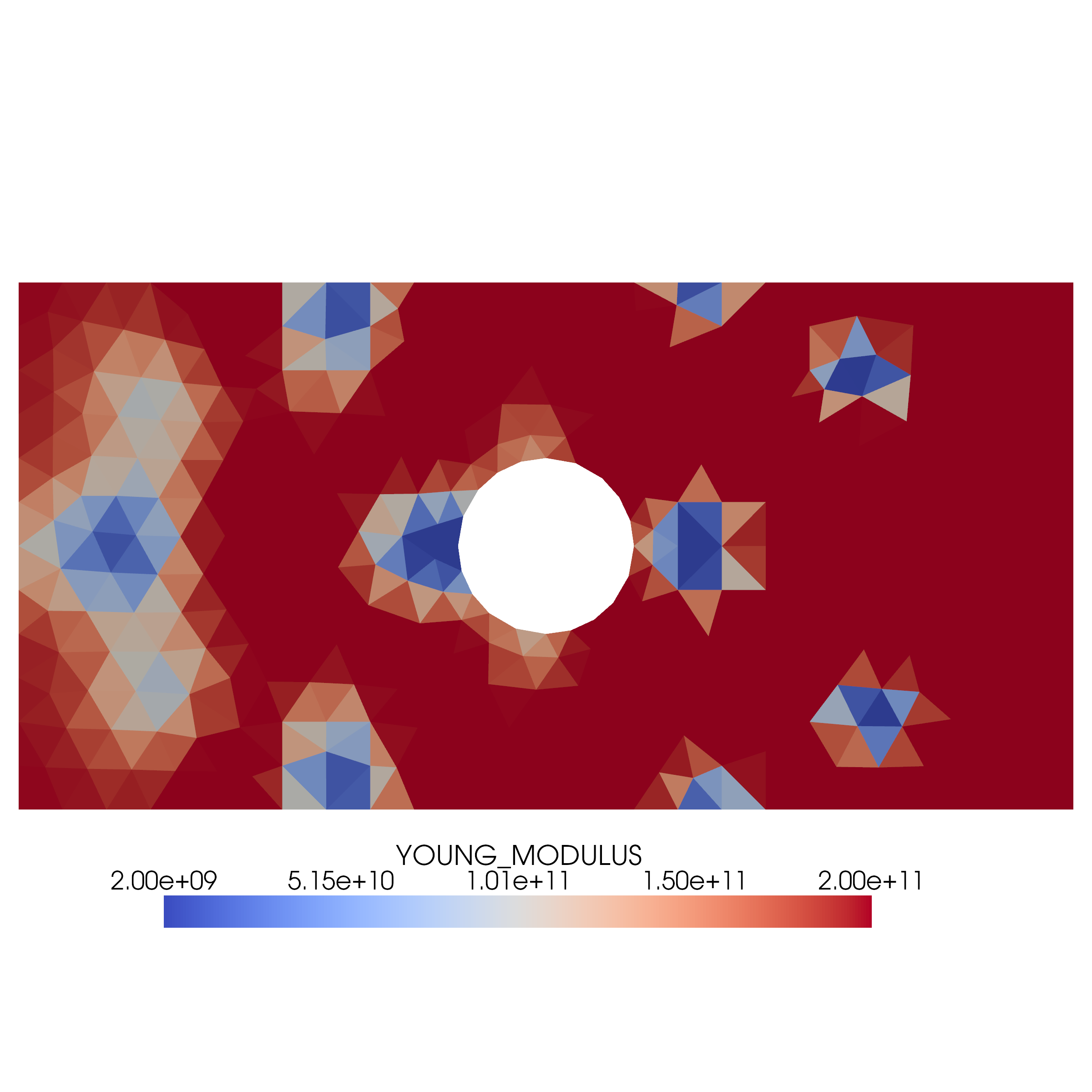}
        \end{minipage}
        \begin{minipage}[t]{0.49\textwidth}
            \centering
            \includegraphics[trim=0 0 0 0, clip,width=\textwidth]{Figure_6b.pdf}
        \end{minipage}  
        \caption{}
        \label{f:plate_w_const_temp_b}
    \end{subfigure}
    \caption{\textbf{Plate with Hole.} Identified Young's modulus distributions when a \textit{constant} temperature distribution of $20 $\degree C is considered during \acrshort{SI}, but the actual structure is subjected to a: (\textbf{a}) Linearly varying thermal field, (\textbf{b}) Localized thermal field.}
\label{f:plate_w_const_temp}
\end{figure}

Scenario 1, as described in Table \ref{tab:Table_1}, is the case where the actual structure is subjected to thermal loads but thermal effects were not considered during \acrshort{SI}. 
The identified Young's modulus distribution for this scenario for the two thermal fields and their convergence plots are shown in Fig. \ref{f:plate_wo_temp}. It can be observed that, in both thermal field cases, damage detection and localization are severely deficient. A lot of false damage is identified, which can be directly attributed to the lack of consideration of the structure's thermal expansion. The convergence plots show that the cost functions decrease by several magnitudes (by compensating for thermal expansion with damage) and are of order $\mathcal{O}(10^{-3})$. 

The spiking behavior of the convergence curve is an inherent behavior observed with the \acrshort{BB} method \citep{fletcher2005barzilai,ansari2025adjoint,antonau2023enhanced} and is inconsequential for this study, which does not focus on the optimization algorithm and step size method. 

Scenario 2, as described in Table \ref{tab:Table_1}, improves on Scenario 1, and accounts for thermal effects via a constant $20$\degree C thermal field applied to the structure in the primal simulation during \acrshort{SI}. 
The identified Young's modulus distribution for this scenario for the two thermal fields and their convergence plots are shown in Fig. \ref{f:plate_w_const_temp}. A slight improvement in damage localization is evident in the results. However, overall, the results are still quite poor, with many false damage detections.
Here also, the convergence plots show reductions in the cost function of a few orders of magnitude, but higher reductions than in Scenario 1, and are of order $\mathcal{O}(10^{-4}-10^{-5})$.

\begin{figure}[!t]
\centering
\subfloat[\centering]{\includegraphics[trim=0 450 0 450, clip, width=0.49\linewidth]{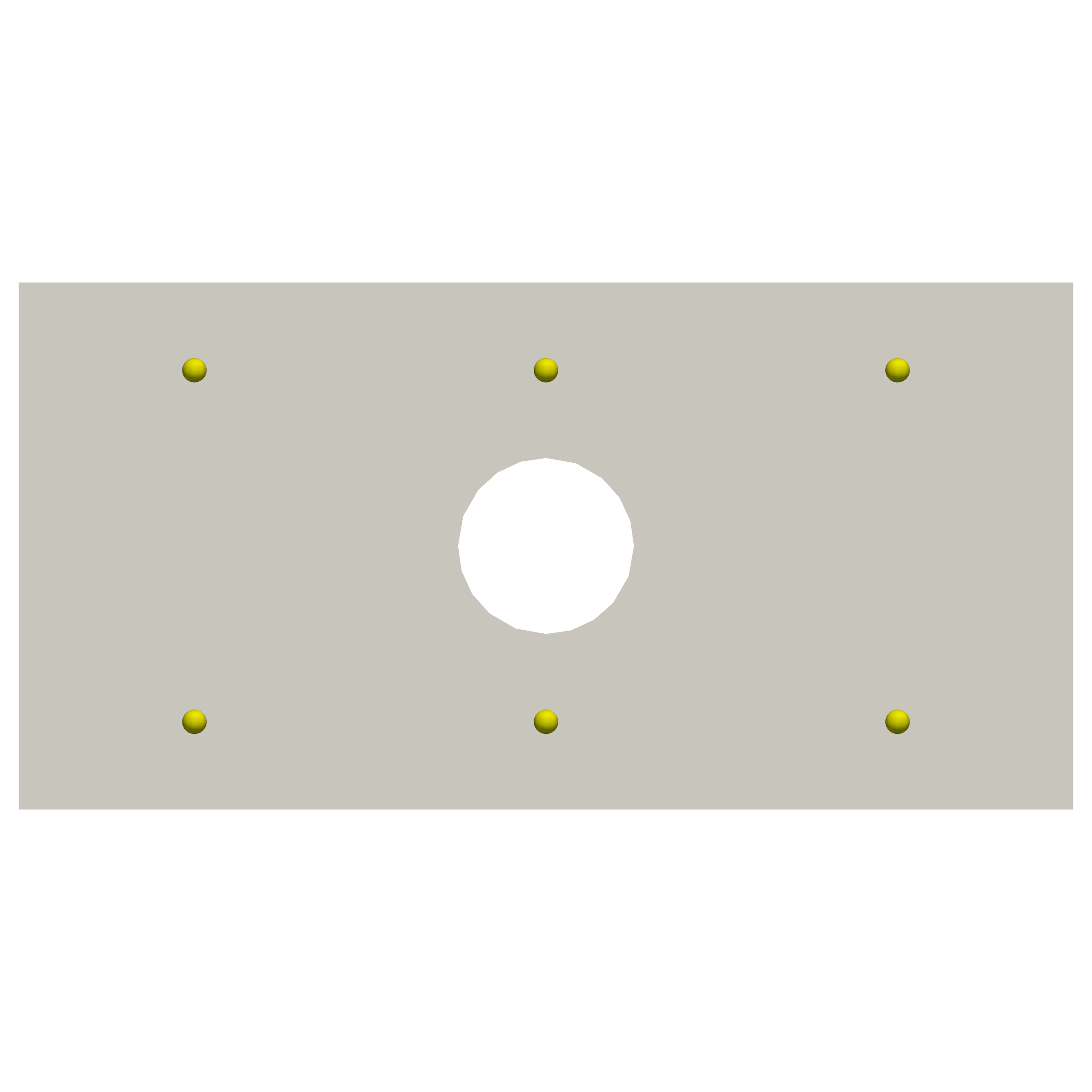}}
\hfill
\subfloat[\centering]{\includegraphics[trim=0 450 0 450, clip, width=0.49\linewidth]{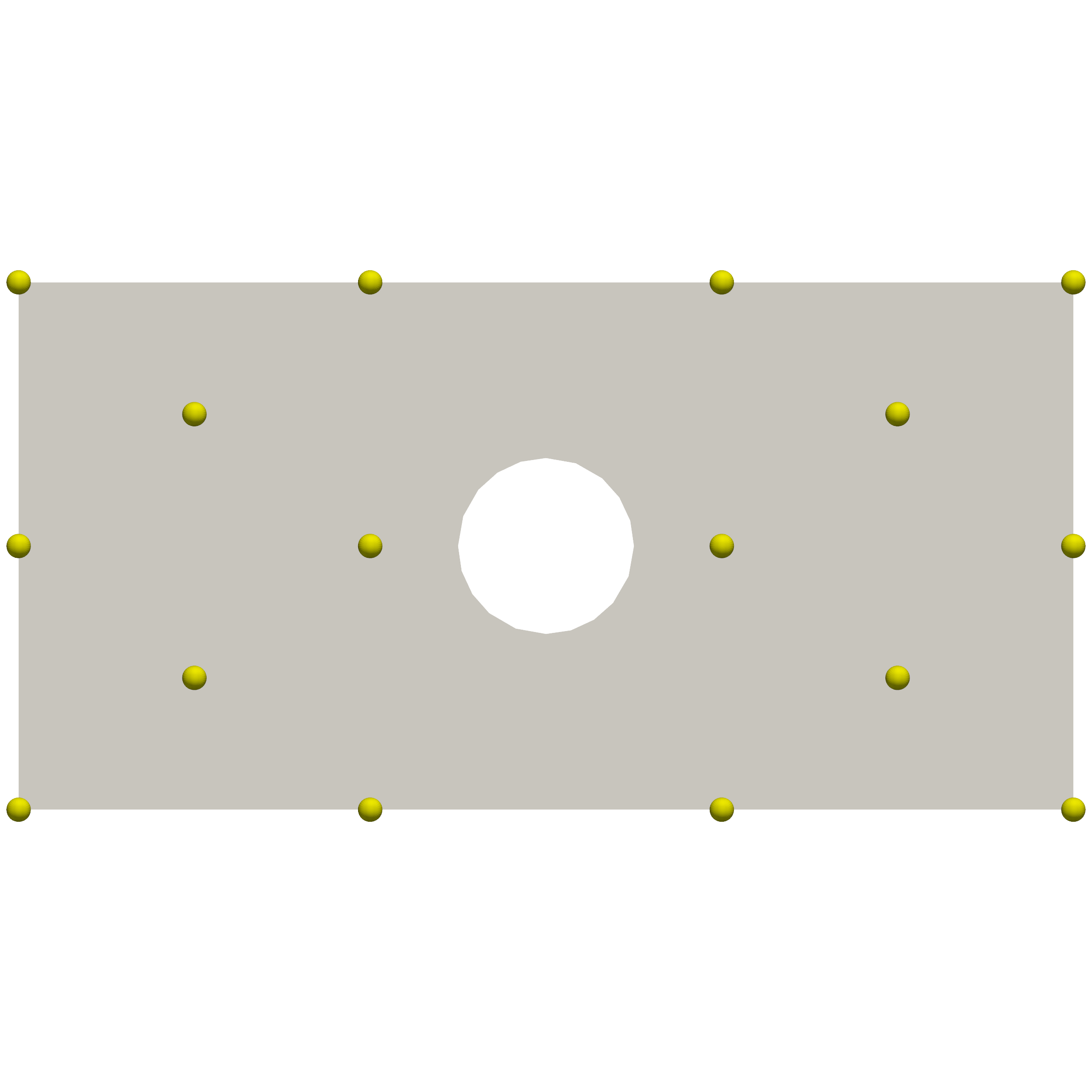}}
\caption{\textbf{Plate with Hole.} Location of the temperature sensors: (\textbf{a}) 6 sensors configuration, (\textbf{b}) 16 sensors configuration.}
\label{f:plate_temp_sensors}
\end{figure}

Analyzing the results from these two scenarios, it is evident that a more accurate estimate of the structure's thermal field is essential for meaningful damage localization. Thus, temperature sensors were placed on the structure. In this example, two sensor configurations consisting of 6 and 16 temperature sensors distributed over the plate were analyzed. These configurations are illustrated in Fig. \ref{f:plate_temp_sensors}. 

Scenario 3, as described in Table \ref{tab:Table_1}, is where the structure's temperature distribution is approximated by interpolating the temperature sensor measurements. In this example, interpolation was performed using \acrshort{kNN} with 3 nearest neighbors and inverse distance weighting. This interpolated temperature distribution was applied as a fixed thermal load on the structure in the primal simulation during \acrshort{SI}. The quantity being identified was still the Young's modulus distribution, as the \acrshort{SI} problem remained unchanged; only the manner in which thermal effects are accounted for had changed.

Scenario 4, as described in Table \ref{tab:Table_1}, is where the temperature distribution of the structure is inferred during \acrshort{SI} alongside the Young's modulus distribution. The proposed methodology presents two approaches to achieve this: the monolithic and the partitioned approaches as described in Section \ref{s:methodology}. 

For this example, the monolithic case was set up as follows: the composite cost function is the sum of the accumulated displacement sensor error and the accumulated temperature sensor error, i.e., $J = J_D + J_T$. The \acrshort{SI} problem is reformulated to identify two fields: the Young's modulus distribution ($\mathbf{E}$) and the temperature distribution ($\boldsymbol{\Delta}\mathbf{T}$). At the start of \acrshort{SI}, the plate had a uniform temperature of $\Delta \text{T}(\mathbf{x}) = 20$ \degree C. The bound for the nodal temperature control variable was set to the range $[-10, 40]$ \degree C. 
The `maximum measured value' sensor normalization was applied separately to each sensor type according to Eqs.~\eqref{eq:sensor_weight_a},~\eqref{eq:sensor_weight_b}.
Similar to Young's modulus distribution, \acrshort{VM} with a linear kernel and radius $r = 20$ was used for the temperature field regularization. The convergence criteria were set to a 6-magnitude reduction in the composite cost function, i.e., $J = (J_D + J_T) \leq (1\cdot 10^{-6} * J_0)$ or $2000$ iterations.  All other settings were kept the same.

\begin{figure}[!b]
    \centering
    \begin{subfigure}[t]{\textwidth}
        \centering
        \begin{minipage}[t]{0.49\textwidth}
            \centering
            \includegraphics[trim=0 300 0 420, clip, width=\textwidth]{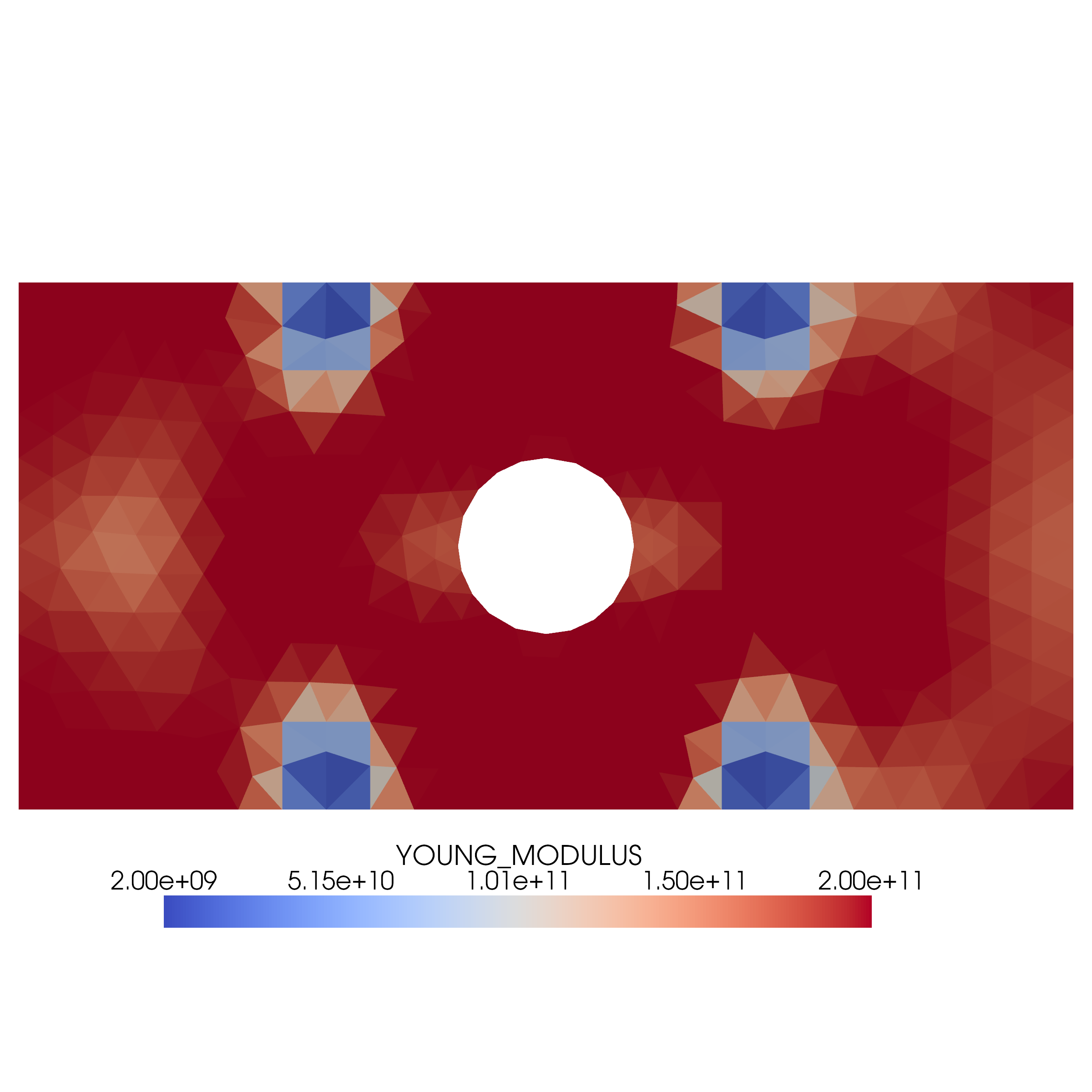}
        \end{minipage}
        \begin{minipage}[t]{0.49\textwidth}
            \centering
            \includegraphics[trim=0 300 0 420, clip, width=\textwidth]{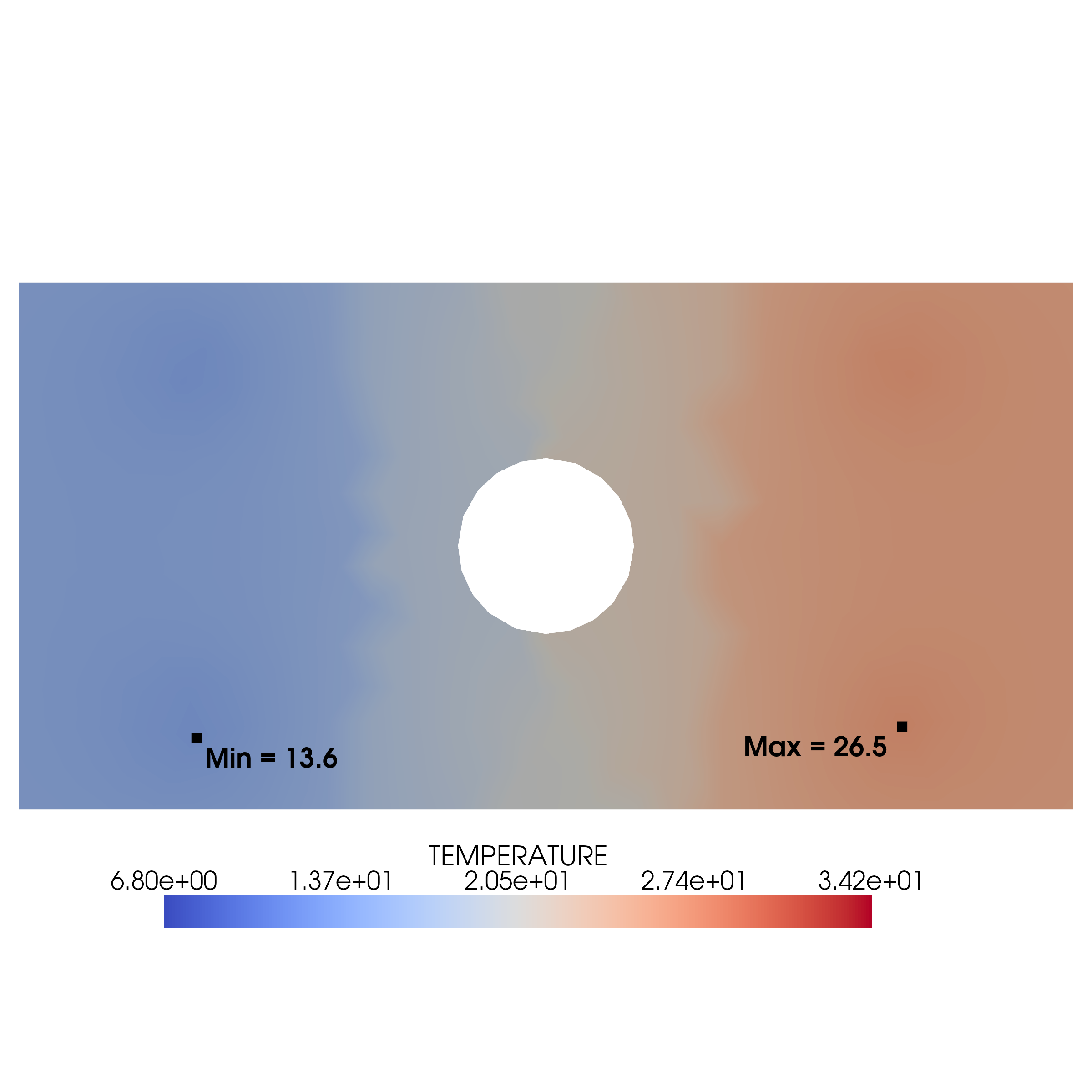}
        \end{minipage}
        \caption{Identified Young's moduli (left) and the temperature distribution (right) when the thermal field is \textit{interpolated}.}
        \label{f:plate_linear_6s_a}
    \end{subfigure}
    \begin{subfigure}[t]{\textwidth}
        \centering
        \begin{minipage}[t]{0.49\textwidth}
            \centering
            \includegraphics[trim=0 300 0 420, clip,width=\textwidth]{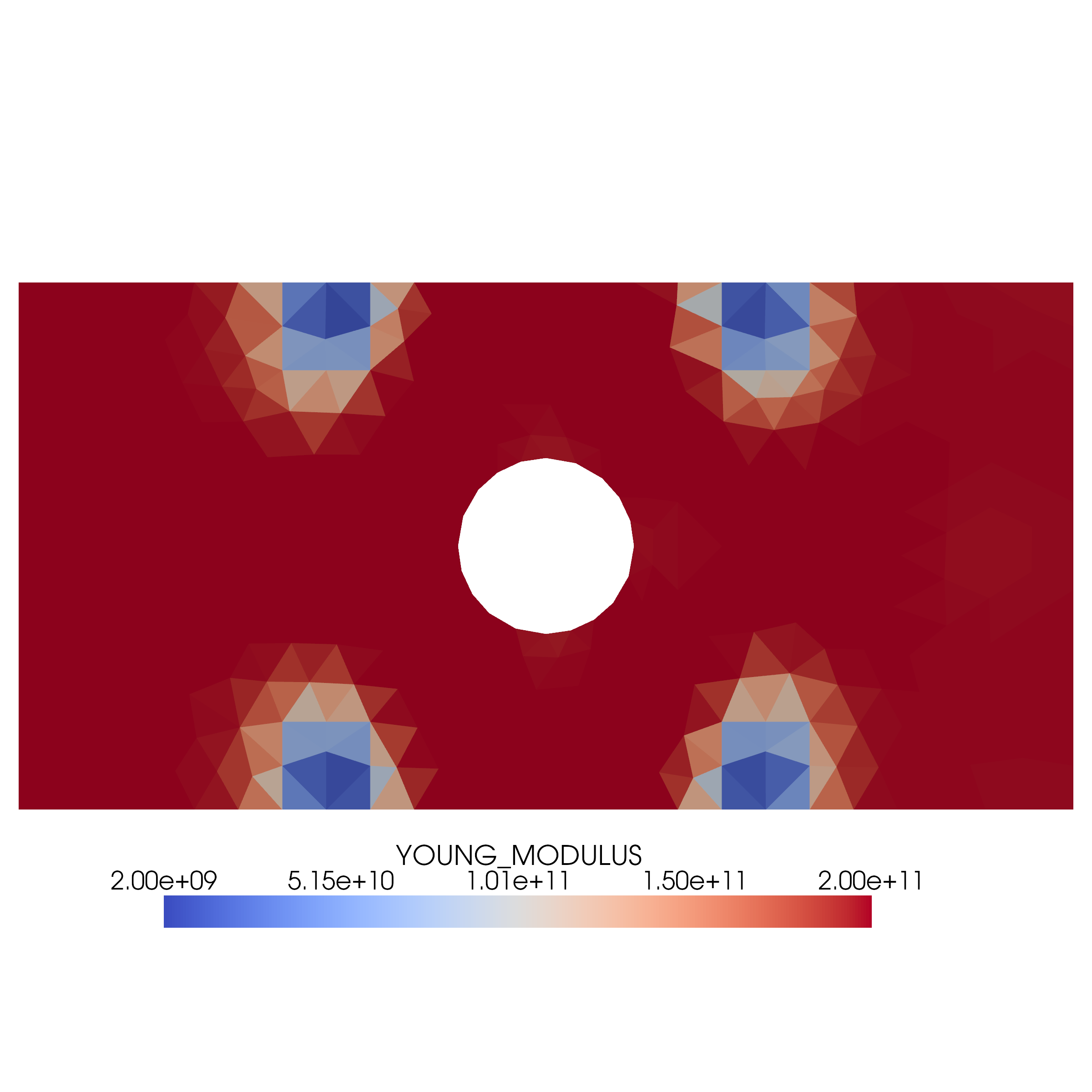}
        \end{minipage}
        \begin{minipage}[t]{0.49\textwidth}
            \centering
            \includegraphics[trim=0 300 0 420, clip,width=\textwidth]{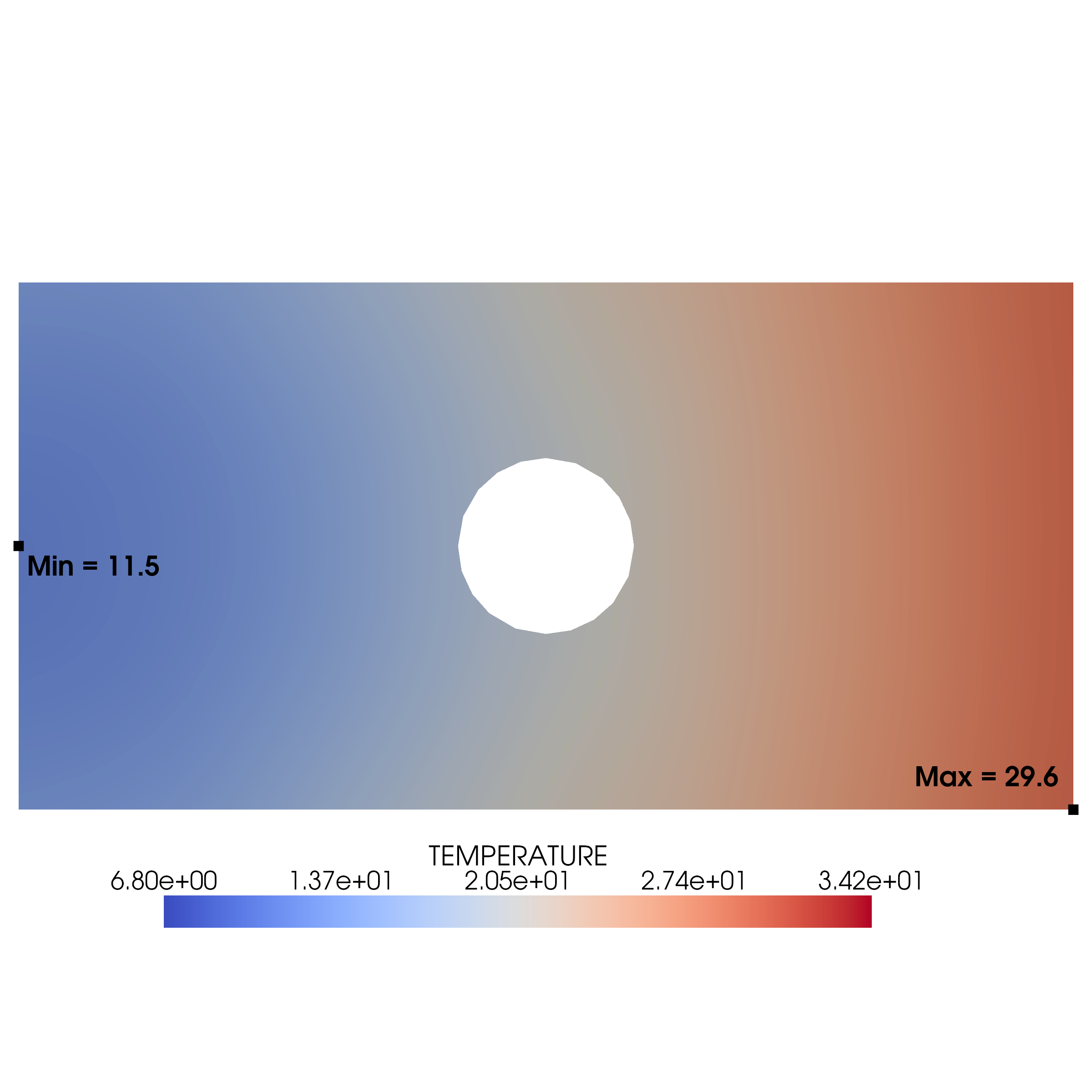}
        \end{minipage}  
        \caption{Identified Young's moduli (left) and \textit{identified} temperature distribution (right) when the thermal field is reconstructed during \acrshort{SI}: \textit{Monolithic} approach.}
        \label{f:plate_linear_6s_b}
    \end{subfigure}
    \begin{subfigure}[t]{\textwidth}
        \centering
        \begin{minipage}[t]{0.49\textwidth}
            \centering
            \includegraphics[trim=0 300 0 420, clip,width=\textwidth]{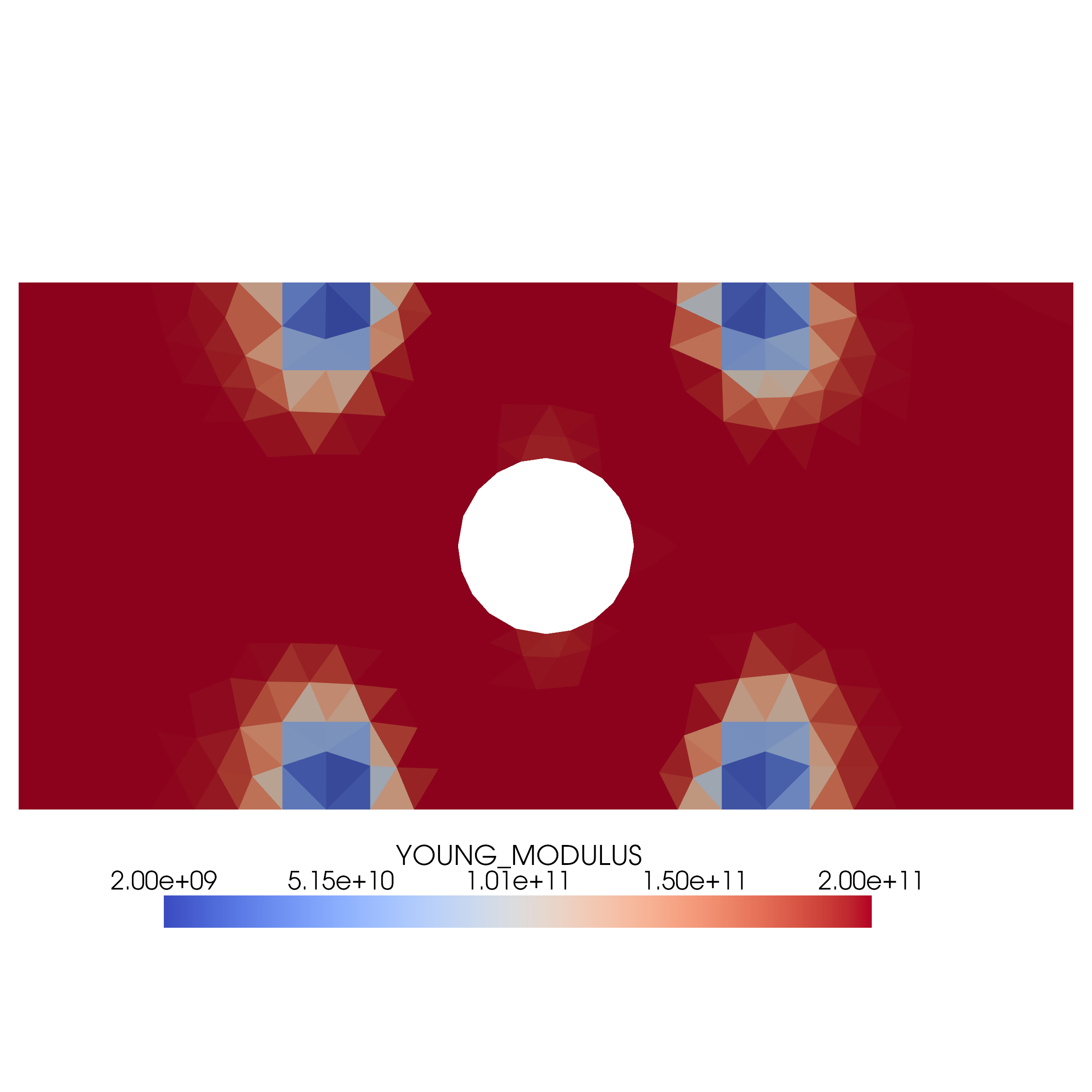}
        \end{minipage}
        \begin{minipage}[t]{0.49\textwidth}
            \centering
            \includegraphics[trim=0 300 0 420, clip,width=\textwidth]{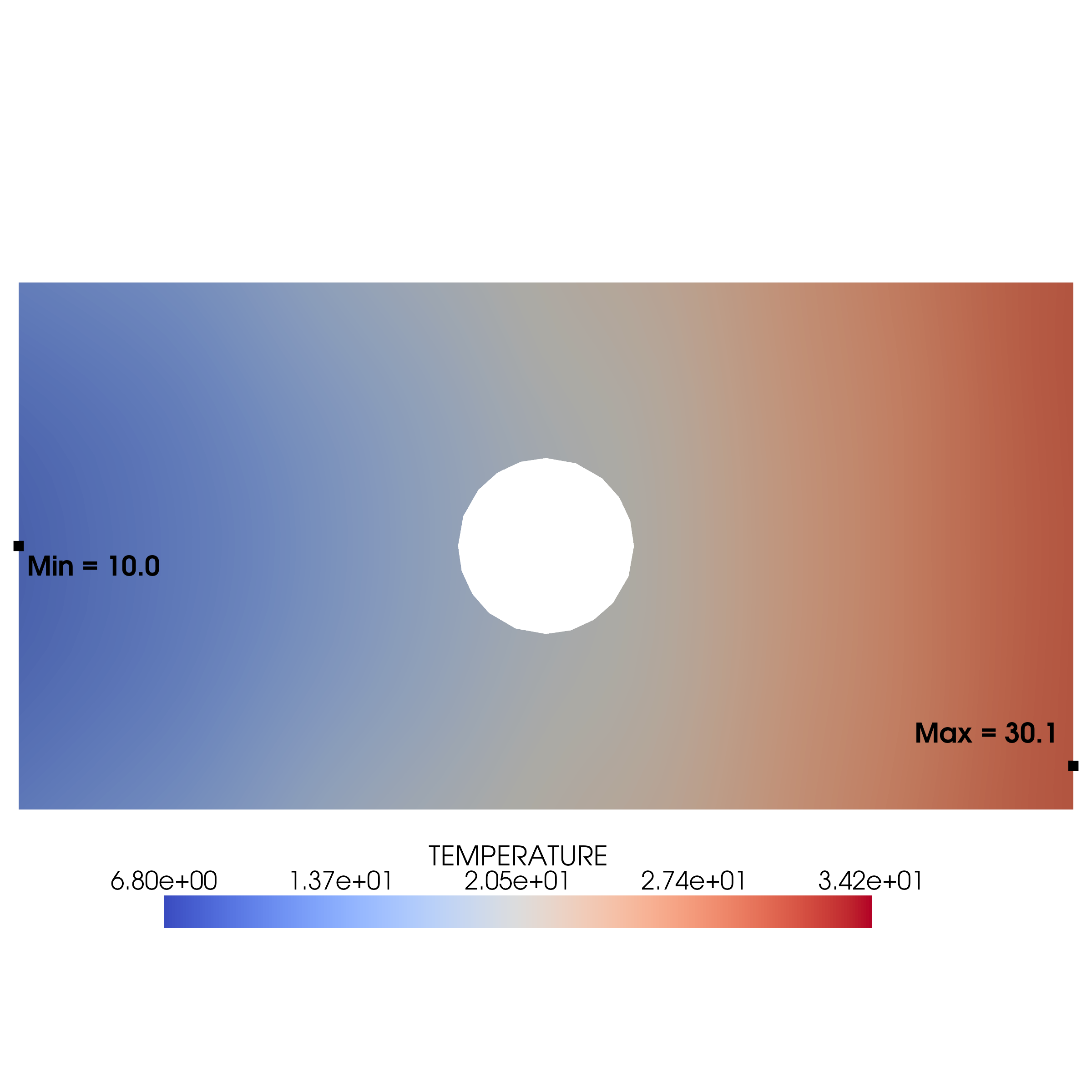}
        \end{minipage}   
        \caption{Identified Young's moduli (left) and \textit{identified} temperature distribution (right) when the thermal field is reconstructed during \acrshort{SI}: \textit{Partitioned} approach.}
        \label{f:plate_linear_6s_c}
    \end{subfigure}
    \caption{\textbf{Plate With Hole. 14 displacement and 6 temperature sensors configuration. Linearly varying thermal field.} Identified Young's moduli and temperature distributions when the thermal load is considered using different approaches during \acrshort{SI}. Peak temperatures are noted in the figures. }
\label{f:plate_linear_6s}
\end{figure} 

\begin{figure}[!t]
\centering
\begin{minipage}[t]{0.49\linewidth}
\centering
\includegraphics[width=\linewidth]{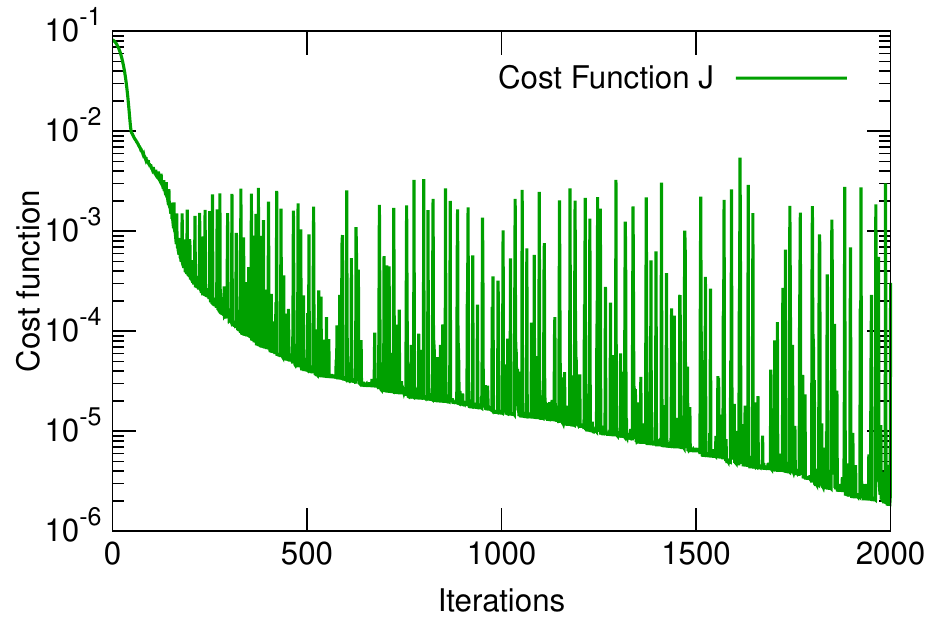}
\subcaption{When the temperature field is \textit{interpolated}.}
\label{f:plate_linear_6s_conv_a}
\end{minipage}
\hfill
\begin{minipage}[t]{0.49\linewidth}
\centering
\includegraphics[width=\linewidth]{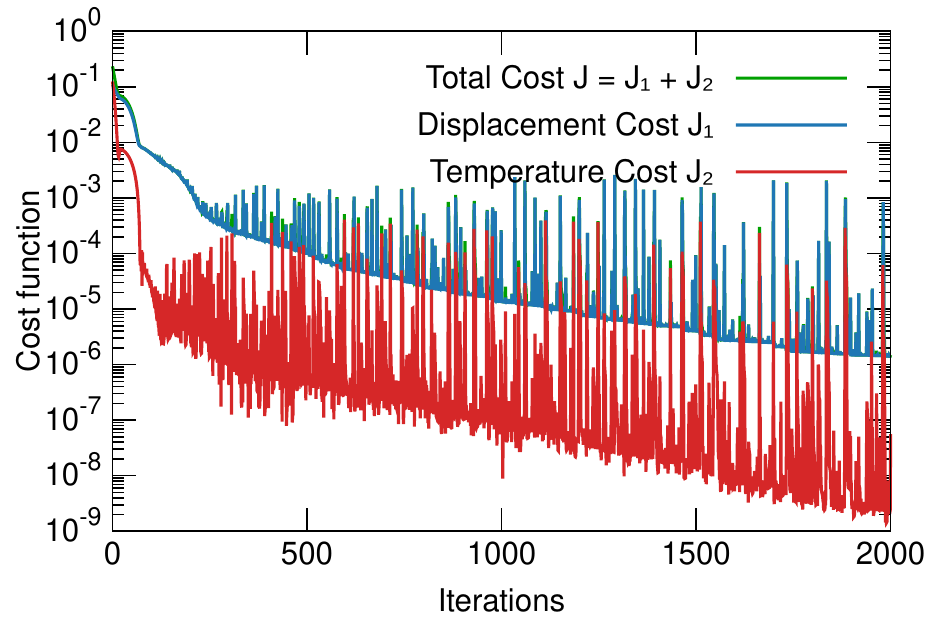}
\subcaption{When the temperature field is \textit{identified}: \textit{Monolithic} approach.}
\label{f:plate_linear_6s_conv_b}
\end{minipage}
\begin{minipage}[t]{0.49\linewidth}
\centering
\includegraphics[width=\linewidth]{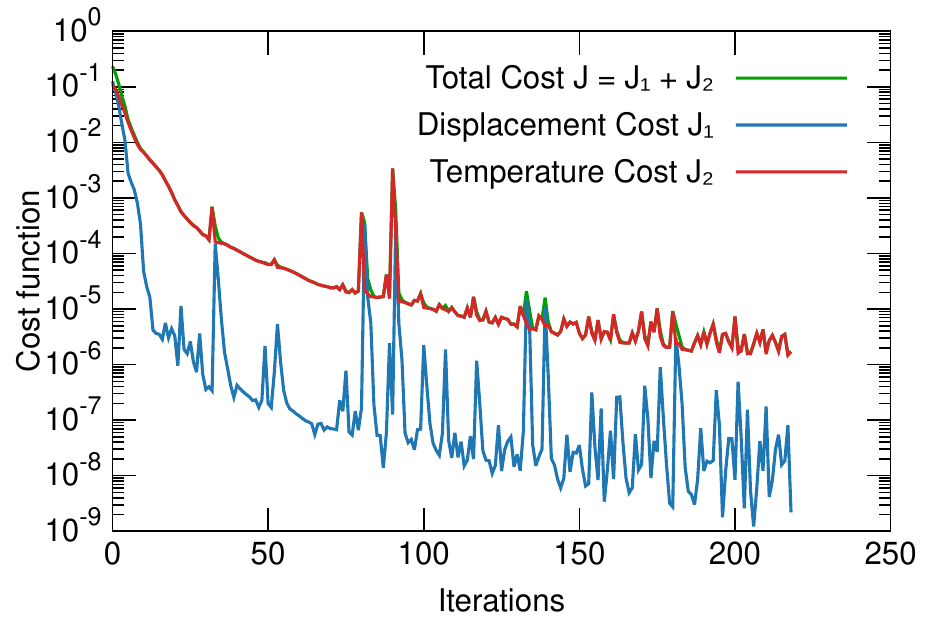}
\subcaption{When the temperature field is \textit{identified}: \textit{Partitioned} approach. Coupling iteration: 218, Overall iteration: 3990.}
\label{f:plate_linear_6s_conv_c}
\end{minipage}
\caption{\textbf{Plate With Hole. 14 displacement and 6 temperature sensors configuration. Linearly varying thermal field.}
Convergence plots when the thermal load is considered using different approaches during \acrshort{SI}.}
\label{f:plate_linear_6s_conv}
\end{figure}

\begin{figure}[!b]
    \centering
    \begin{subfigure}[t]{\textwidth}
        \centering
        \begin{minipage}[t]{0.49\textwidth}
            \centering
            \includegraphics[trim=0 300 0 420, clip, width=\textwidth]{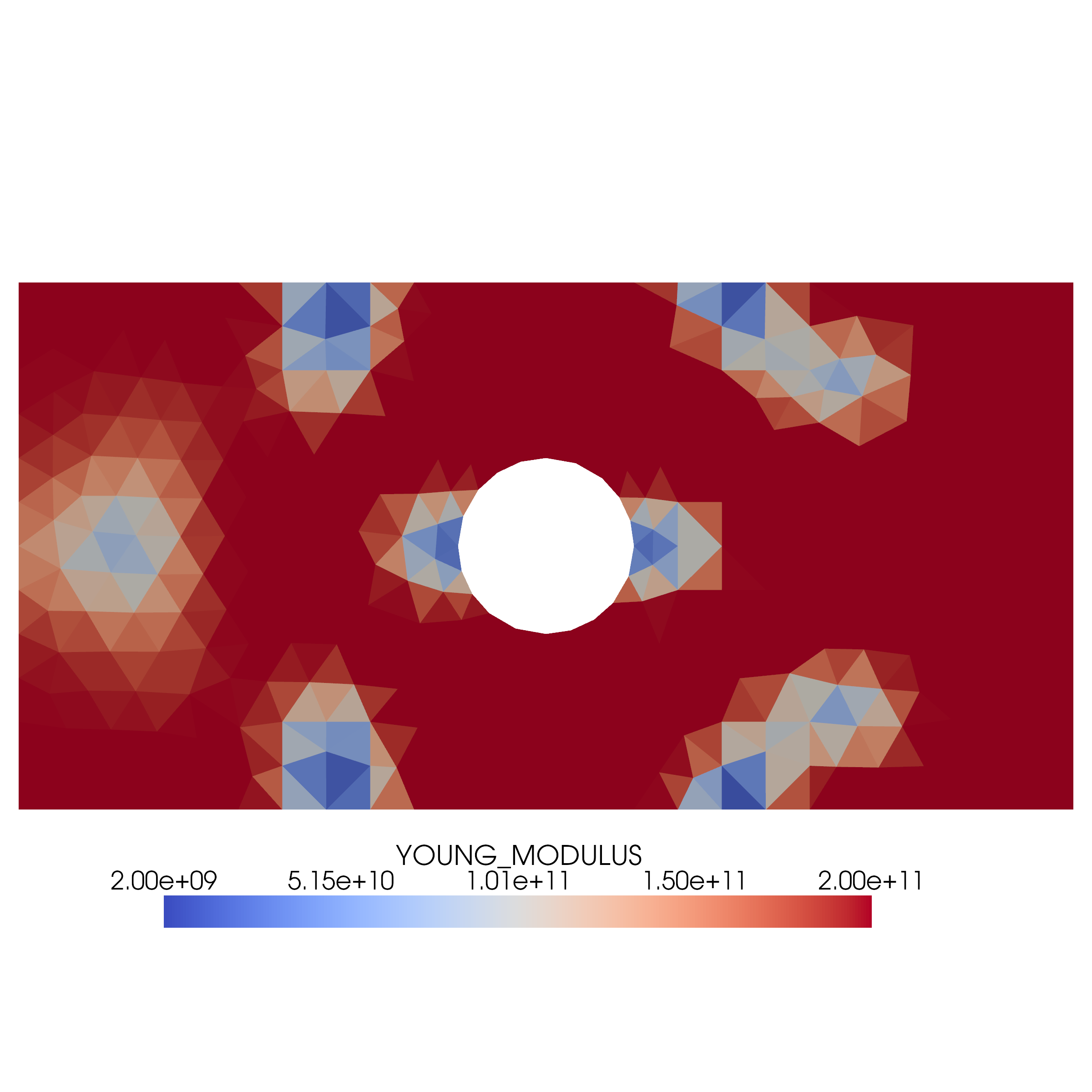}
        \end{minipage}
        \hfill
        \begin{minipage}[t]{0.49\textwidth}
            \centering
            \includegraphics[trim=0 300 0 420, clip, width=\textwidth]{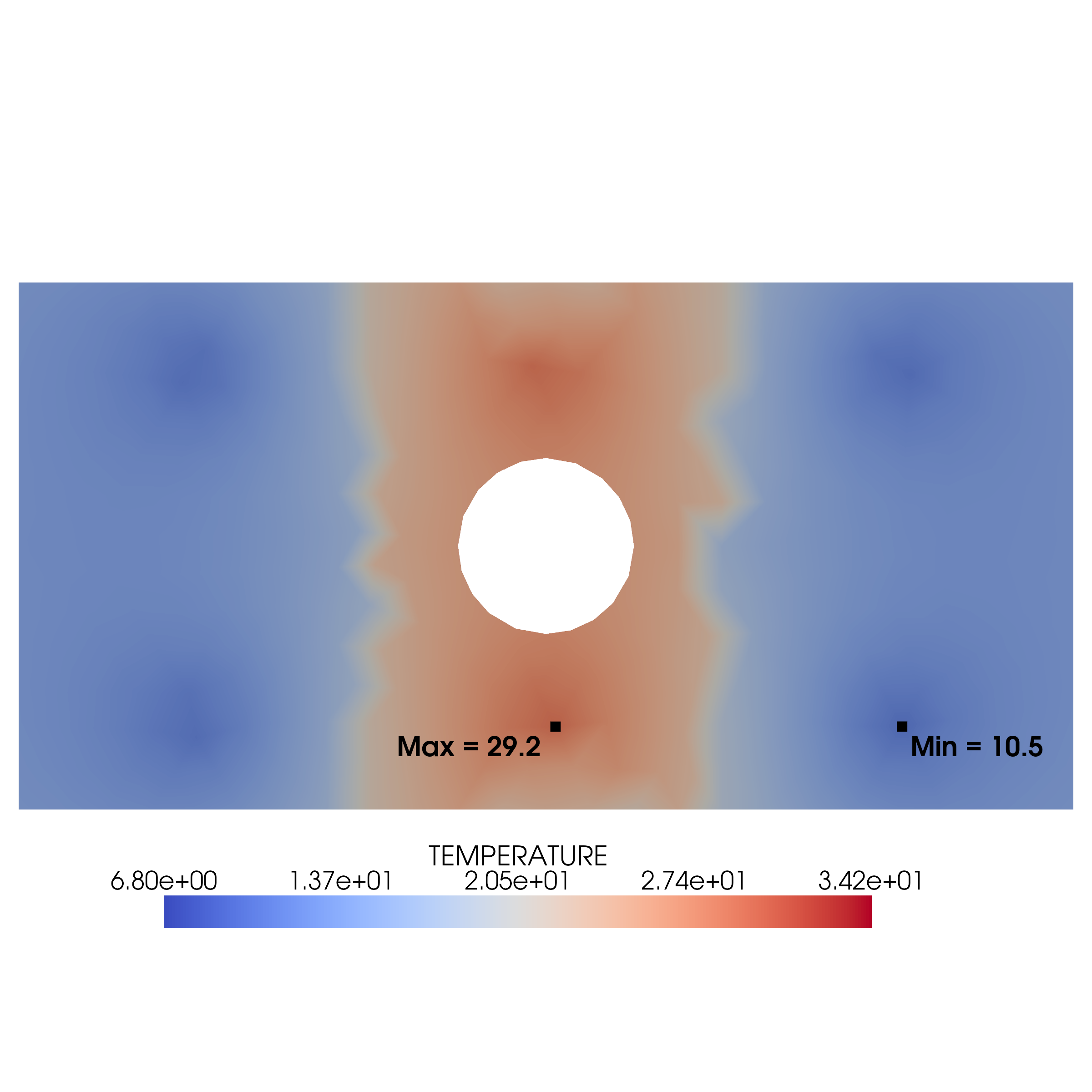}
        \end{minipage}
        \caption{Identified Young's moduli (left) and the temperature distribution (right) when the thermal field is \textit{interpolated}.}
        \label{f:plate_local_6s_a}
    \end{subfigure}
    \begin{subfigure}[t]{\textwidth}
        \centering
        \begin{minipage}[t]{0.49\textwidth}
            \centering
            \includegraphics[trim=0 300 0 420, clip,width=\textwidth]{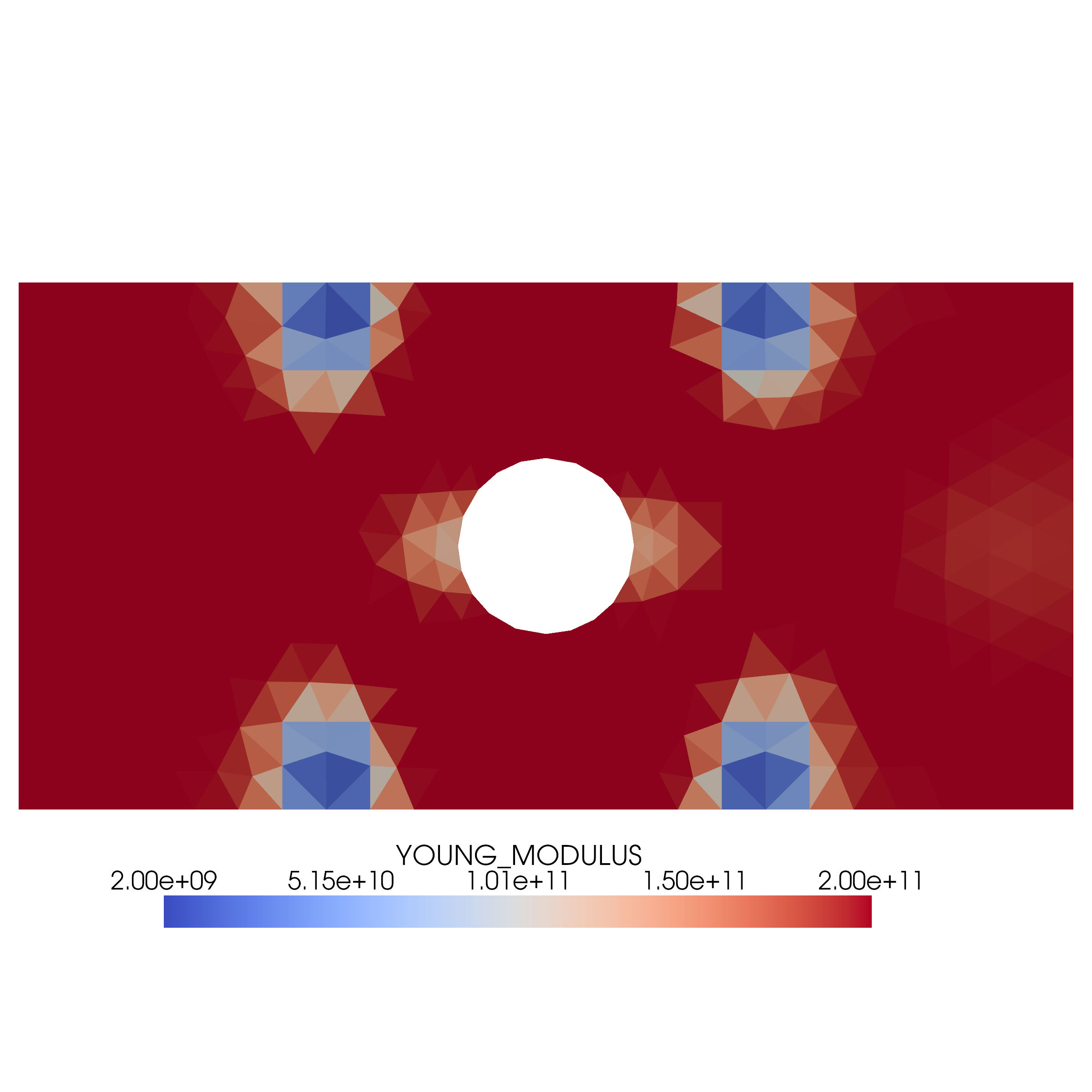}
        \end{minipage}
        \begin{minipage}[t]{0.49\textwidth}
            \centering
            \includegraphics[trim=0 300 0 420, clip,width=\textwidth]{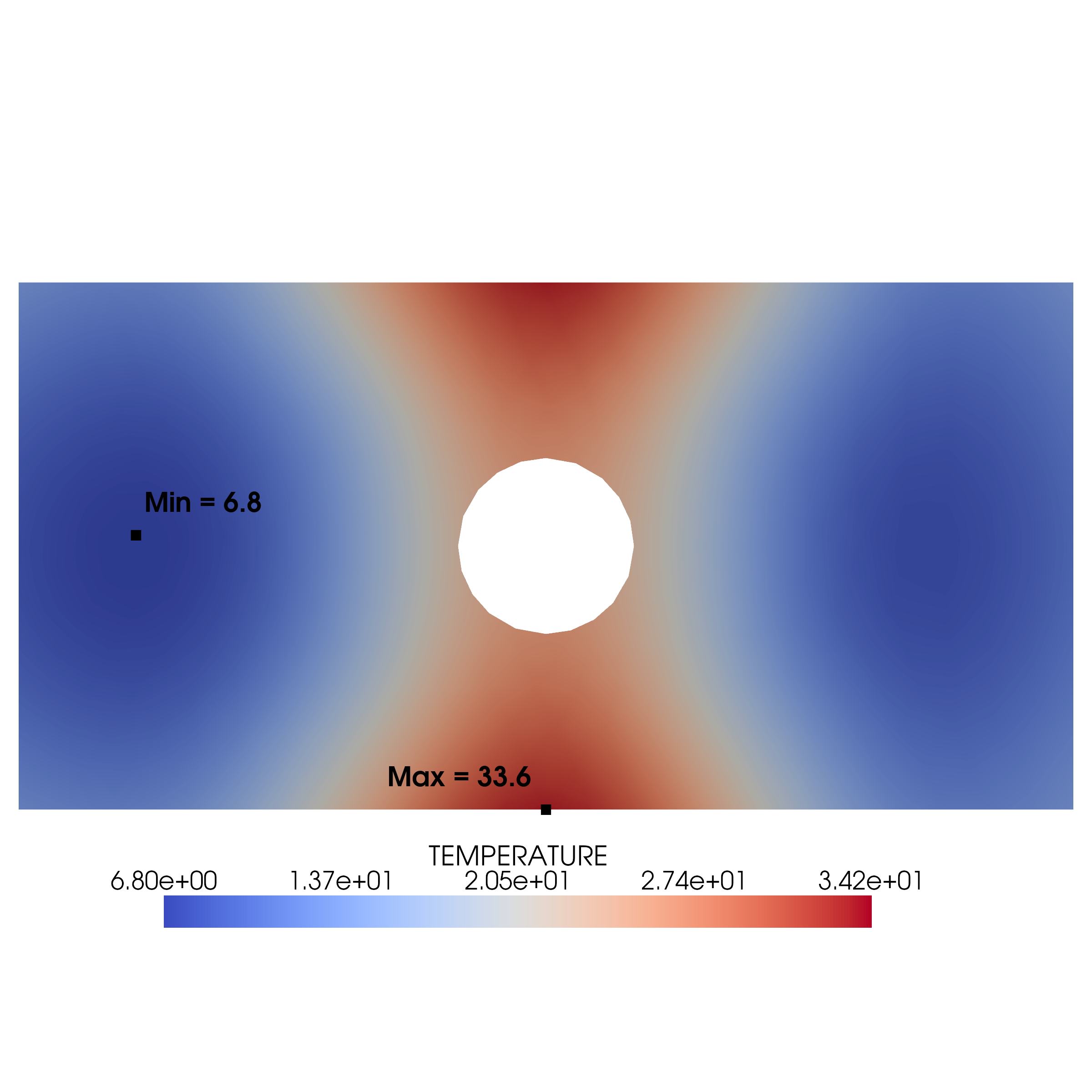}
        \end{minipage}  
        \caption{Identified Young's moduli (left) and \textit{identified} temperature distribution (right) when the thermal field is reconstructed during \acrshort{SI}: \textit{Monolithic} approach.}
        \label{f:plate_local_6s_b}
    \end{subfigure}
    \begin{subfigure}[t]{\textwidth}
        \centering
        \begin{minipage}[t]{0.49\textwidth}
            \centering
            \includegraphics[trim=0 300 0 420, clip,width=\textwidth]{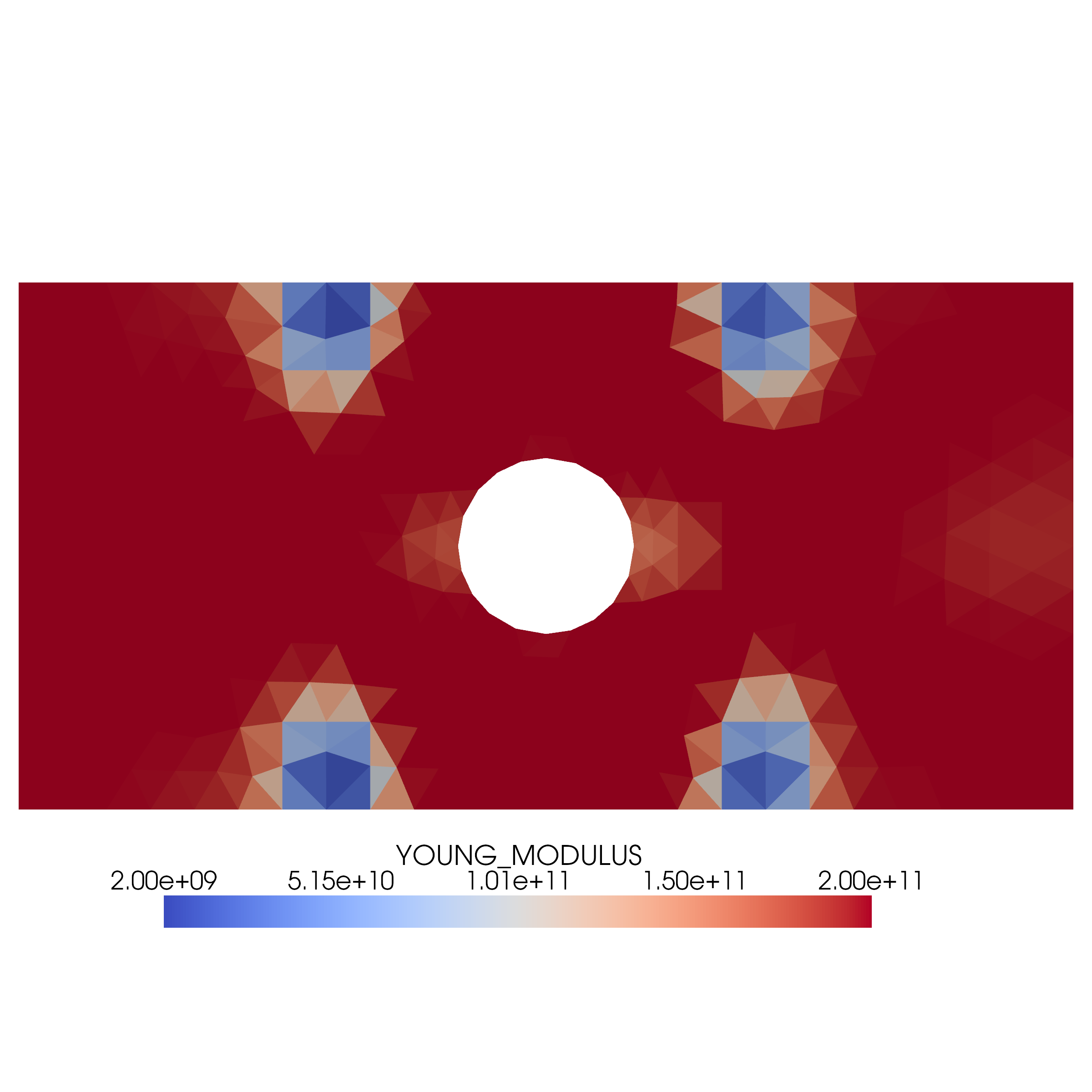}
        \end{minipage}
        \begin{minipage}[t]{0.49\textwidth}
            \centering
            \includegraphics[trim=0 300 0 420, clip,width=\textwidth]{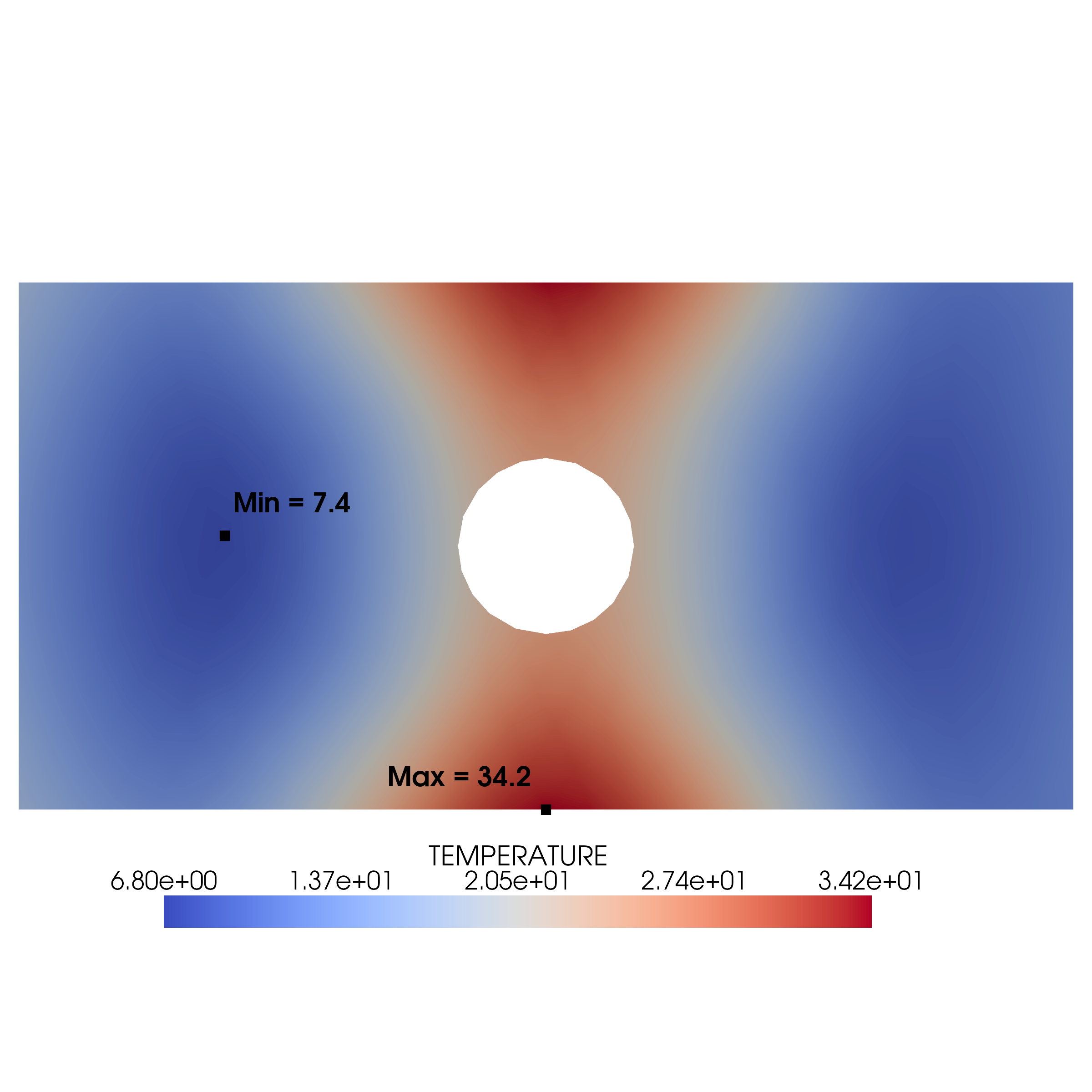}
        \end{minipage}   
        \caption{Identified Young's moduli (left) and \textit{identified} temperature distribution (right) when the thermal field is reconstructed during \acrshort{SI}: \textit{Partitioned} approach.}
        \label{f:plate_local_6s_c}
    \end{subfigure}
    \caption{\textbf{Plate With Hole. 14 displacement and 6 temperature sensors configuration. Localized thermal field.} Identified Young's moduli and temperature distributions when the thermal load is considered using different approaches during \acrshort{SI}. Peak temperatures are noted in the figures. }
\label{f:plate_local_6s}
\end{figure} 

\begin{figure}[!t]
\centering
\begin{minipage}[t]{0.49\linewidth}
\centering
\includegraphics[width=\linewidth]{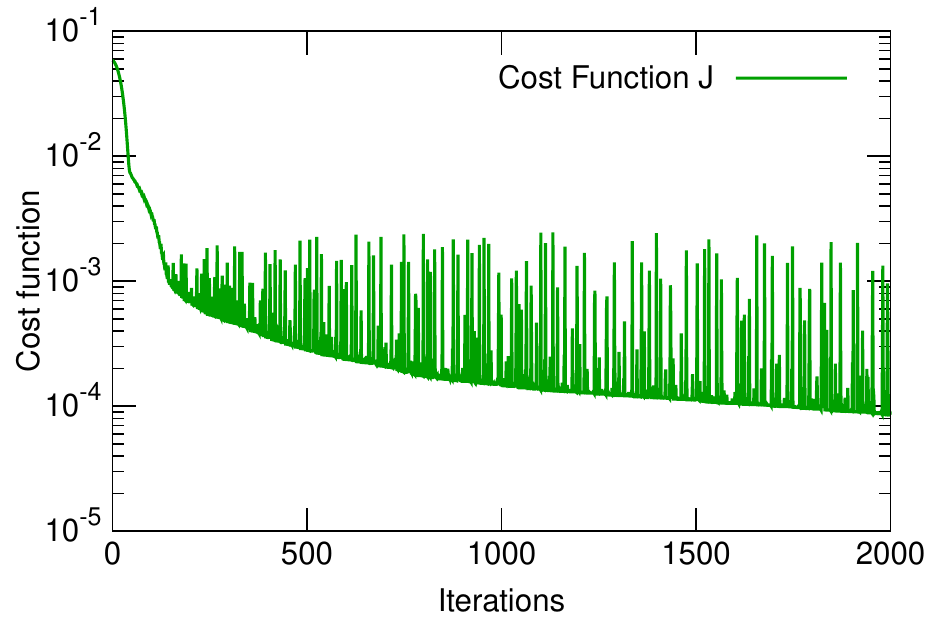}
\subcaption{When the temperature field is \textit{interpolated}.}
\label{f:plate_local_6s_conv_a}
\end{minipage}
\hfill
\begin{minipage}[t]{0.49\linewidth}
\centering
\includegraphics[width=\linewidth]{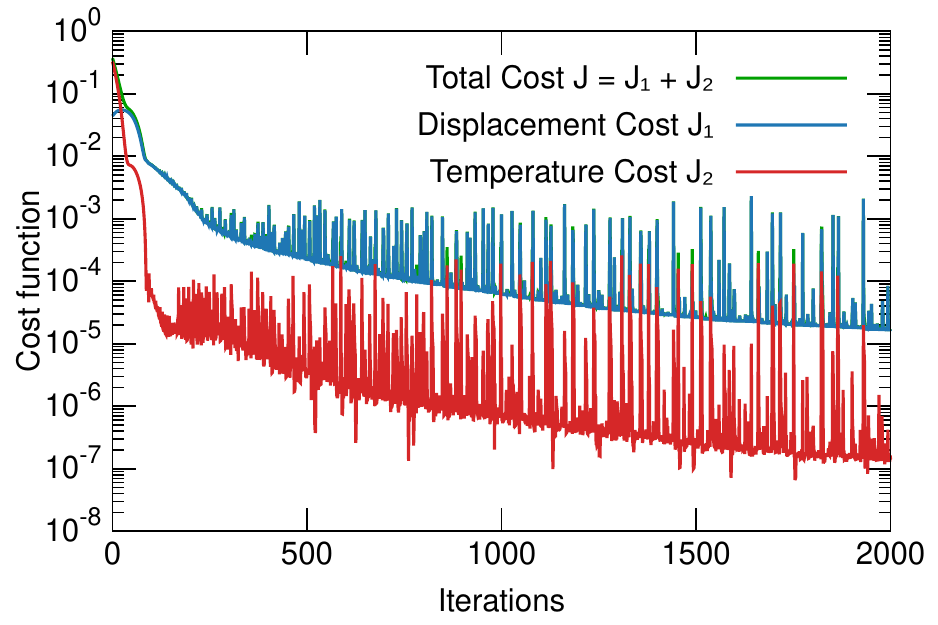}
\subcaption{When the temperature field is \textit{identified}: \textit{Monolithic} approach.}
\label{f:plate_local_6s_conv_b}
\end{minipage}
\begin{minipage}[t]{0.49\linewidth}
\centering
\includegraphics[width=\linewidth]{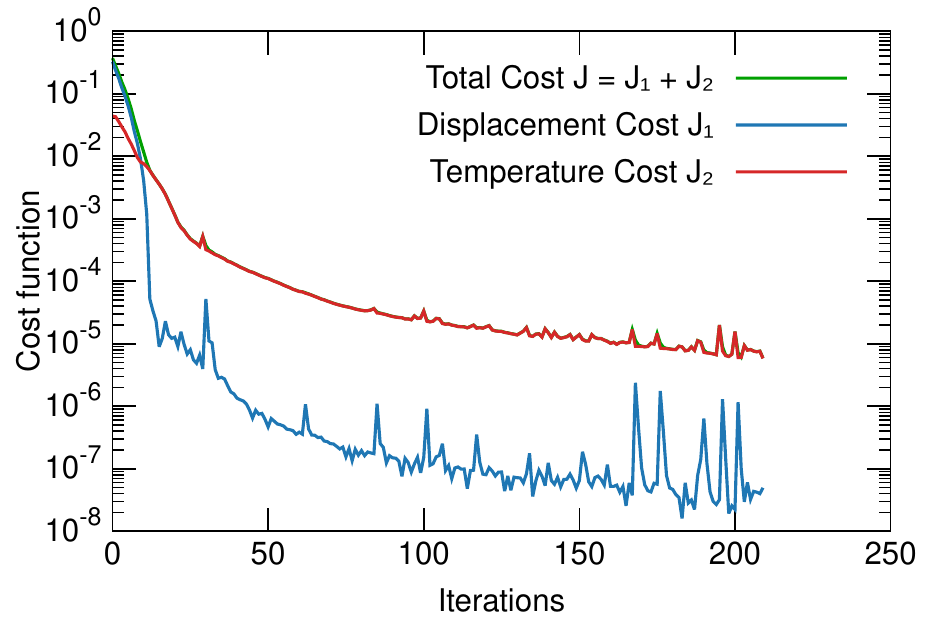}
\subcaption{When the temperature field is \textit{identified}: \textit{Partitioned} approach. Coupling iteration: 209, Overall iteration: 4009.}
\label{f:plate_local_6s_conv_c}
\end{minipage}
\caption{\textbf{Plate With Hole. 14 displacement and 6 temperature sensors configuration. Localized thermal field.}
Convergence plots when the thermal load is considered using different approaches during \acrshort{SI}.}
\label{f:plate_local_6s_conv}
\end{figure}

\begin{figure}[!b]
    \centering
    \begin{subfigure}[t]{\textwidth}
        \centering
        \begin{minipage}[t]{0.49\textwidth}
            \centering
            \includegraphics[trim=0 300 0 420, clip, width=\textwidth]{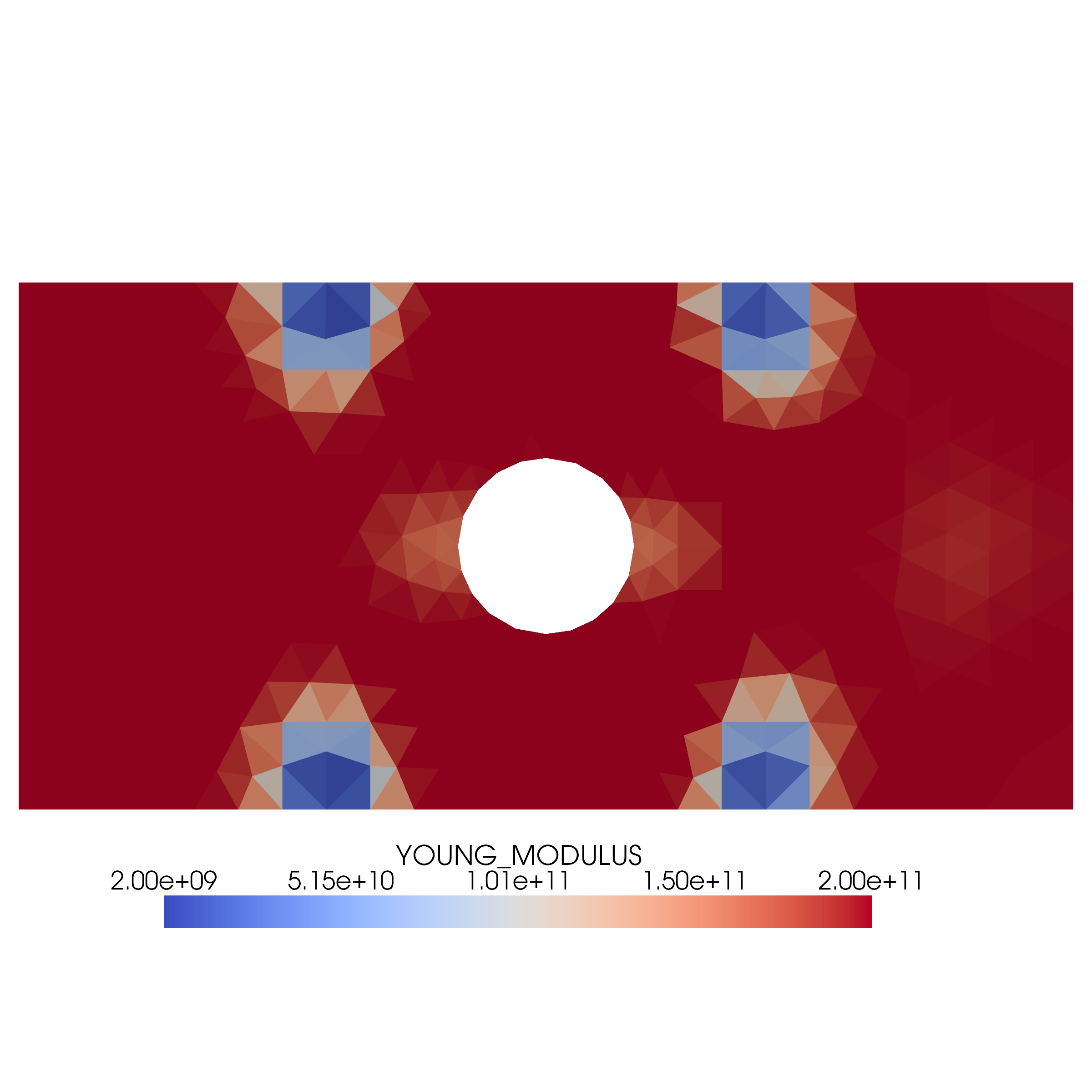}
        \end{minipage}
        \hfill
        \begin{minipage}[t]{0.49\textwidth}
            \centering
            \includegraphics[trim=0 300 0 420, clip, width=\textwidth]{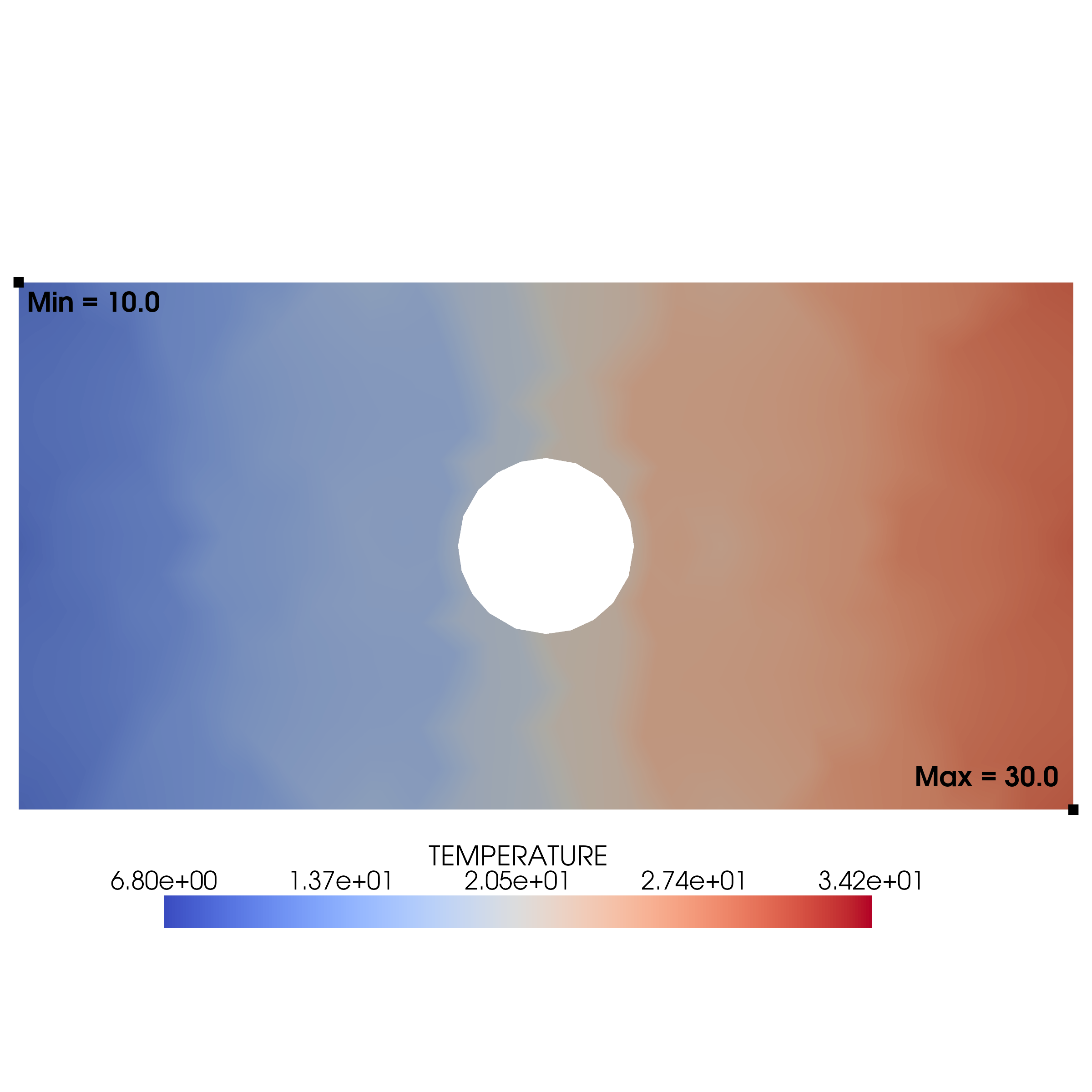}
        \end{minipage}
        \caption{Identified Young's moduli (left) and the temperature distribution (right) when the thermal field is \textit{interpolated}.}
        \label{f:plate_linear_16s_a}
    \end{subfigure}
    \begin{subfigure}[t]{\textwidth}
        \centering
        \begin{minipage}[t]{0.49\textwidth}
            \centering
            \includegraphics[trim=0 300 0 420, clip,width=\textwidth]{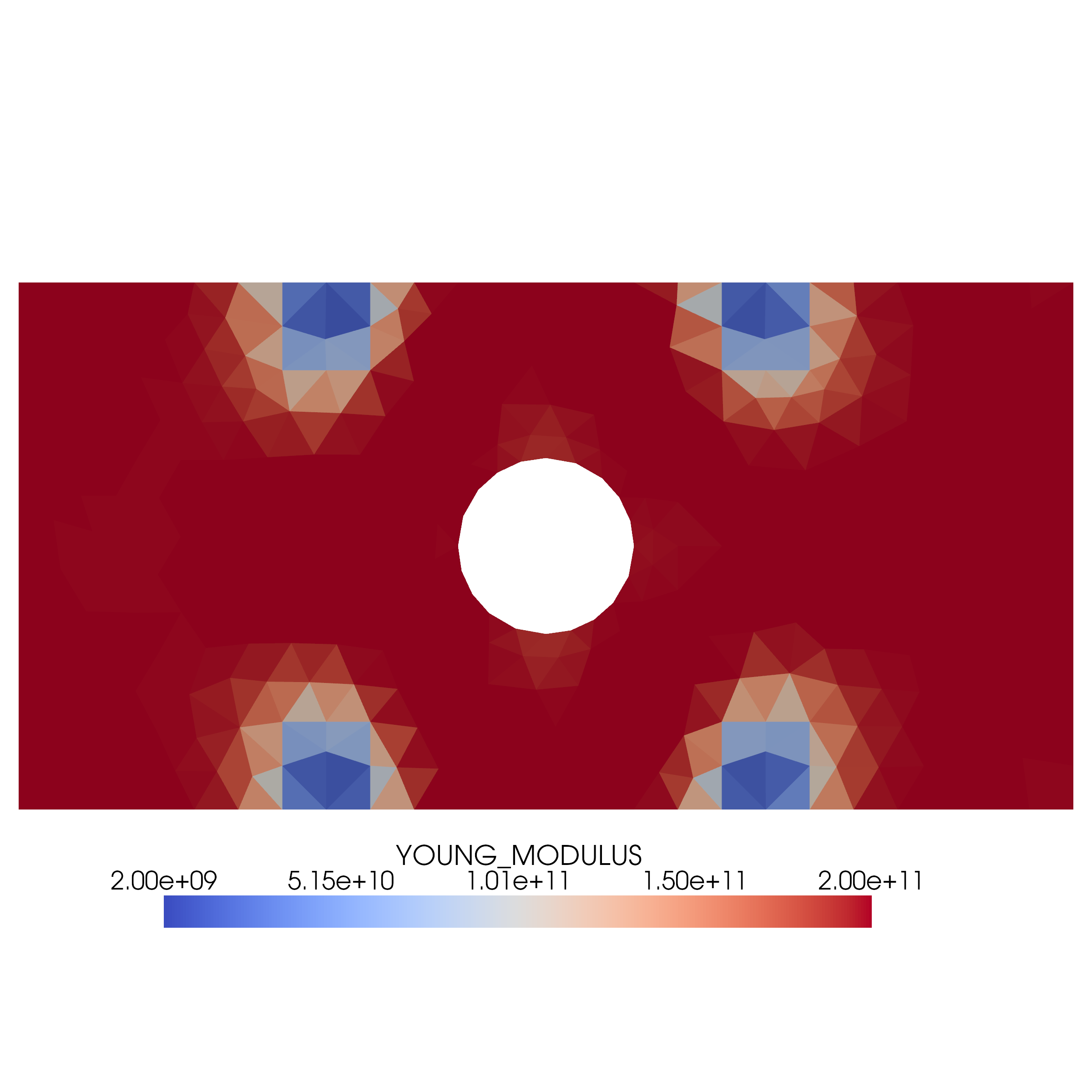}
        \end{minipage}
        \begin{minipage}[t]{0.49\textwidth}
            \centering
            \includegraphics[trim=0 300 0 420, clip,width=\textwidth]{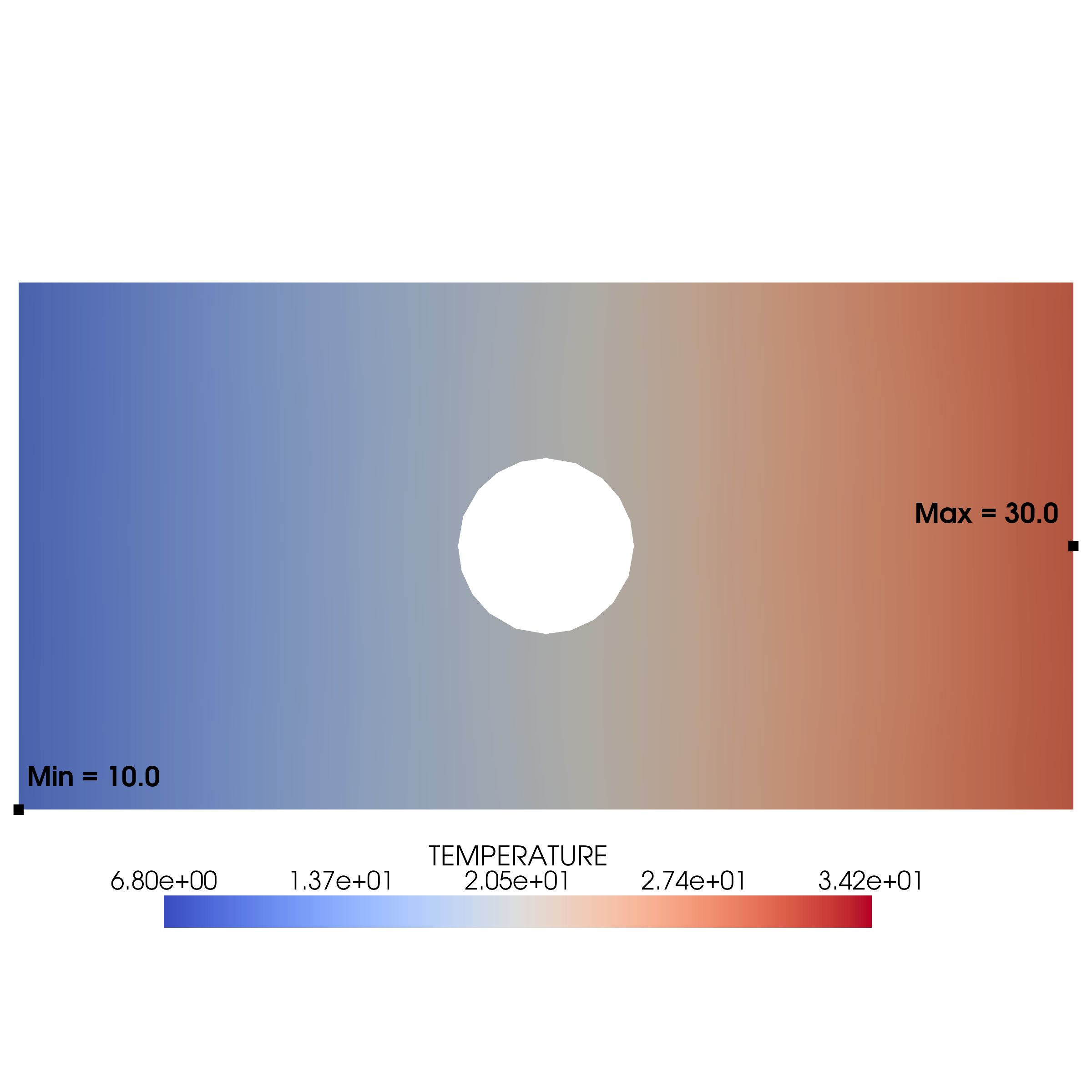}
        \end{minipage}  
        \caption{Identified Young's moduli (left) and \textit{identified} temperature distribution (right) when the thermal field is reconstructed during \acrshort{SI}: \textit{Monolithic} approach.}
        \label{f:plate_linear_16s_b}
    \end{subfigure}
    \begin{subfigure}[t]{\textwidth}
        \centering
        \begin{minipage}[t]{0.49\textwidth}
            \centering
            \includegraphics[trim=0 300 0 420, clip,width=\textwidth]{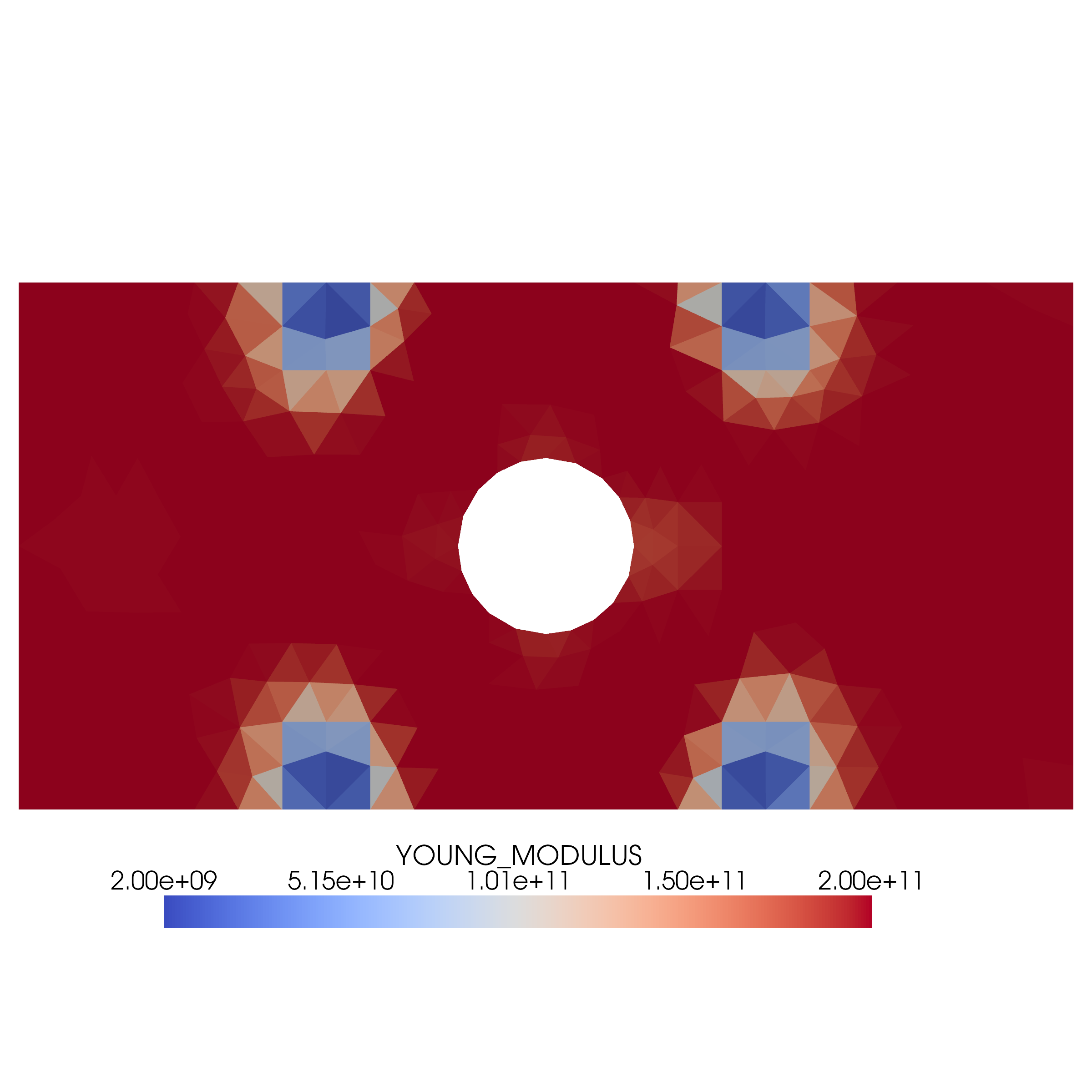}
        \end{minipage}
        \begin{minipage}[t]{0.49\textwidth}
            \centering
            \includegraphics[trim=0 300 0 420, clip,width=\textwidth]{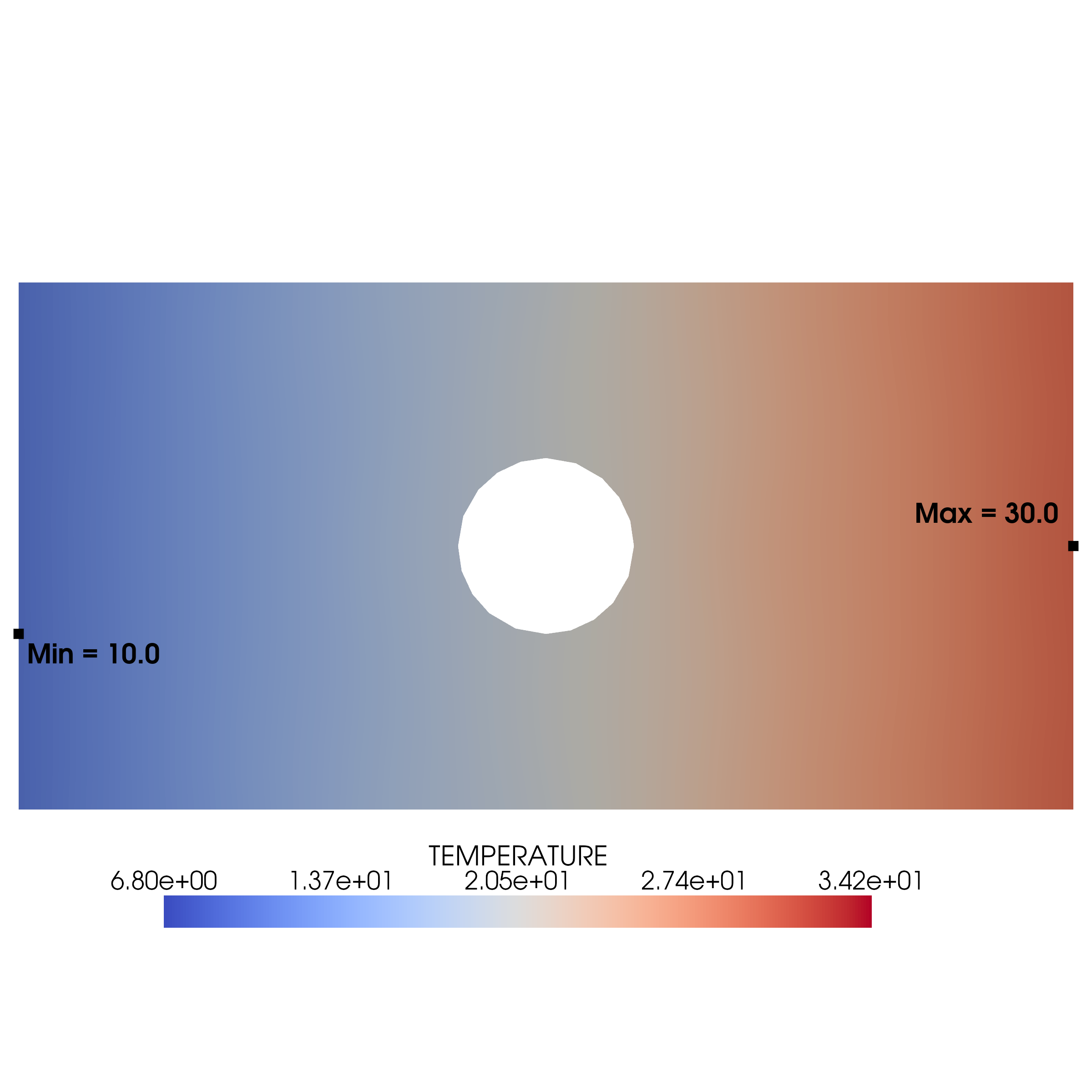}
        \end{minipage}  
        \caption{Identified Young's moduli (left) and \textit{identified} temperature distribution (right) when the thermal field is reconstructed during \acrshort{SI}: \textit{Partitioned} approach.}
        \label{f:plate_linear_16s_c}
    \end{subfigure}
    \caption{\textbf{Plate With Hole. 14 displacement and 16 temperature sensors configuration. Linearly varying thermal field.} Identified Young's moduli and temperature distributions when the thermal load is considered using different approaches during \acrshort{SI}. Peak temperatures are noted in the figures. }
\label{f:plate_linear_16s}
\end{figure} 

\begin{figure}[!t]
\centering
\begin{minipage}[t]{0.49\linewidth}
\centering
\includegraphics[width=\linewidth]{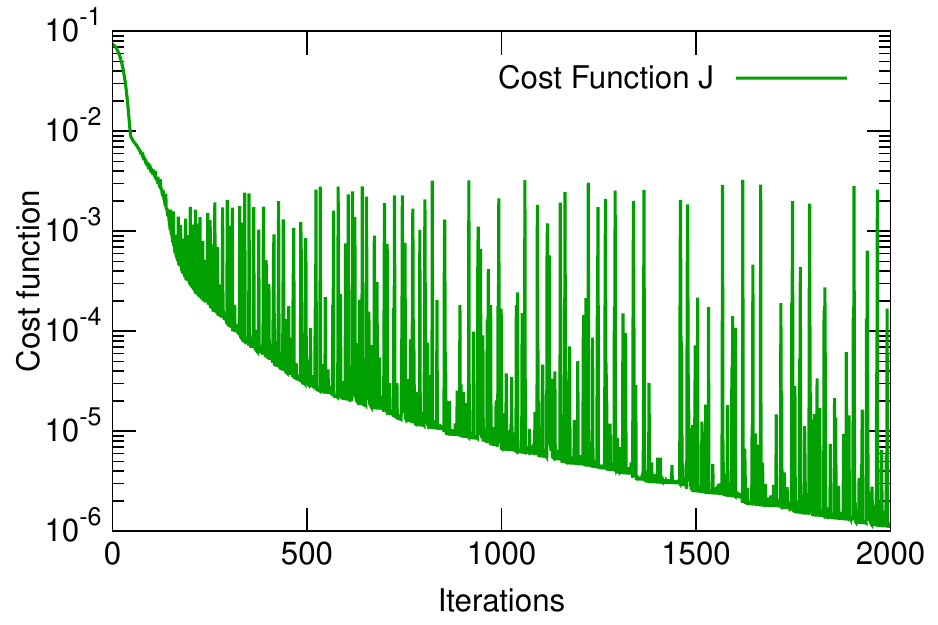}
\subcaption{When the temperature field is \textit{interpolated}.}
\label{f:plate_linear_16s_conv_a}
\end{minipage}
\hfill
\begin{minipage}[t]{0.49\linewidth}
\centering
\includegraphics[width=\linewidth]{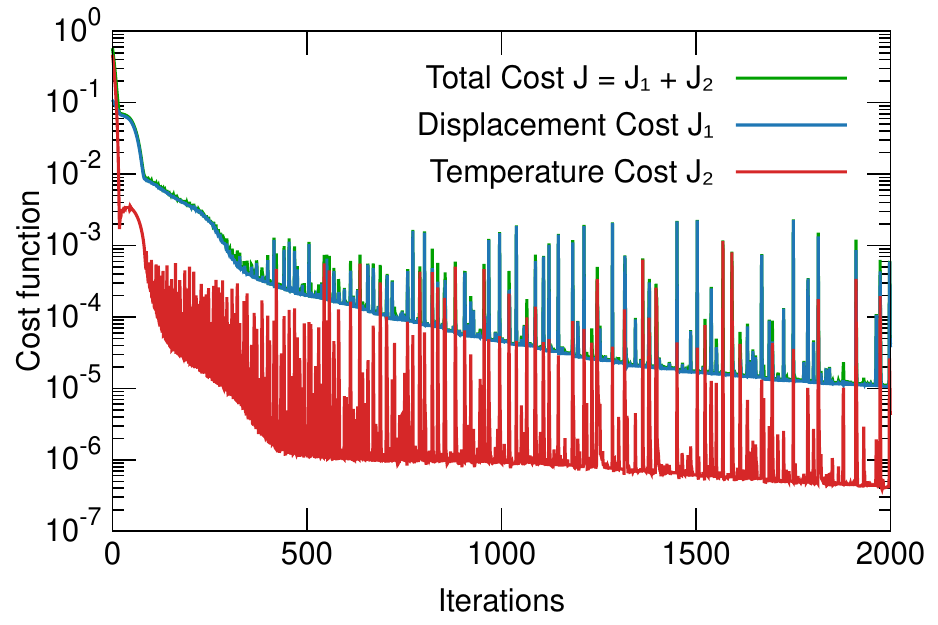}
\subcaption{When the temperature field is \textit{identified}: \textit{Monolithic} approach.}
\label{f:plate_linear_16s_conv_b}
\end{minipage}
\begin{minipage}[t]{0.49\linewidth}
\centering
\includegraphics[width=\linewidth]{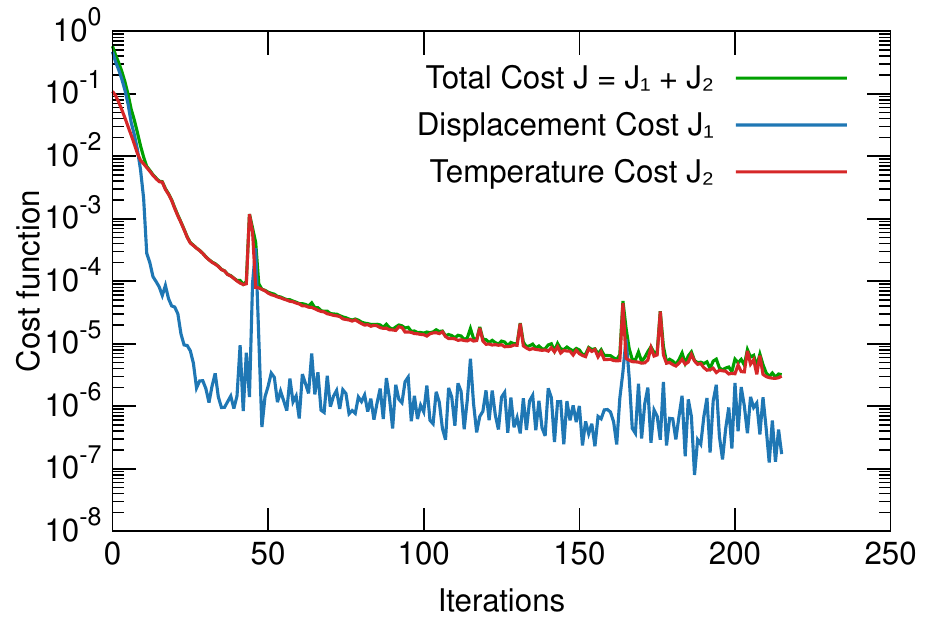}
\subcaption{When the temperature field is \textit{identified}: \textit{Partitioned} approach. Coupling iteration: 215, Overall iteration: 3997.}
\label{f:plate_linear_16s_conv_c}
\end{minipage}
\caption{\textbf{Plate With Hole. 14 displacement and 16 temperature sensors configuration. Linearly varying thermal field.}
Convergence plots when the thermal load is considered using different approaches during \acrshort{SI}.}
\label{f:plate_linear_16s_conv}
\end{figure}

\begin{figure}[!b]
    \centering
    \begin{subfigure}[t]{\textwidth}
        \centering
        \begin{minipage}[t]{0.49\textwidth}
            \centering
            \includegraphics[trim=0 300 0 420, clip, width=\textwidth]{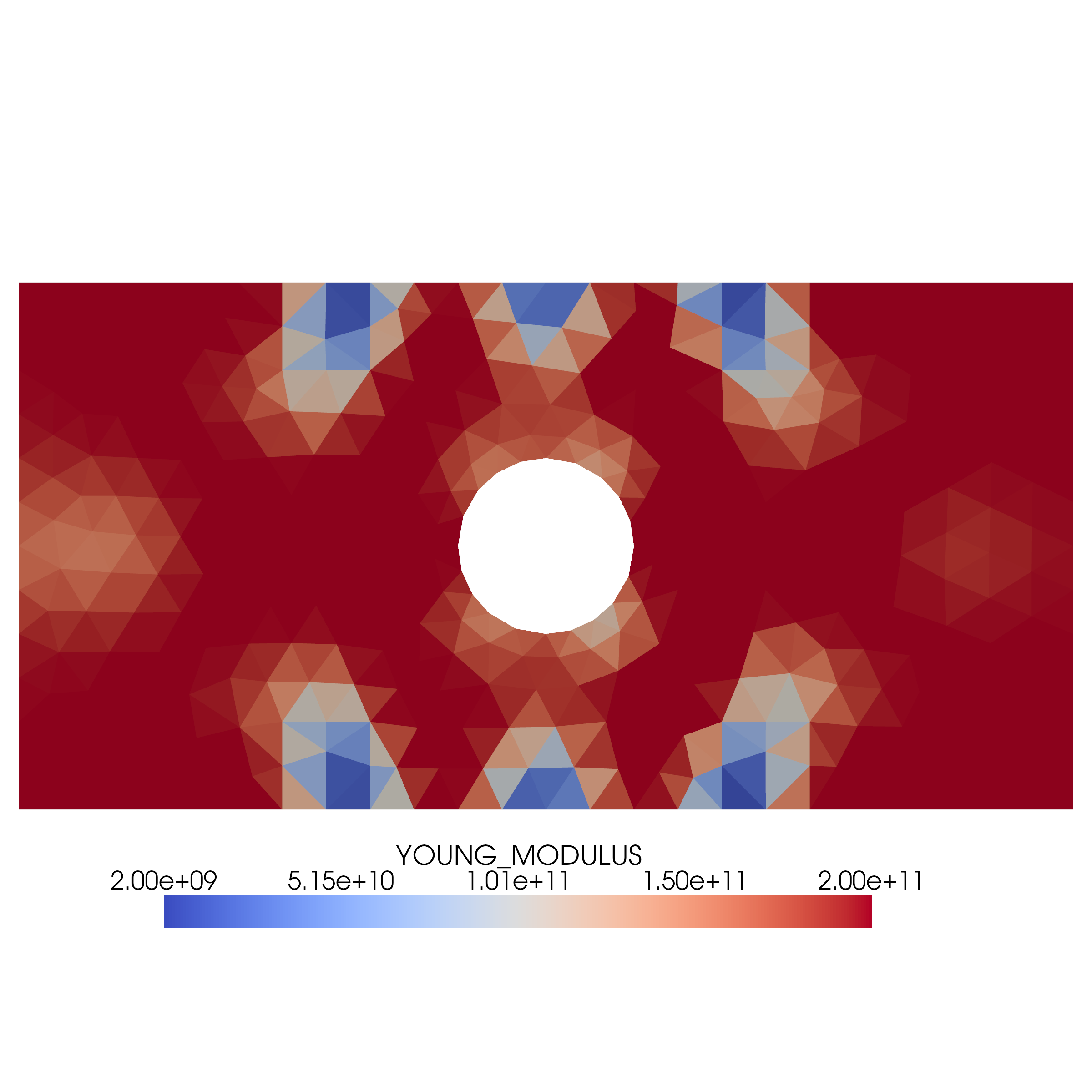}
        \end{minipage}
        \hfill
        \begin{minipage}[t]{0.49\textwidth}
            \centering
            \includegraphics[trim=0 300 0 420, clip, width=\textwidth]{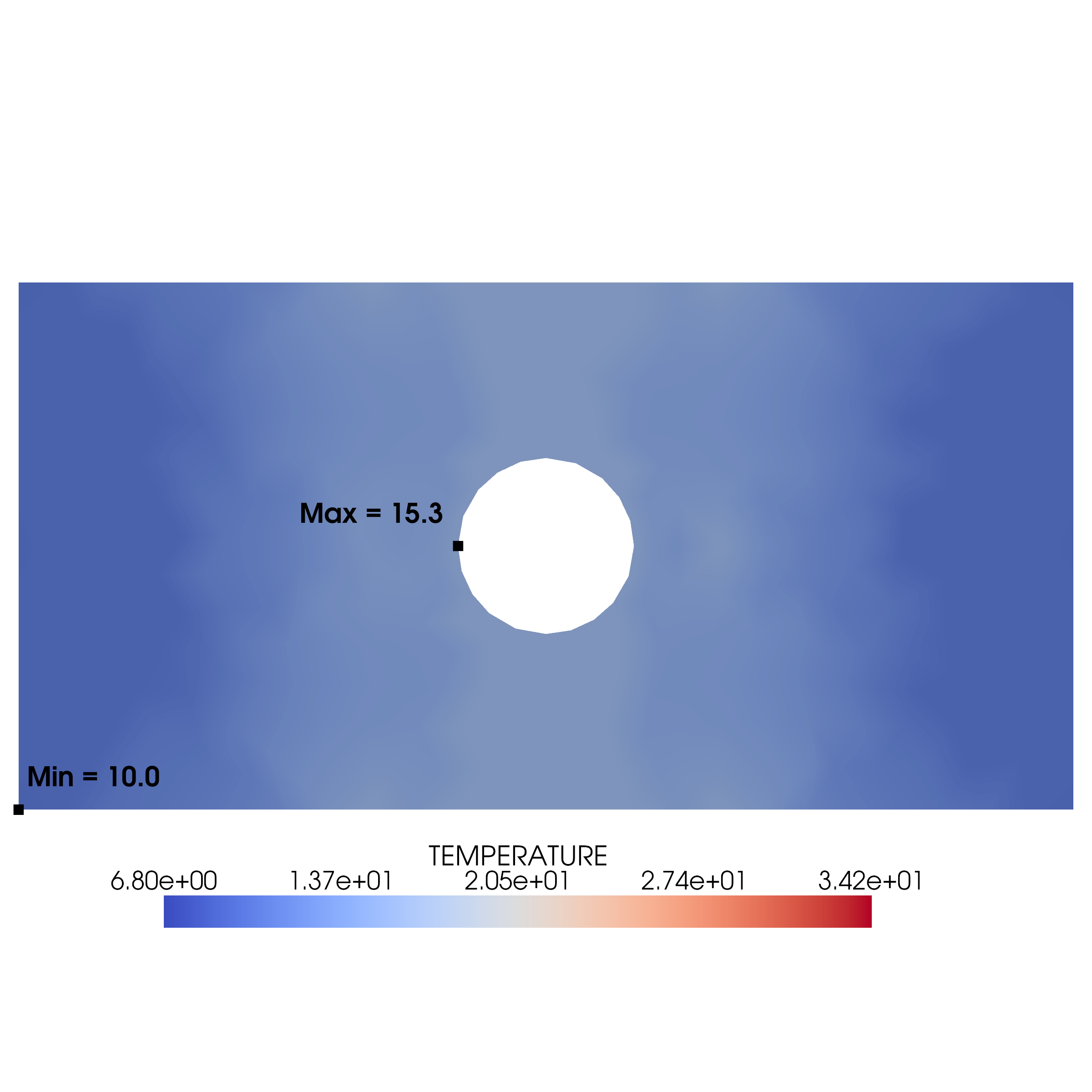}
        \end{minipage}
        \caption{Identified Young's moduli (left) and the temperature distribution (right) when the thermal field is \textit{interpolated}.}
        \label{f:plate_local_16s_a}
    \end{subfigure}
    \begin{subfigure}[t]{\textwidth}
        \centering
        \begin{minipage}[t]{0.49\textwidth}
            \centering
            \includegraphics[trim=0 300 0 420, clip,width=\textwidth]{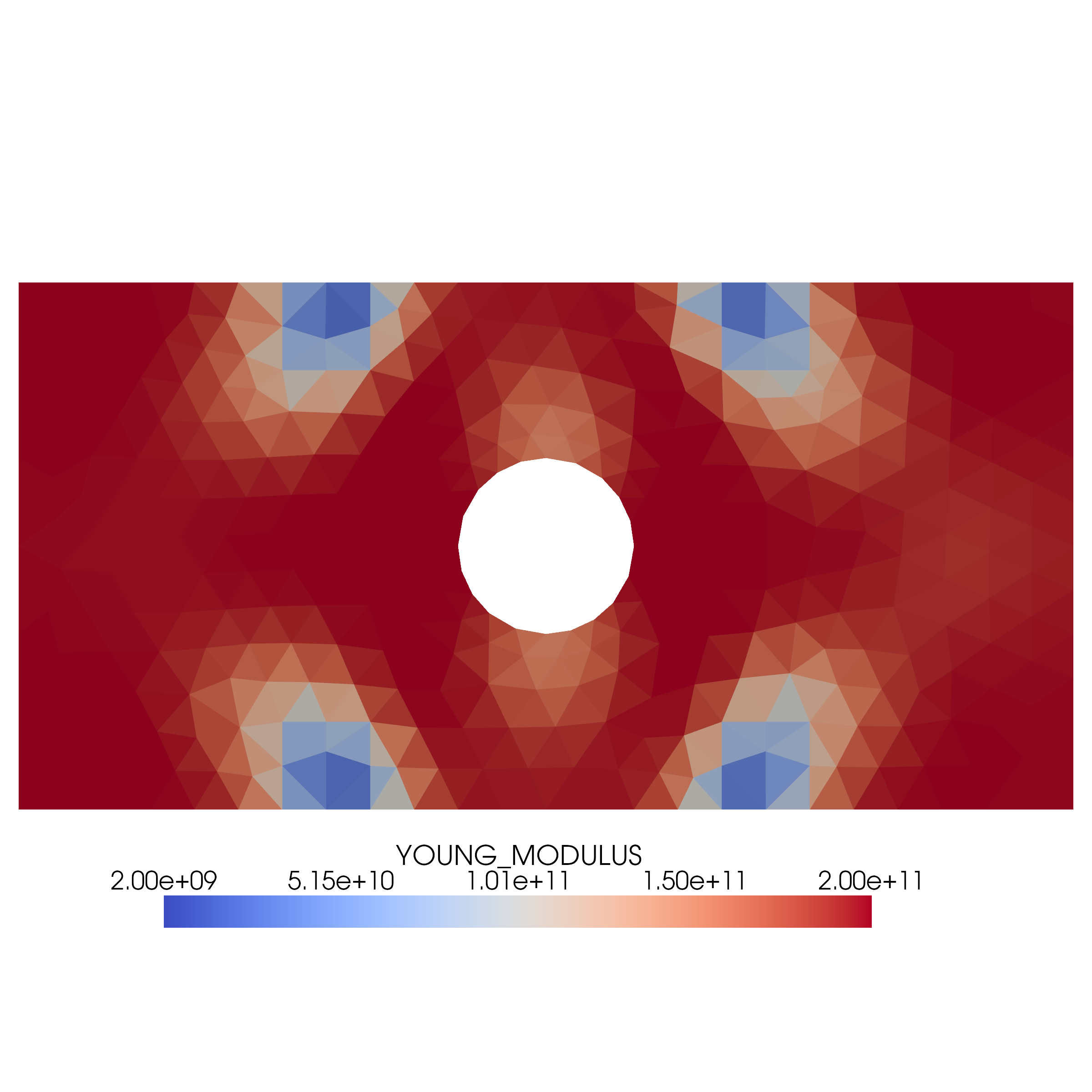}
        \end{minipage}
        \begin{minipage}[t]{0.49\textwidth}
            \centering
            \includegraphics[trim=0 300 0 420, clip,width=\textwidth]{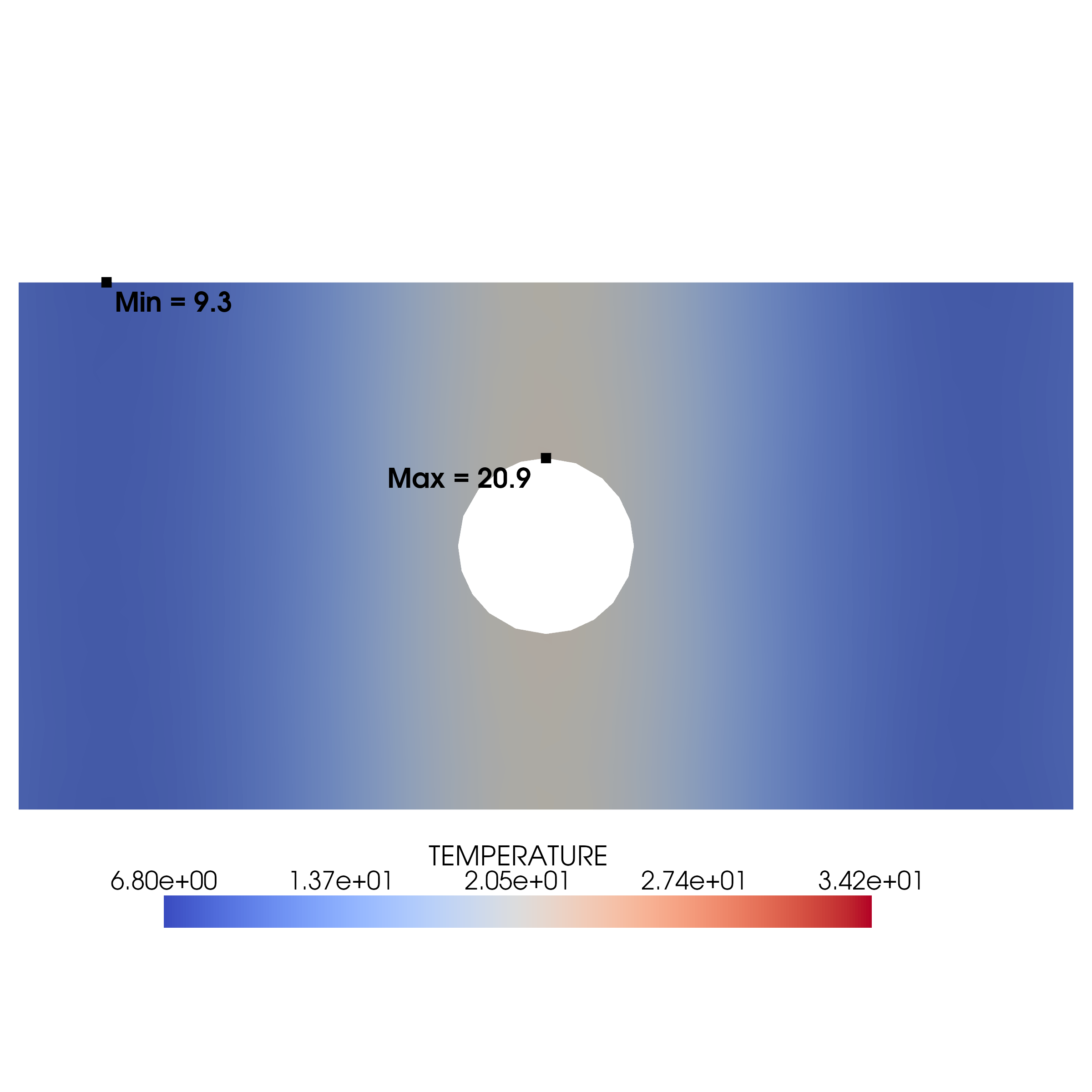}
        \end{minipage}  
        \caption{Identified Young's moduli (left) and \textit{identified} temperature distribution (right) when the thermal field is reconstructed during \acrshort{SI}: \textit{Monolithic} approach.}
        \label{f:plate_local_16s_b}
    \end{subfigure}
    \begin{subfigure}[t]{\textwidth}
        \centering
        \begin{minipage}[t]{0.49\textwidth}
            \centering
            \includegraphics[trim=0 300 0 420, clip,width=\textwidth]{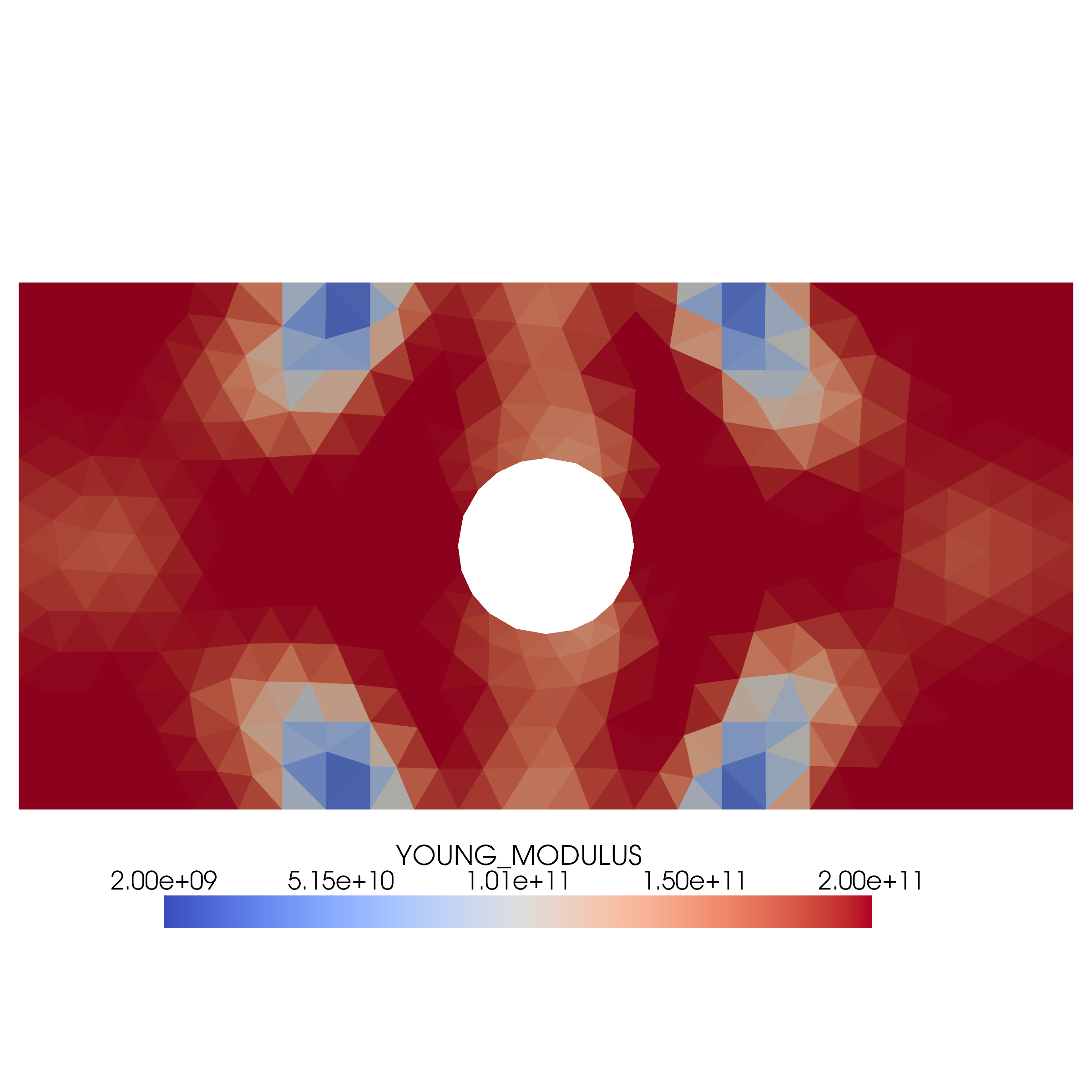}
        \end{minipage}
        \begin{minipage}[t]{0.49\textwidth}
            \centering
            \includegraphics[trim=0 300 0 420, clip,width=\textwidth]{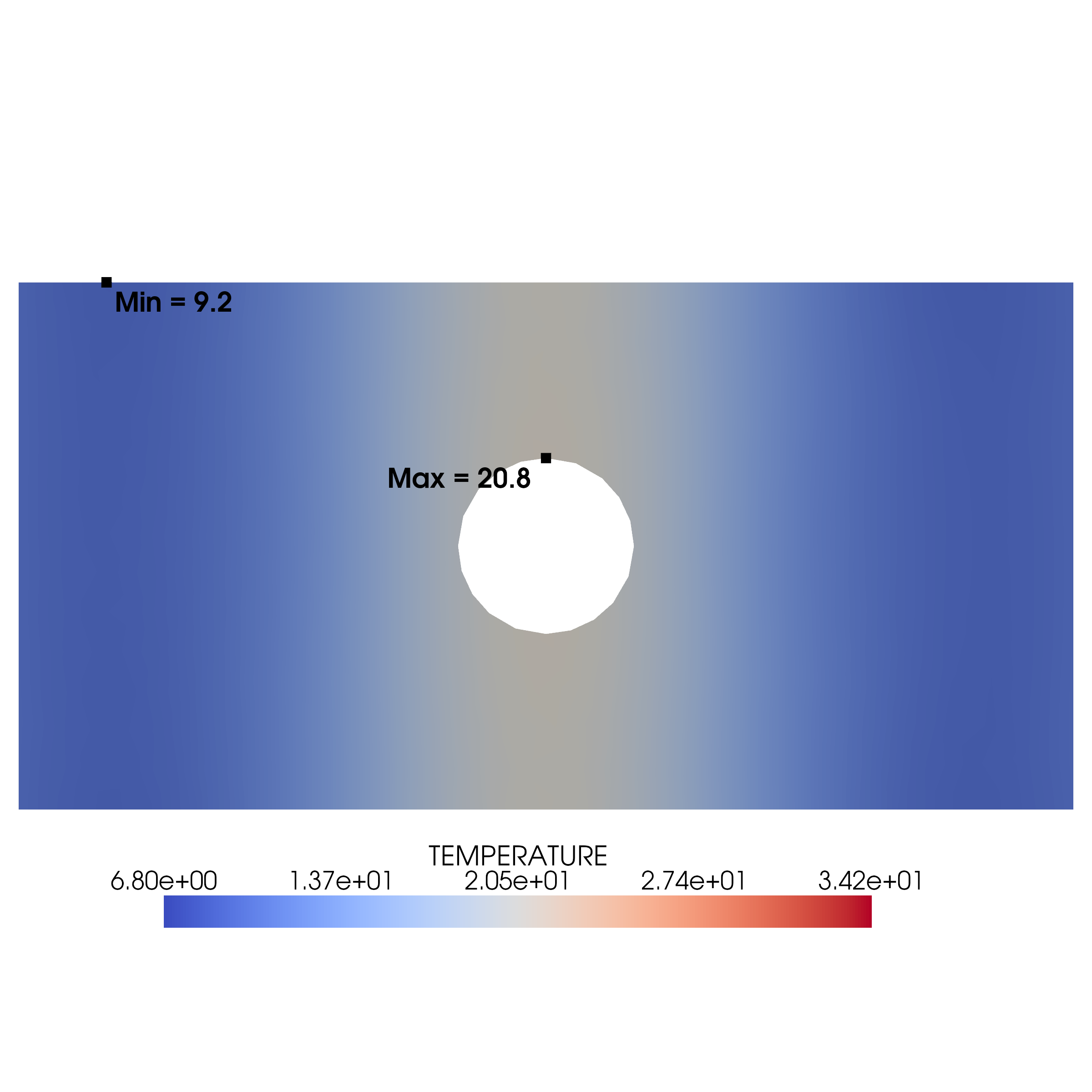}
        \end{minipage}   
        \caption{Identified Young's moduli (left) and \textit{identified} temperature distribution (right) when the thermal field is reconstructed during \acrshort{SI}: \textit{Partitioned} approach.}
        \label{f:plate_local_16s_c}
    \end{subfigure}
    \caption{\textbf{Plate With Hole. 14 displacement and 16 temperature sensors configuration. Localized thermal field.} Identified Young's moduli and temperature distributions when the thermal load is considered using different approaches during \acrshort{SI}. Peak temperatures are noted in the figures. }
\label{f:plate_local_16s}
\end{figure}

\begin{figure}[!t]
\centering
\begin{minipage}[t]{0.49\linewidth}
\centering
\includegraphics[width=\linewidth]{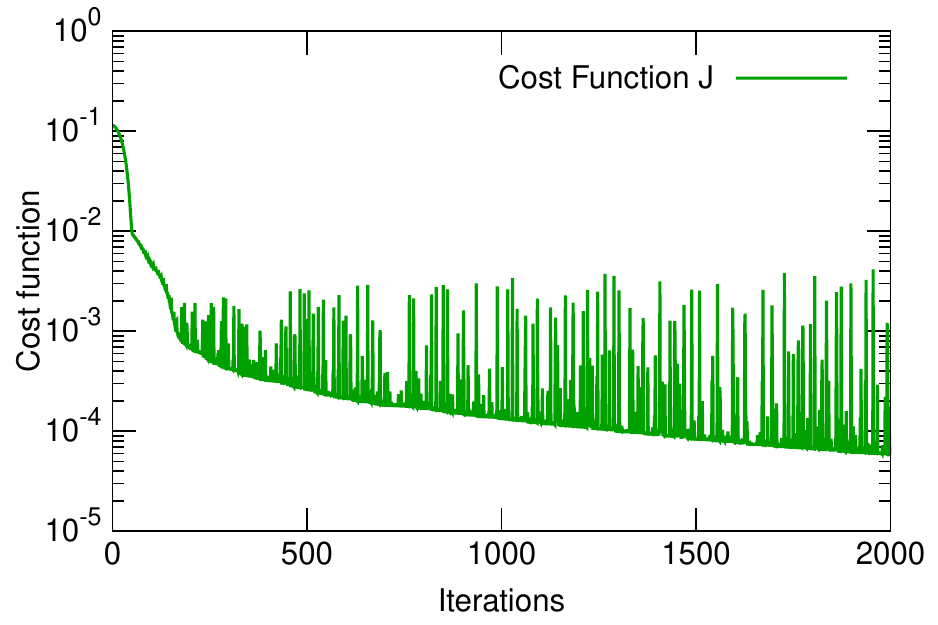}
\subcaption{When the temperature field is \textit{interpolated}.}
\label{f:plate_local_16s_conv_a}
\end{minipage}
\hfill
\begin{minipage}[t]{0.49\linewidth}
\centering
\includegraphics[width=\linewidth]{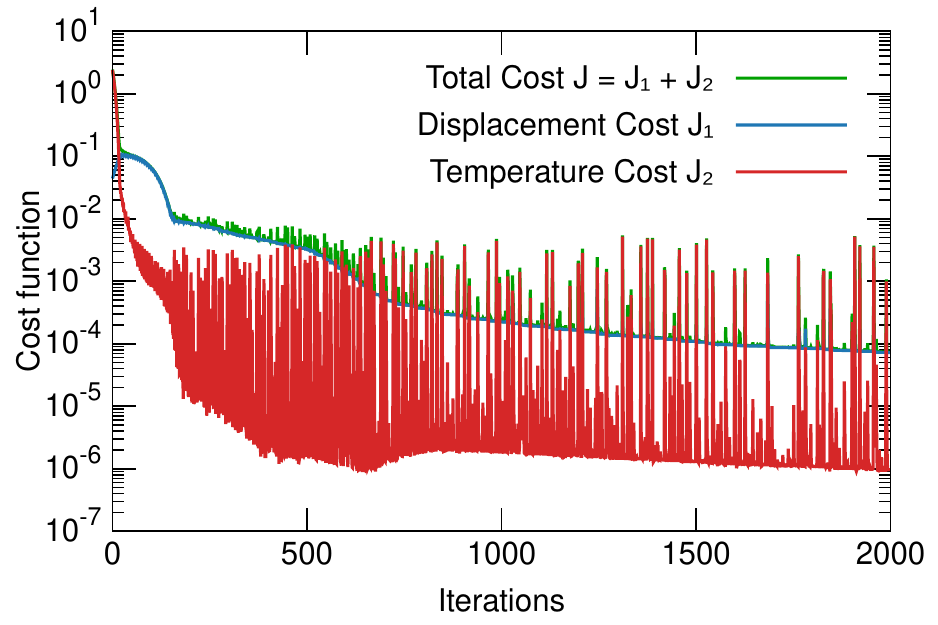}
\subcaption{When the temperature field is \textit{identified}: \textit{Monolithic} approach.}
\label{f:plate_local_16s_conv_b}
\end{minipage}
\begin{minipage}[t]{0.49\linewidth}
\centering
\includegraphics[width=\linewidth]{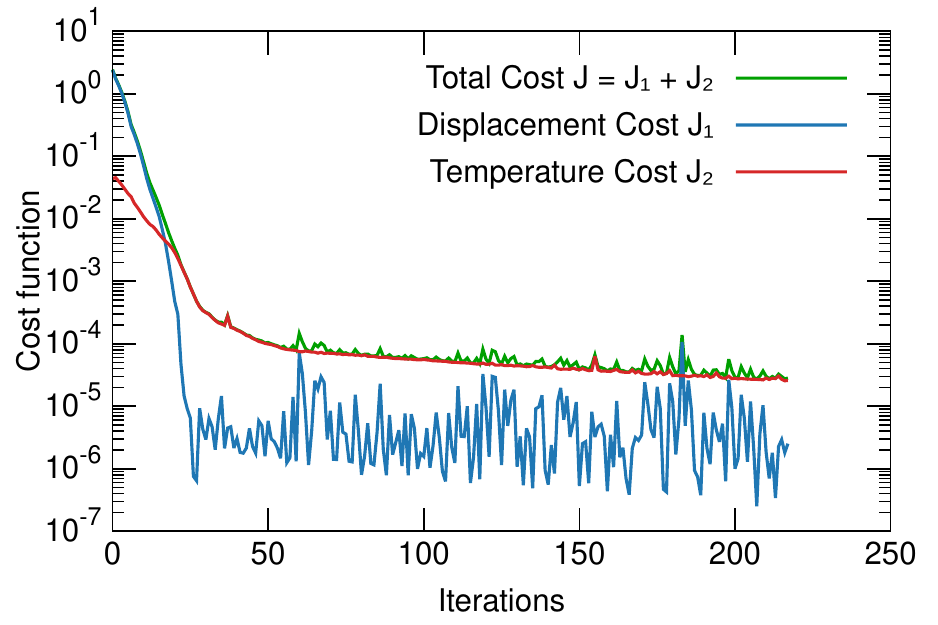}
\subcaption{When the temperature field is \textit{identified}: \textit{Partitioned} approach. Coupling iteration: 217, Overall iteration: 3992.}
\label{f:plate_local_16s_conv_c}
\end{minipage}
\caption{\textbf{Plate With Hole. 14 displacement and 16 temperature sensors configuration. Localized thermal field.}
Convergence plots when the thermal load is considered using different approaches during \acrshort{SI}.}
\label{f:plate_local_16s_conv}
\end{figure}

On the other hand, the partitioned case was set up as follows: the main \acrshort{SI} problem is reformulated into a Gauss-Seidel-type fixed-point iteration loop consisting of two sub-optimization problems. The cost function of the sub-optimization `A'  is similar to the composite function in the monolithic case (since the thermal distribution affects both $J_D$ and $J_T$) and controls the nodal temperatures. The cost function of the sub-optimization `B' is associated with the accumulated displacement sensors error and controls the elemental Young's modulus. The `maximum measured value' sensor normalization was applied separately to each sub-optimization. The relaxation factor was set to $\beta=1$, i.e., no relaxation was applied.
Per coupling iteration, the two sub-optimizations were solved inexactly by defining loose convergence criteria of $20$\% reduction in the cost function or a maximum of 10 iterations. One solve of each sub-optimization constituted one coupling iteration. The convergence criteria for the outer coupling loop were set to a 6-magnitude reduction in the composite cost function, i.e., $J = (J_D + J_T) \leq (1\cdot 10^{-6} * J_0)$ or approximately $4000$ optimization iterations when summing the iterations of the two sub-optimizations. 
In the monolithic approach, each optimization iteration updated both design variable fields ($\mathbf{E}$ and $\boldsymbol{\Delta \mathbf{T}}$). In contrast, in the partitioned approach, each sub-optimization iteration updates only one design variable field, either $\mathbf{E}$ or $\boldsymbol{\Delta \mathbf{T}}$.
To enable a fair comparison between the two approaches, the maximum number of optimization iterations for the partitioned approach was therefore set to approximately twice that of the monolithic case, so that each design variable field undergoes a comparable number of updates.
The term approximately is used because coupling iterations were not interrupted once started. As a result, the total number of optimization iterations may differ slightly from the nominal value, with a tolerance of $\pm10$ iterations. All other aspects of the optimization algorithm, including step size and update strategy, were kept identical across all cases.

The results for Scenario 3 (temperature interpolation) and Scenario 4 (temperature identification) are presented together for easy comparison. They are grouped by the number of temperature sensors and the target temperature distribution, and depicted in the following order: i) 6 temperature sensors, Linearly varying thermal field; ii) 6 temperature sensors, Localized thermal field;
iii) 16 temperature sensors, Linearly varying thermal field; iv) 16 temperature sensors, Localized thermal field.


\textbf{Case 1: 6 temperature sensors, linearly varying thermal field}

Figure \ref{f:plate_linear_6s_a} shows the identified Young’s modulus and temperature fields for Scenario 3 (kNN‑interpolated temperature) with 6 temperature sensors and a linearly varying thermal field. Relative to Scenario 2 (constant temperature field; Fig. \ref{f:plate_w_const_temp_a}), Young’s modulus identification improves markedly: most false positives are removed, and damage localization is clear. Because no temperature sensors are placed on the left/right edges, the interpolated temperature range is $[13.6,26.5]$\degree C. The convergence history (Fig. \ref{f:plate_linear_6s_conv_a}) shows a larger reduction than in Scenario 2 (Fig. \ref{f:plate_w_const_temp_a}), with the cost reaching $\mathcal{O}(10^{-5})$.

Figures \ref{f:plate_linear_6s_b} and \ref{f:plate_linear_6s_c} present the identified fields for Scenario 4 (temperature identified) for the monolithic and partitioned approaches, respectively. Both approaches further improve damage localization compared with the interpolation case, eliminating all false-positive damage detections. The identified temperature ranges are $[11.5,29.6]$ \degree C (monolithic) and $[10.0, 30.1]$ \degree C (partitioned). The corresponding convergence plots (Figs. \ref{f:plate_linear_6s_conv_b}, \ref{f:plate_linear_6s_conv_c}) differ because the partitioned case is plotted over \textit{outer coupling} iterations while the monolithic case is plotted over \textit{optimization} iterations. Despite visually similar fields, the cost function components show a contrast: the displacement cost term dominates in the monolithic case, whereas the temperature cost term dominates in the partitioned case. In both, the composite cost is $\mathcal{O}(10^{-5})$. Note that $J_D$ and $J_T$ differ by several orders of magnitude, so the total cost $J$ may be overlapped by the dominant sub‑cost function, making it difficult to see in these plots.


\textbf{Case 2: 6 temperature sensors, localized thermal field}

Figure \ref{f:plate_local_6s_a} shows the identified Young’s modulus and temperature fields for Scenario 3 (kNN interpolation) with 6 temperature sensors and a localized thermal field. As in Case 1, Young’s modulus improves over Scenario 2 (constant temperature field; Fig. \ref{f:plate_w_const_temp_b}); the main damage regions are detected, though several false positives remain. With two sensors near the plate center, the peak temperature is close to the target $30$ \degree C, and the interpolated range is $[10.5,29.2]$ \degree C. The convergence plot (Fig. \ref{f:plate_local_6s_conv_a}) shows a larger reduction than Scenario 2 (Fig. \ref{f:plate_w_const_temp_b}), with the cost of order $\mathcal{O}(10^{-4})$.

Figures \ref{f:plate_local_6s_b} and \ref{f:plate_local_6s_c} show the identified fields for Scenario 4 (temperature identified) for the monolithic and partitioned approaches, respectively. Both approaches notably enhance the Young’s modulus field and strongly attenuate false damage identification observed in Scenario 3. The identified temperatures span $[6.8,33.6]$ \degree C (monolithic) and $[7.4,34.2]$ \degree C (partitioned), which are wider than those in Scenario 3. Convergence plots (Figs. \ref{f:plate_local_6s_conv_b}, \ref{f:plate_local_6s_conv_c}) again show contrasting dominance: displacement error in the monolithic approach and temperature error in the partitioned approach, with composite costs in the order $\mathcal{O}(10^{-4}\text{–}10^{-5})$.


\textbf{Case 3: 16 temperature sensors, linearly varying thermal field}

Figure \ref{f:plate_linear_16s_a} presents Scenario 3 (kNN interpolation) with 16 temperature sensors and a linearly varying thermal field. The identified Young’s modulus and temperature fields are very close to the targets, with only minor artificial damage near the hole. The improvement is attributed to better spatial sampling - especially sensors on the left and right edges that capture the bounds of the linear field. The interpolated temperature range is $[10.0,30.0]$ \degree C. The convergence history (Fig. \ref{f:plate_linear_16s_conv_a}) shows several orders‑of‑magnitude reduction in the cost, approaching $10^{-6}$. Compared with the corresponding thermal field in the 6‑sensor case (Figs. \ref{f:plate_linear_6s_a}, \ref{f:plate_linear_6s_conv_a}), the distributions are improved, and the final cost is slightly lower.

Figures \ref{f:plate_linear_16s_b} and \ref{f:plate_linear_16s_c} show the identified fields for Scenario 4 (temperature identified) for the monolithic and partitioned approaches, respectively. Both approaches further reduce the spurious damage near the hole observed in Scenario 3. The identified temperature ranges match the target bounds: $[10.0,30.0]$ \degree C. Convergence plots (Figs. \ref{f:plate_linear_16s_conv_b}, \ref{f:plate_linear_16s_conv_c}) again show displacement- dominant cost in monolithic case and temperature‑dominant cost in partitioned case, with composite costs in the order $\mathcal{O}(10^{-4}\text{–}10^{-5})$. Compared to the same thermal field with 6‑sensor counterpart case (Figs. \ref{f:plate_linear_6s_b}, \ref{f:plate_linear_6s_c}, \ref{f:plate_linear_6s_conv_b}, \ref{f:plate_linear_6s_conv_c}), Young’s modulus is similarly accurate, while the temperature field shows a noticeable improvement.

\textbf{Case 4: 16 temperature sensors, localized thermal field}

Figure \ref{f:plate_local_16s_a} shows Scenario 3 (kNN interpolation) with 16 temperature sensors and a localized thermal field. Both fields appear deficient at first glance: false damage identification is widespread, and the maximum temperature is only $15.3$ \degree C (target $30$ \degree C). Nevertheless, the Young’s modulus field is still better than Scenario 2 (constant temperature field; Fig. \ref{f:plate_w_const_temp_b}). The interpolated temperatures lie in $[10.0,15.3]$ \degree C. The convergence plot (Fig. \ref{f:plate_local_16s_conv_a}) indicates a cost of the order of $\mathcal{O}(10^{-4})$.

Figures \ref{f:plate_local_16s_b} and \ref{f:plate_local_16s_c} display the identified fields for Scenario 4 (temperature identified) for the monolithic and partitioned approaches, respectively. Although not perfect, both approaches improve upon Scenario 3 with reduced false damage identification; the monolithic approach yields a slightly better Young’s modulus distribution. The identified temperature ranges increase to $[9.3,20.9]$ \degree C (monolithic) and $[9.2,20.8]$ \degree C (partitioned). Convergence plots (Figs. \ref{f:plate_local_16s_conv_b}, \ref{f:plate_local_16s_conv_c}) again show displacement‑dominant cost for the monolithic case and temperature‑dominant for the partitioned case, with composite costs of the order of $\mathcal{O}(10^{-4})$.

Comparing these results with the corresponding 6-sensor case results (Figs. \ref{f:plate_local_6s}, \ref{f:plate_local_6s_conv}), the identified distributions are substandard. This is directly related to the positions of the 16 temperature sensors. From Fig. \ref{f:plate_temp_sensors}, it can be seen that for the 6-sensor configuration, the positioning of the sensors is such that they capture the main features (peaks) of the localized thermal field. On the other hand, even with 16 sensors, they are unable to capture the localized behavior of the thermal field, resulting in subpar performance compared to the 6-sensor configuration.
This observation is contrary to the natural intuition that adding more sensors improves \acrshort{SI}. However, both the number and location of the sensors are important factors influencing the accuracy of the system identification.


\textbf{Quantitative Comparison}

\begin{table}[!b]
\caption{\textbf{Plate With Hole.} Identification errors in Young’s modulus distributions for scenarios \textbf{ (1 and 2)} where thermal effects are ignored or accounted for using a constant temperature field. 
Percentage changes are computed relative to the errors at optimization start. }
\label{tab:Table_2}
\begin{tabularx}{\textwidth}{Xccc}
\hline \hline
\textbf{\begin{tabular}[c]{@{}c@{}} Approach for \\ system identification    \end{tabular}} 
& \textbf{\begin{tabular}[c]{@{}c@{}} Type of \\ thermal field \end{tabular}}  & \textbf{\begin{tabular}[c]{@{}c@{}} \boldsymbol{$\epsilon_{L_2}$}   {[}-{]} \end{tabular}}  &  \textbf{\begin{tabular}[c]{@{}c@{}} \boldsymbol{$\delta\epsilon_{L_2}$}   {[}\%{]} \end{tabular}} \\ 
\hline
\multicolumn{4}{c}{\textbf{\begin{tabular}[c]{@{}c@{}}Scenario 1: Thermal load not accounted for during \acrshort{SI}\end{tabular}}}      \\ \hline
\multirow{2}{=}{\acrshort{SI} without considering temperature}          & Linearly varying      & $4.009\cdot10^{-1}$   & $95.384$   \\ \cline{2-4}    & Localized     & $3.290\cdot10^{-1}$      & $60.363$     \\ \hline        &     &    &    \\ \hline
\multicolumn{4}{c}{\textbf{\begin{tabular}[c]{@{}c@{}}Scenario 2: Thermal load accounted as a constant  \\ temperature field during \acrshort{SI}\end{tabular}}}     \\ \hline
\multirow{2}{=}{\acrshort{SI} considering constant $20$\degree C temperature field} & Linearly varying   & $2.164\cdot10^{-1}$       & $5.482 $     \\ \cline{2-4} & Localized   & $3.041\cdot10^{-1}$    & $48.194$     \\ \hline \hline
\end{tabularx}
\end{table}

Beyond visual inspection, a quantitative comparison based on the relative discrete \acrshort{$L_2$} error ($\epsilon_{L_2}$) and its percentage change ($\delta \epsilon_{L_2}$) has also been conducted. 
At the start of the system identification, the relative discrete \acrshort{$L_2$} error in Young's modulus distribution was $2.052\cdot10^{-1}$, in the temperature field for the linearly varying thermal field was $2.859\cdot10^{-1}$, while for the localized thermal field was $4.973\cdot10^{-1}$. These errors are considered as the reference errors to evaluate the performance of different approaches and also serve as the baseline to compute percentage change in the relative discrete \acrshort{$L_2$} errors ($\delta \epsilon_{L_2}$) for each case.

Table \ref{tab:Table_2} tabulates the $\epsilon_{L_2}$ and $\delta\epsilon_{L_2}$ for the Young’s modulus distributions for Scenarios 1 and 2, where thermal effects were unaccounted and where they were accounted for using a constant $20$ \degree C temperature field. 
The superscripts `$E$' and `$T$' in $\epsilon_{L_2}$ and $\delta\epsilon_{L_2}$ refer to errors in the Young's modulus and temperature fields.

From Scenario 1 errors, it can be observed that failing to account for thermal effects during \acrshort{SI} leads to catastrophic results, and the errors are actually higher than at the start of the optimization ($\delta \epsilon_{L_2}^{E}$ are positive and high). This is in line with the visual findings in Fig. \ref{f:plate_wo_temp}.

In Scenario 2, where a simplistic (constant thermal field) approach to account for thermal effects was applied, the errors still increased compared to the initial values; however, they were better than those in Scenario 1 ($\delta \epsilon_{L_2}^{E}$ are still positive but lower magnitude), indicating that the constant $20$ \degree C field provided some improvement, which was also confirmed visually in Fig. \ref{f:plate_w_const_temp}.

Table \ref{tab:Table_3} tabulates the $\epsilon_{L_2}$ and $\delta\epsilon_{L_2}$ for the Young’s modulus and temperature distributions, where the temperature sensor information was incorporated to account for thermal effects, i.e., Scenario 3 (via temperature interpolation) and Scenario 4 (via temperature identification). The table is grouped based on the number of temperature sensors and the type of thermal field. 

For case 1 (6-sensors and linearly varying thermal field), the interpolation case improved the Young's modulus identification ($\delta \epsilon_{L_2}^{E}=-48.531$\%), which can be attributed to the better approximation of the thermal field ($\delta \epsilon_{L_2}^{T}=-67.66$\%). The monolithic and partitioned approaches in the identification case demonstrated marginally better accuracy than the interpolation case. These values are consistent with the visual examination in Fig. \ref{f:plate_linear_6s}.

For case 2 (6-sensors and localized thermal field), the interpolation case yielded subpar results, with an increase in Young's modulus error but a significant reduction in temperature error. On the other hand, both the monolithic and partitioned approaches in the identification case outperformed the interpolation case with substantially better identification of Young's modulus distribution and slightly better identification of the thermal field. These values corroborate the observations in Fig. \ref{f:plate_local_6s}.

For case 3 (16-sensors and linearly varying thermal field), the interpolation and the identification cases yielded accurate reconstructions (as confirmed earlier by the visual examination in Fig. \ref{f:plate_linear_16s}), with the interpolation case performing marginally better for Young's modulus distribution and the identification case performing marginally better for the temperature distribution. 

For case 4 (16-sensors and localized thermal field), both the interpolation and identification cases yielded better results than the constant thermal field case (Scenario 2), but were suboptimal compared to the corresponding 6-sensors configuration (case 2 above). As mentioned earlier, this was due to the positioning of the sensors, which did not capture the localized trend in the thermal field effectively. Interestingly, the identification approach outperformed the interpolation approach. This trend aligns with the observation in Fig. \ref{f:plate_local_16s}, where better maximum temperatures were seen in the thermal fields, and a reduction in the false damage identification was observed.
The reductions in $\delta\epsilon_{L_2}^{E},\delta\epsilon_{L_2}^{T}$ for the identification approach were almost double that of the interpolation approach. 

\begin{table}[!t]
\caption{\textbf{Plate With Hole.} Identification Errors in Young's modulus (\textbf{E}) and temperature ($\boldsymbol{\Delta \mathbf{T}}$) fields for Scenarios (\textbf{3 and 4}) \textbf{incorporating temperature sensor information} to account for thermal effects. Percentage changes are computed relative to the errors at optimization start. }
\label{tab:Table_3}
\begin{tabularx}{\textwidth}{Xcccc}
\hline \hline
\multirow{2}{=}{\parbox{\linewidth}{\centering
\textbf{Approach accounting}\\
\textbf{for thermal effects}}}                                        & \multicolumn{2}{c}{\textbf{\begin{tabular}[c]{@{}c@{}} \boldsymbol{$\epsilon_{L_2}$}   {[}-{]}\end{tabular}}} & \multicolumn{2}{c}{\textbf{\begin{tabular}[c]{@{}c@{}} \boldsymbol{$\delta\epsilon_{L_2}$}   {[}\%{]}\end{tabular}}} \\ \cline{2-5}   & \textbf{E}                                      & $\boldsymbol{\Delta \mathbf{T}}$                                     & \textbf{E}                                      & $\boldsymbol{\Delta \mathbf{T}}$                                    \\ \hline
\multicolumn{5}{c}{\textbf{Case 1: 6 temperature sensors, Linearly varying thermal field}}                                                              \\ \hline
temperature interpolation  &$ 1.056\cdot10^{-1}$                                       & $9.248\cdot10^{-2}$                                      & $-48.531$                                         & $-67.660$                                        \\ \hline
temperature identification & \multicolumn{1}{l}{}                            & \multicolumn{1}{l}{}                           & \multicolumn{1}{l}{}                            & \multicolumn{1}{l}{}                           \\
\multicolumn{1}{c}{Monolithic approach}                                       & $1.028\cdot10^{-1}$                                       & $3.939\cdot10^{-2}$                                      & $-49.902$                                         & $-86.223$                                        \\ 
\multicolumn{1}{c}{Partitioned approach}                                      & $1.027\cdot10^{-1}$                                       & $3.287\cdot10^{-2}$                                      & $-49.941$                                         & $-88.504$                                        \\ \hline
   &    &    &     &      \\ \hline
   \multicolumn{5}{c}{\textbf{Case 2: 6 temperature sensors, Localized thermal field}}     \\ \hline
temperature interpolation   & $2.064\cdot10^{-1}$                                       & $2.055\cdot10^{-1} $                                     & $0.599$                                           & $-58.667$                                        \\ \hline
temperature identification & \multicolumn{1}{l}{}                            & \multicolumn{1}{l}{}                           & \multicolumn{1}{l}{}                            & \multicolumn{1}{l}{}                           \\
\multicolumn{1}{c}{Monolithic approach}                                       & $1.041\cdot10^{-1} $                                      & $1.811\cdot10^{-1}$                                      & $-49.280$                                         & $-63.575$                                        \\ 
\multicolumn{1}{c}{Partitioned approach}                                      & $9.635\cdot10^{-2}$                                       & $1.879\cdot10^{-1}$                                      & $-53.043$                                         & $-62.222$                                        \\ \hline &                                                 &                                                &                                                 &                                                \\ \hline
\multicolumn{5}{c}{\textbf{Case 3: 16 temperature sensors, Linearly varying thermal field}}    \\ \hline
temperature interpolation   & $9.631\cdot10^{-2}$                                       & $3.766\cdot10^{-2}$                                      & $-53.064$                                         & $-86.830$                                        \\ \hline
temperature identification & \multicolumn{1}{l}{}                            & \multicolumn{1}{l}{}                           & \multicolumn{1}{l}{}                            & \multicolumn{1}{l}{}                           \\
\multicolumn{1}{c}{Monolithic approach}                                       & $1.071\cdot10^{-1}$                                       & $7.614\cdot10^{-3}$                                      & $-47.811$                                         & $-97.337$                                        \\ 
\multicolumn{1}{c}{Partitioned approach}                                      & $9.956\cdot10^{-2}$                                       & $8.843\cdot10^{-3}$                                      & $-51.478$                                         & $-96.907$                                        \\ \hline   &                                                 &                                                &                                                 &                                                \\ \hline

\multicolumn{5}{c}{\textbf{Case 4: 16 temperature sensors, Localized thermal field}}    \\ \hline
temperature interpolation  & $1.802\cdot10^{-1} $                                      & $3.500\cdot10^{-1}$                                      & $-12.200$                                         & $-29.620$                                        \\ \hline
temperature identification & \multicolumn{1}{l}{}                            & \multicolumn{1}{l}{}                           & \multicolumn{1}{l}{}                            & \multicolumn{1}{l}{}                           \\
\multicolumn{1}{c}{Monolithic approach}                                       & $1.484\cdot10^{-1}$                                       & $2.054\cdot10^{-1}$                                      & $-27.669$                                         & $-58.692$                                        \\ 
\multicolumn{1}{c}{Partitioned approach}                                      & $1.569\cdot10^{-1} $                                      & $2.077\cdot10^{-1}$                                      & $-23.532$                                         & $-58.241$                                        \\ \hline  \hline
\end{tabularx}
\end{table}

\subsection{Footbridge}
\label{s:footbridge}

\begin{figure}[!t]
\centering
\includegraphics[trim=0 0 0 0, clip, width=0.99\linewidth]{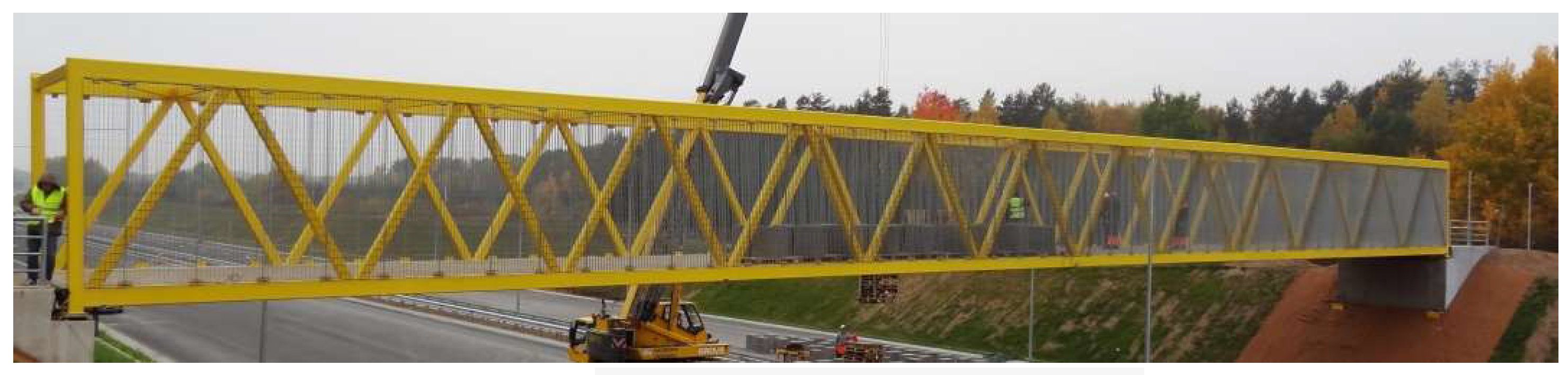}
\caption{\textbf{Footbridge.} General view of the footbridge located above the western bypass of Vilnius, near V. Maciulevičius Street (Vilnius, Lithuania). Reproduced from \citep{kilikevivcius2020influence} under the Creative Commons CC BY 4.0 license.}
\label{f:footbridge_vilinius}
\end{figure}

\begin{figure}[!t]
\centering
\includegraphics[trim=100 760 100 760, clip, width=0.95\linewidth]{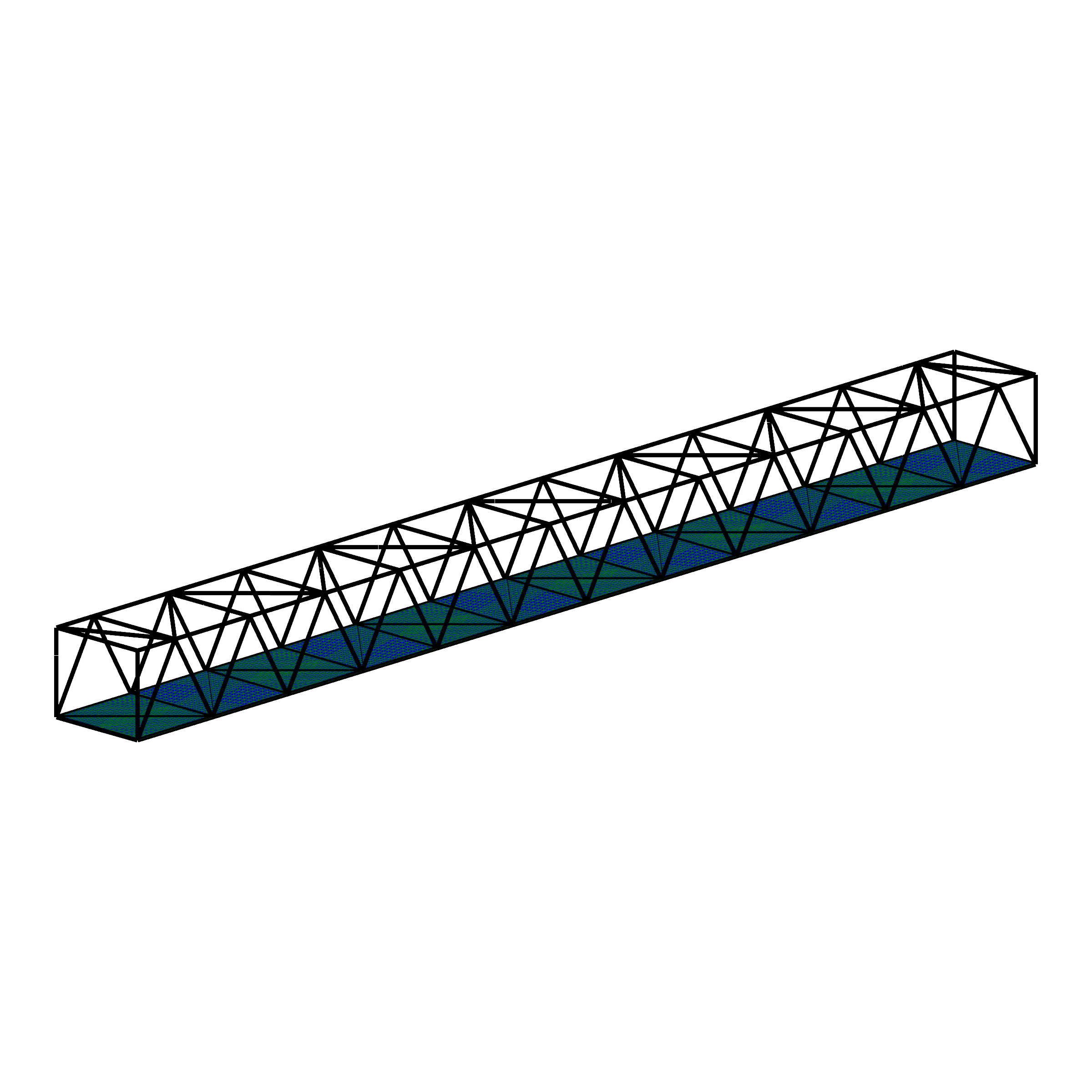}
\caption{\textbf{Footbridge.} Mesh.}
\label{f:footbridge_mesh}
\end{figure}

 \begin{figure}[!b]
\centering
\includegraphics[trim=0 530 0 760, clip, width=0.85\linewidth]{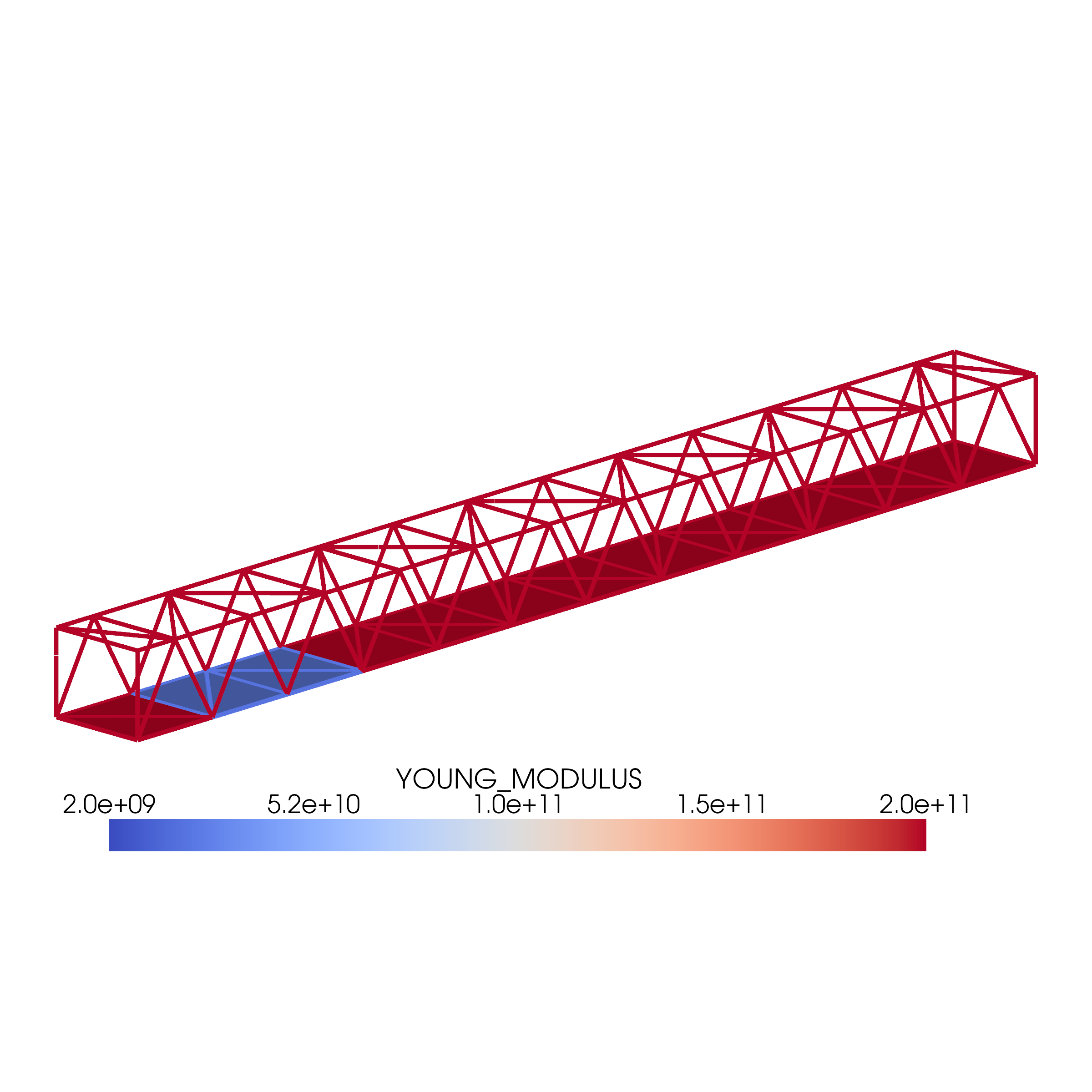}
\caption{\textbf{Footbridge.} Target Young's modulus distribution, i.e., localized damages.}
\label{f:footbridge_target_damage}
\end{figure}  

In this section, a more realistic example of a footbridge structure shown in Fig. \ref{f:footbridge_vilinius} is presented. In this work, a simplified model of the bridge is used. The bridge has dimensions: $0 \le x \le 48.3$, $0 \le y \le 4.2$, and $0 \le z \le 3.5$. All quantities are in SI units, unless explicitly mentioned. The footbridge is discretized using a fine mesh of $7680$ linear, triangular, plane stress shell elements for the deck with a thickness of $0.01$, and $2632$ two-noded linear beam elements for the rest of the structure with cross-area of $1.1424\cdot10^{-2}$, second moments of area $I_{22}=7.6\cdot10^{-5}$, $I_{33}=1.5\cdot10^{-4}$, and torsional inertia of $1.0\cdot10^{-4}$. The structure is illustrated in Fig. \ref{f:footbridge_mesh}.

The material properties of the beams and shells were: density $\rho = 7800$, Young's modulus (in pristine state) $\text{E} = 2\cdot 10^{11}$, Poisson's ratio $\nu = 0.3$, and coefficient of thermal expansion $\alpha = 1.0 \cdot 10^{-5}$ /\degree C.  The bridge is fixed on the left bottom and right bottom edges ($\mathbf{u}(x=0,z=0)=\mathbf{u}(x=48.3,z=0)= \mathbf{0}$), and a static downward uniform pressure of $P=1\cdot10^{4}$ $ \text{N/m}^2$ is applied on the bridge deck shell elements.

\begin{figure}[!t]
\centering
\subfloat[\centering]{\includegraphics[trim=0 745 0 760, clip, width=0.85\linewidth]{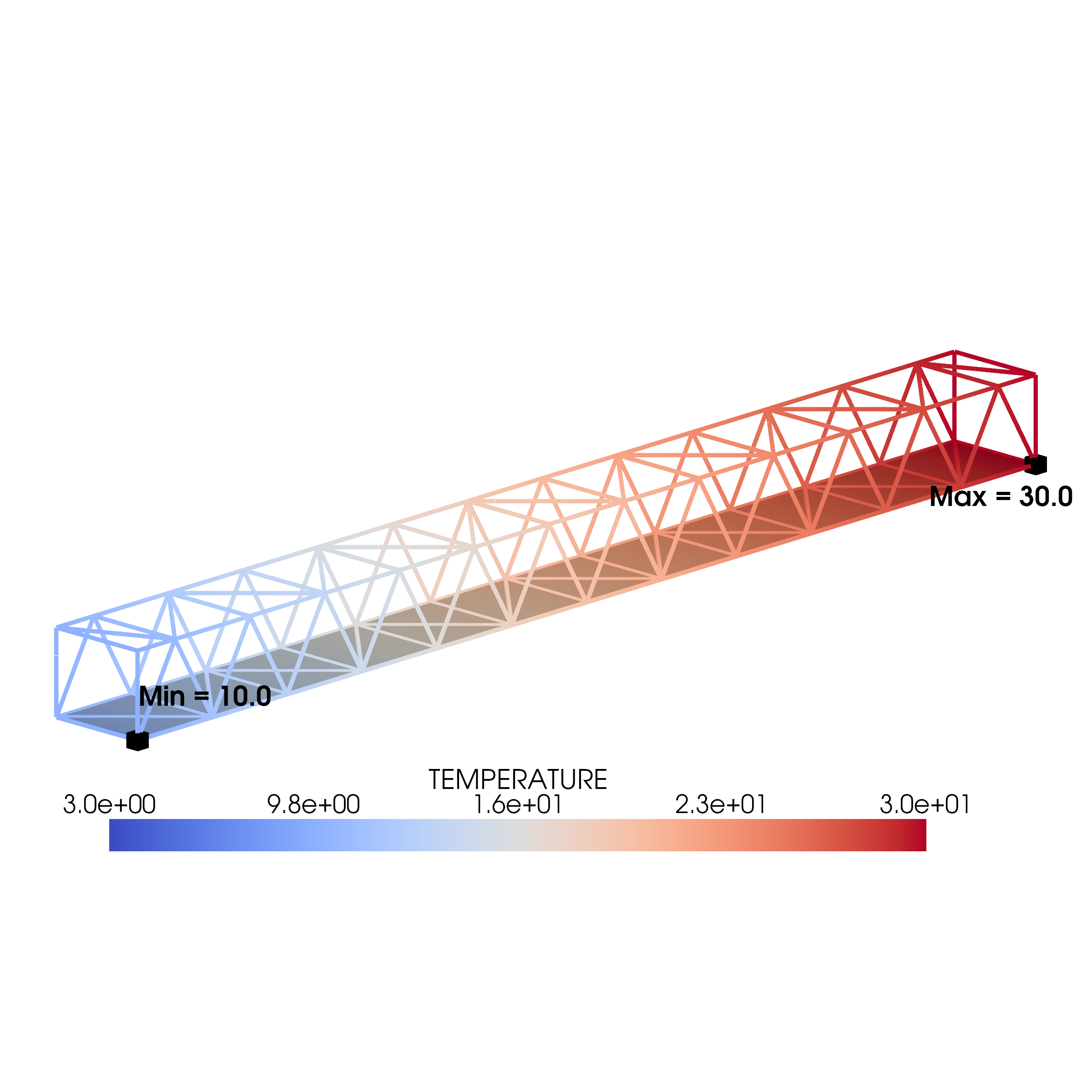}
\label{f:footbridge_target_temp_a}}
\hfill
\subfloat[\centering]{\includegraphics[trim=0 530 0 760, clip, width=0.85\linewidth]{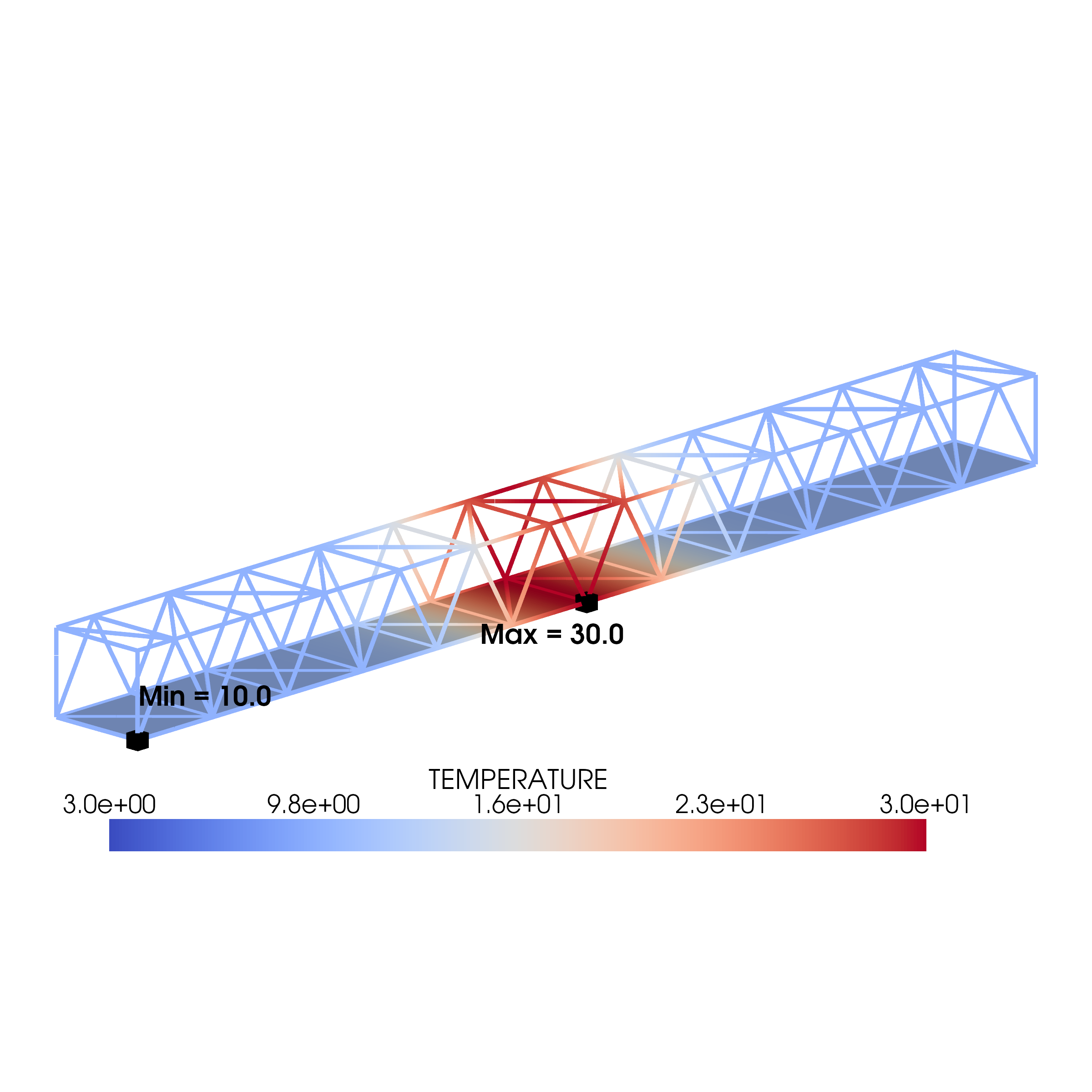}
\label{f:footbridge_target_temp_b}}
\caption{\textbf{Footbridge.} Target temperature distributions: (\textbf{a}) Linearly varying thermal field, and (\textbf{b}) Localized thermal field.}
\label{f:footbridge_target_temp}
\end{figure}

In the damaged state, the Young's modulus was prescribed to be reduced to $\text{E}=2\cdot10^{10}$ for the beam and shell elements in the locally damaged region of $4.025 \leq x \leq 12.075$, $z=0$. This `target' Young's modulus distribution is shown in Fig. \ref{f:footbridge_target_damage}.
The actual structure was subjected to a thermal load. Similar to the plate with hole example, two different types of thermal fields were investigated. The first was a linearly varying thermal field with $\Delta \text{T}(x=0)=10 $\degree C at the left edge and $\Delta \text{T}(x=48.3)=30 $\degree C at the right edge. The second was a localized Gaussian type thermal field defined as $\Delta \text{T}(\mathbf{x}) = \Delta \text{T}_{\text{min}} + (\Delta \text{T}_{\text{max}}- \Delta \text{T}_{\text{min}})e^{\frac{-(x-24.15)^{2}}{30}}$, $\Delta \text{T}_{\text{min}}=10 $\degree C, $\Delta \text{T}_{\text{max}}=30 $\degree C. These `target' temperature distributions are illustrated in Fig. \ref{f:footbridge_target_temp}. The color bar is shown in the range $[3.0, 30.0]$ \degree C for consistent visualization and better comparison with later results.

\begin{figure}[!b]
\centering
\includegraphics[trim=100 760 100 760, clip, width=0.95\linewidth]{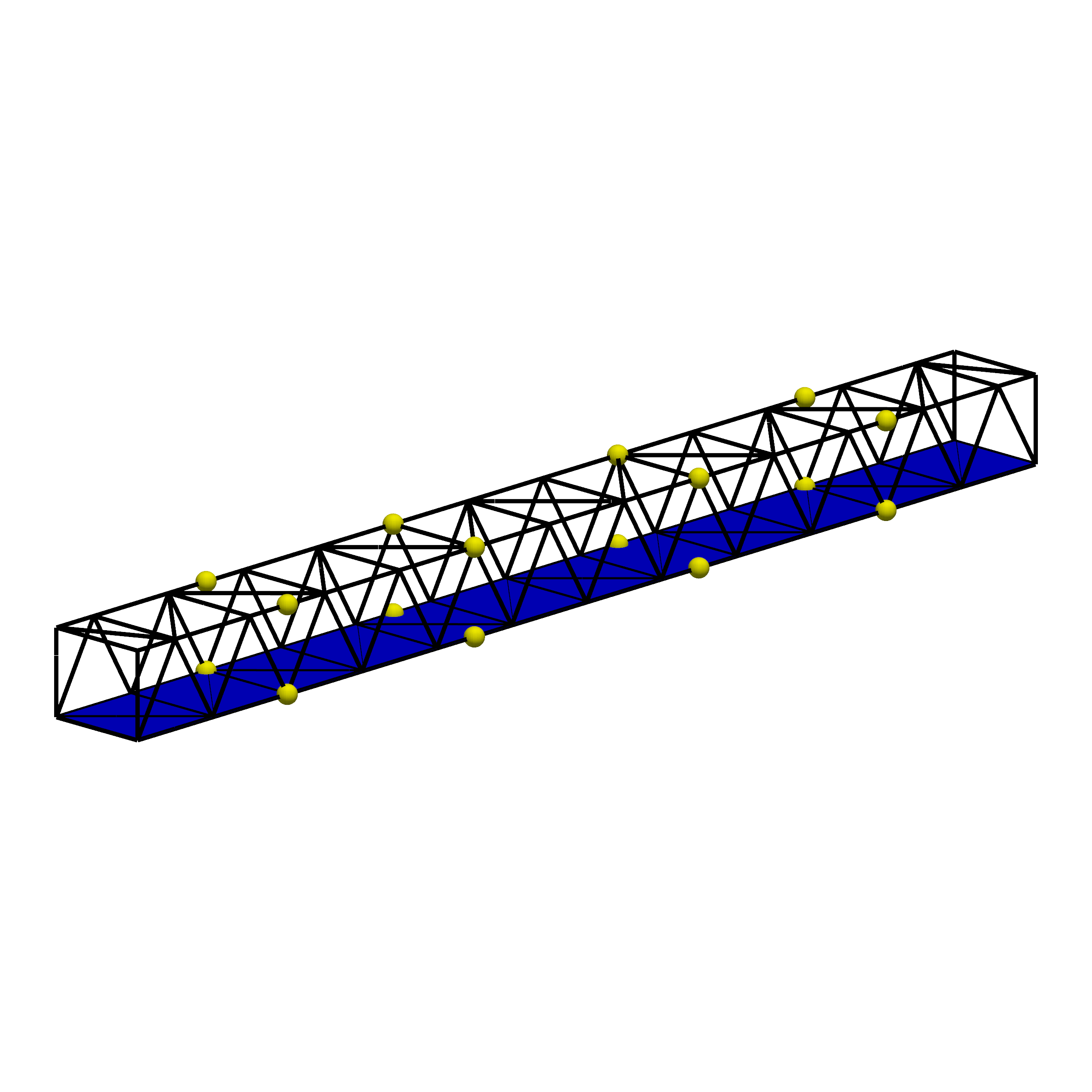}
\caption{\textbf{Footbridge.} Location of the 16 displacement sensors.}
\label{f:footbridge_sensors}
\end{figure}

\begin{figure}[!b]
\centering
\subfloat[\centering]{\includegraphics[trim=0 730 0 730, clip, width=0.85\linewidth]{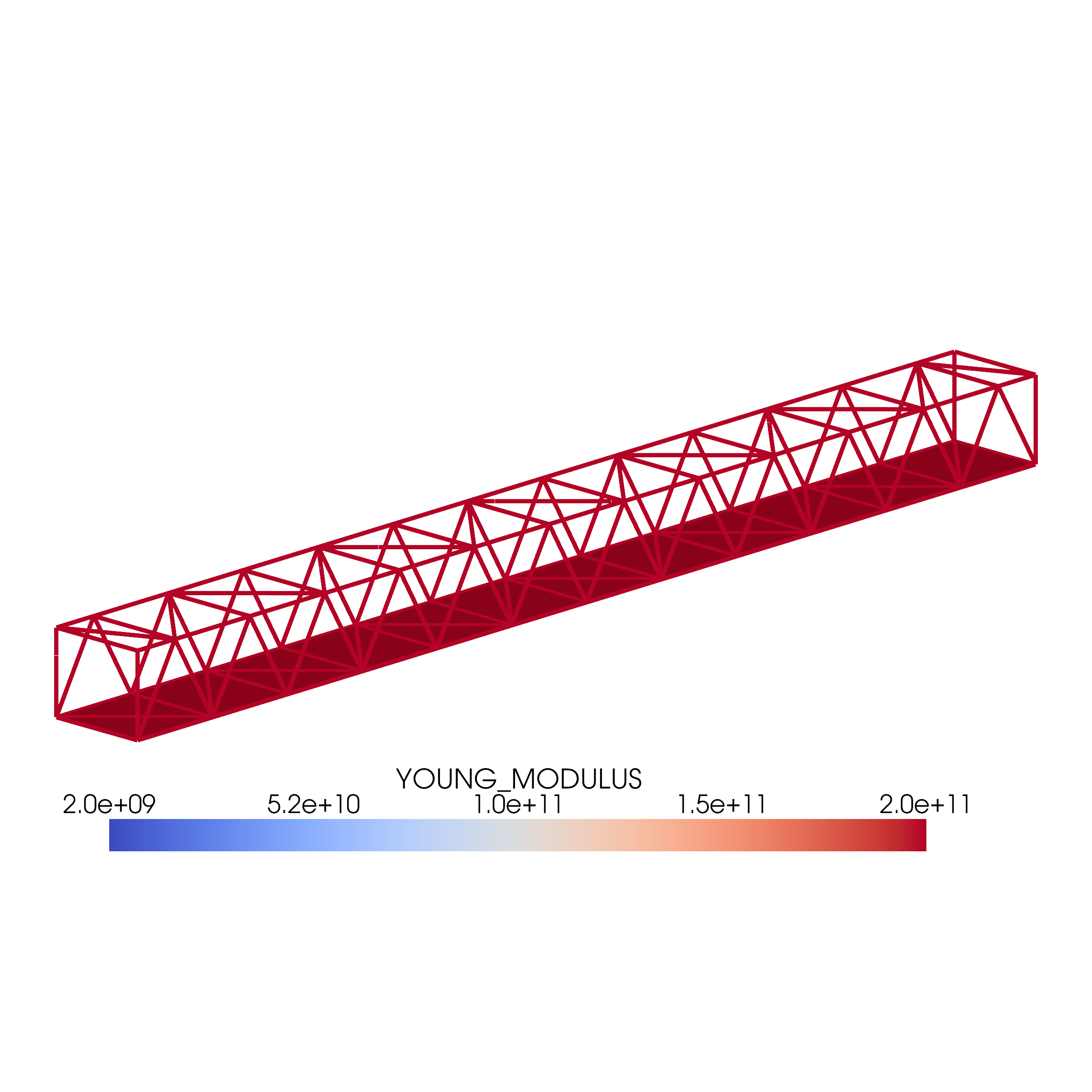}
\label{f:footbridge_wo_temp_a}}
\hfill
\subfloat[\centering]{\includegraphics[trim=0 530 0 760, clip, width=0.85\linewidth]{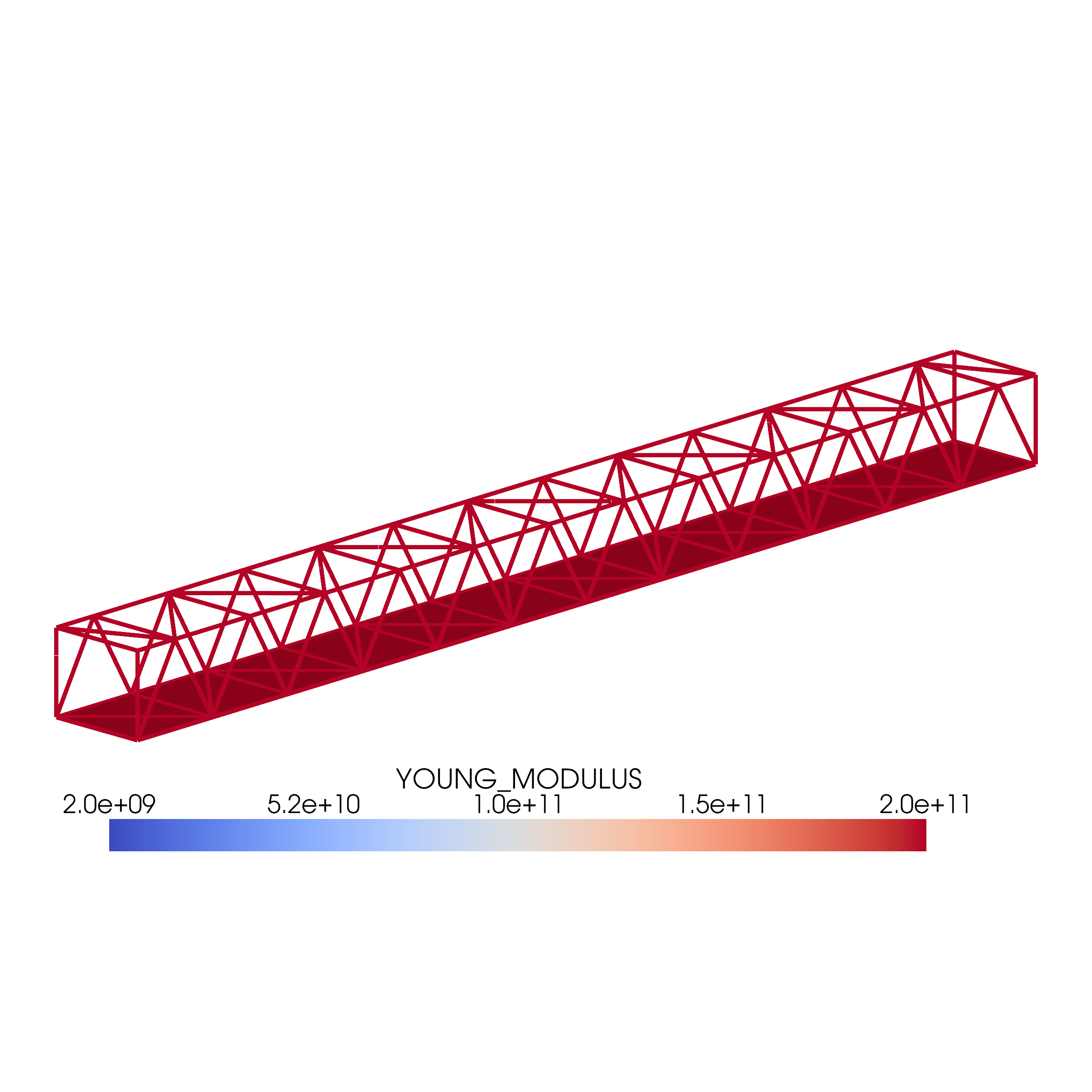}
\label{f:footbridge_wo_temp_b}}
\hfill
\subfloat[\centering Convergence plots for the linearly varying (left) and the localized (right) thermal field cases.]{\includegraphics[trim=0 0 0 0, clip, width=0.49\linewidth]{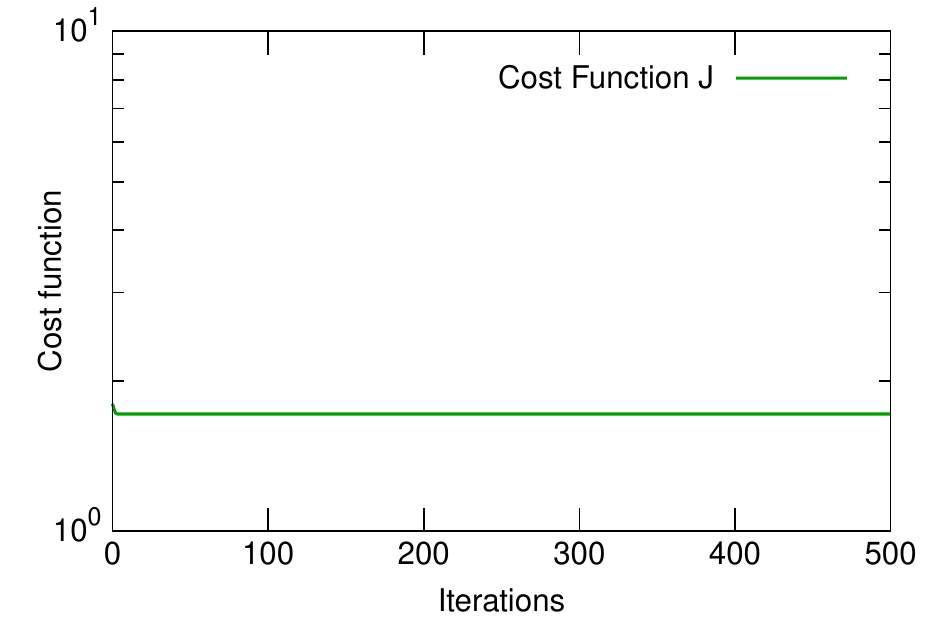}
\includegraphics[trim=0 0 0 0, clip, width=0.49\linewidth]{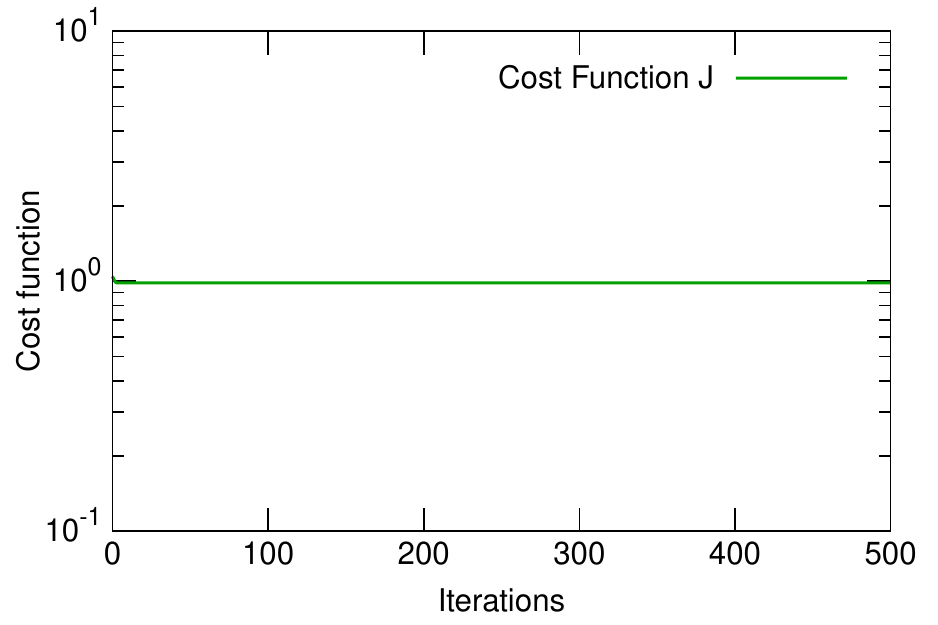}
\label{f:footbridge_wo_temp_c}}
\caption{\textbf{Footbridge.} Identified Young's modulus distributions when thermal load is not considered during \acrshort{SI}, but the actual structure is subjected to a: (\textbf{a}) Linearly varying thermal field, (\textbf{b}) Localized thermal field.}
\label{f:footbridge_wo_temp}
\end{figure}

\begin{figure}[!b]
\centering
\subfloat[\centering]{\includegraphics[trim=0 730 0 730, clip, width=0.85\linewidth]{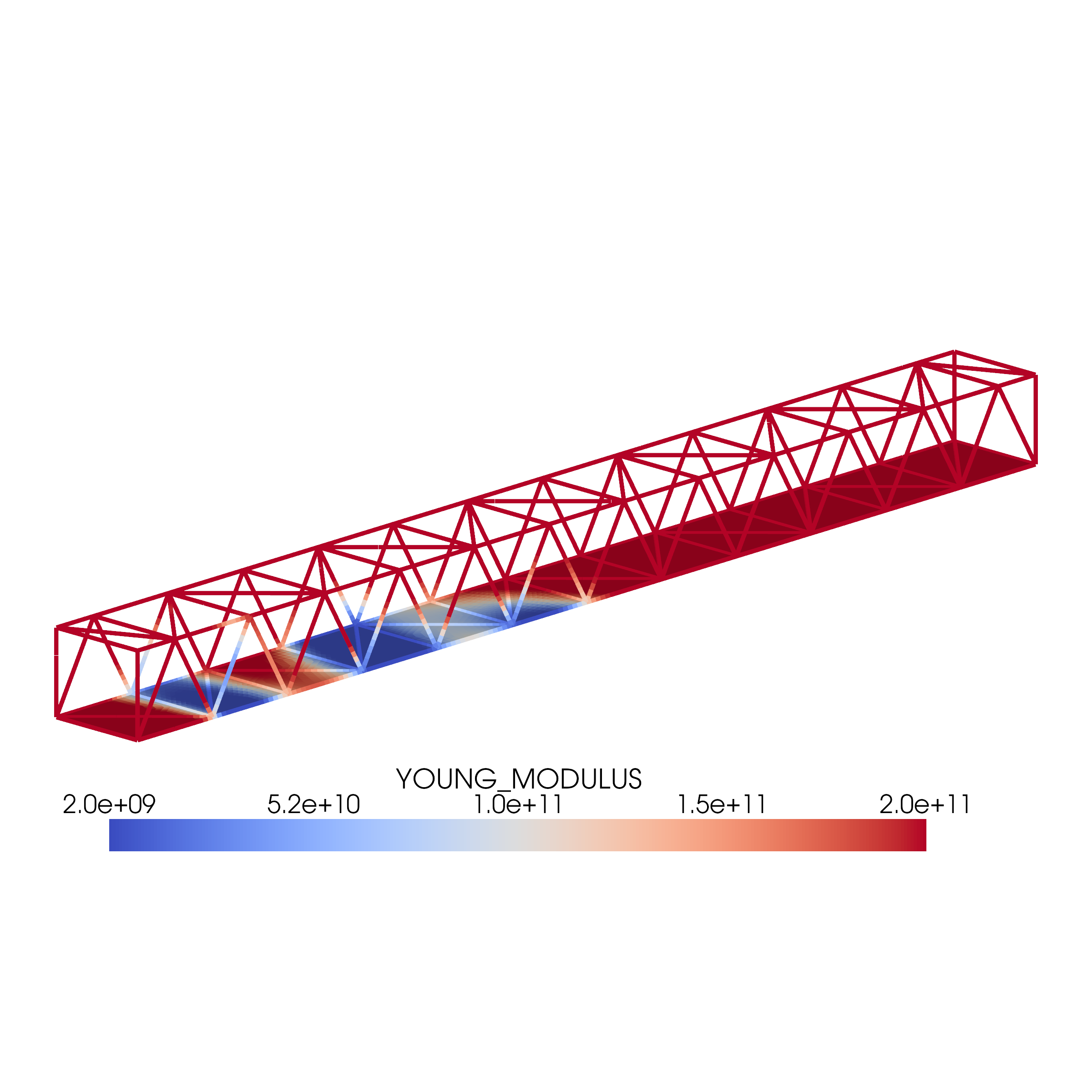}
\label{f:footbridge_w_const_temp_a}
}
\hfill
\subfloat[\centering]{\includegraphics[trim=0 530 0 760, clip, width=0.85\linewidth]{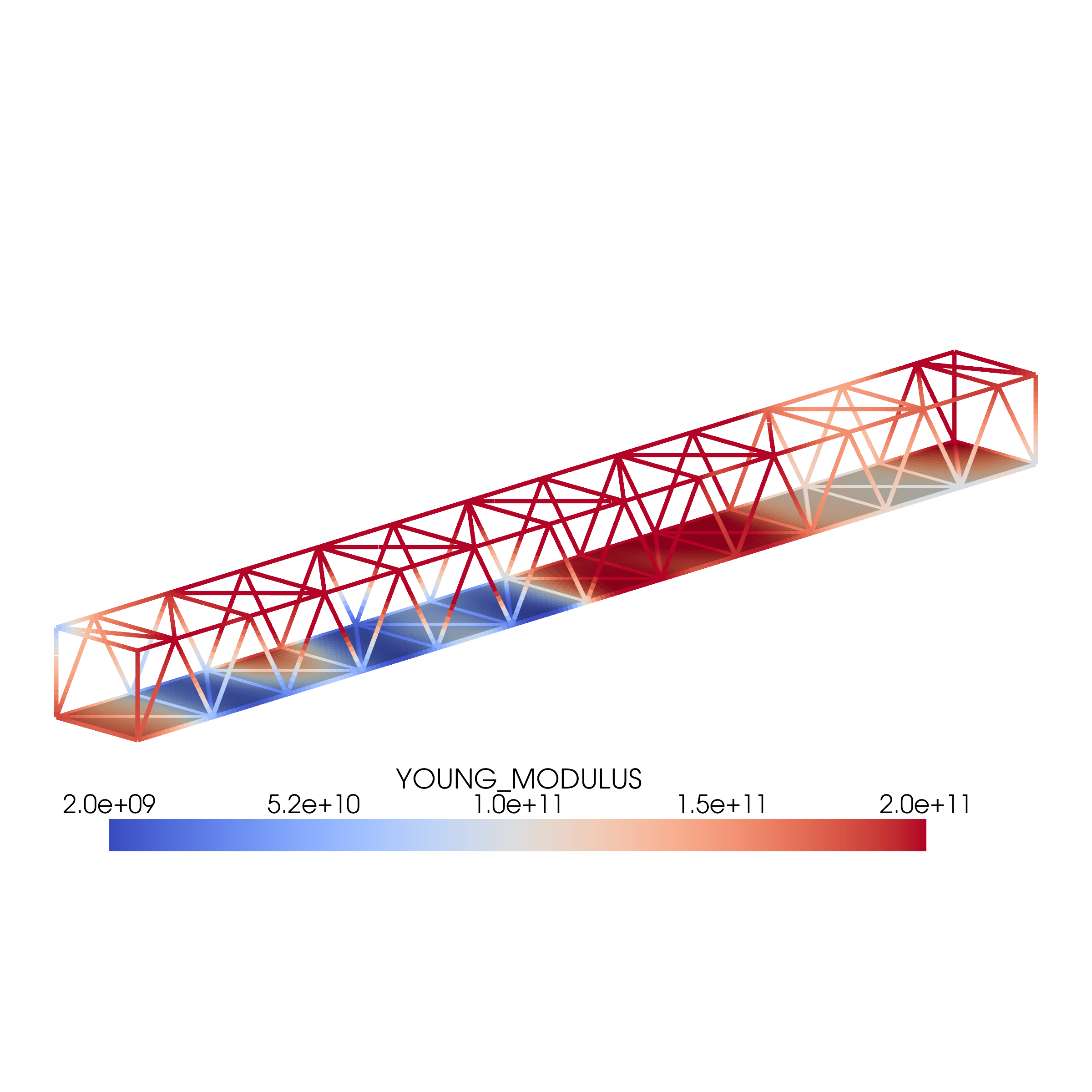}
\label{f:footbridge_w_const_temp_b}}
\hfill
\subfloat[\centering Convergence plots for the linearly varying (left) and the localized (right) thermal field cases.]{\includegraphics[trim=0 0 0 0, clip, width=0.49\linewidth]{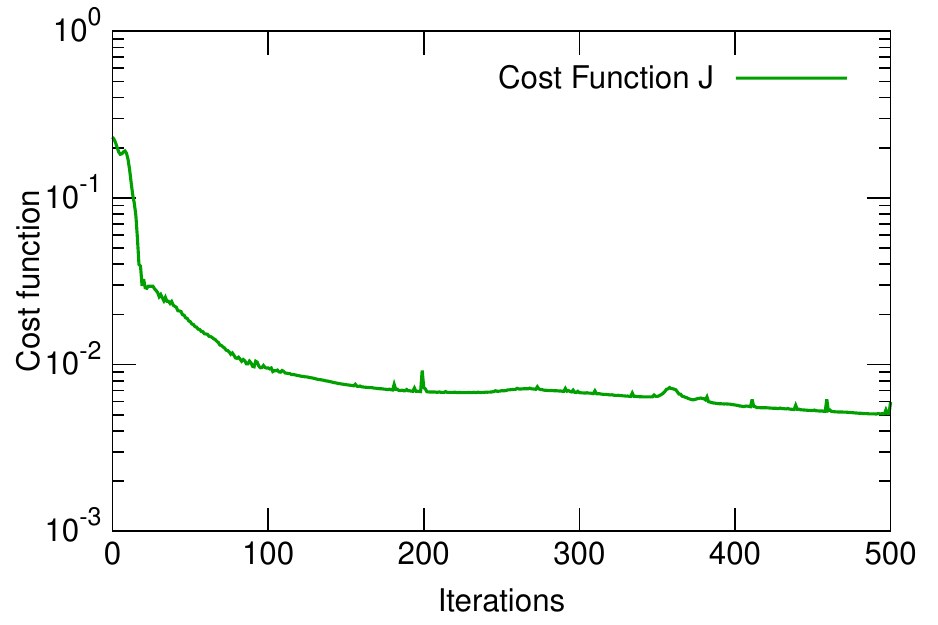}
\includegraphics[trim=0 0 0 0, clip, width=0.49\linewidth]{Figure_22c.pdf}
\label{f:footbridge_w_const_temp_c}}
\caption{\textbf{Footbridge.} Identified Young's modulus distributions when a \textit{constant} temperature distribution of $20$\degree C is considered during \acrshort{SI}, but the actual structure is subjected to a: (\textbf{a}) Linearly varying thermal field, (\textbf{b}) Localized thermal field.}
\label{f:footbridge_w_const_temp}
\end{figure}

To measure the deformations, sixteen displacement sensors were distributed over the footbridge as depicted in Fig. \ref{f:footbridge_sensors}. At the start of system identification, the footbridge was assumed to be in near-pristine condition with $\mathbf{E}(\mathbf{x}) = 1.98\cdot 10^{11}$. The `maximum measured value' sensor weighting (Eq.\eqref{eq:sensor_weight_a}) was used for normalization. 
The Nesterov accelerated gradient descent algorithm with \acrshort{QN-BB} method and a maximum step size of $5\cdot 10^{-1}$ was used for the optimization.
The convergence criteria were set to a target cost function value of $J \leq 1\cdot 10^{-7}$ or a maximum of 500 iterations. The Young's modulus control variables were bounded in the range of $[2\cdot10^{9},2\cdot10^{11}]$.
To help with the ill-conditioned \acrshort{SI}, \acrshort{VM} with a linear kernel and radius $r = 5$ was employed to regularize gradients and updates.

To compare results, the damage location and intensity, as well as the optimization and \acrshort{VM} settings, were kept the same across all cases, except in Scenario 4, where the overall maximum number of iterations was $1000$ for the partitioned case, to have approximately the same number of design variable updates as the monolithic approach.

Scenario 1, as described in Table \ref{tab:Table_1}, is the case where the actual structure is subjected to thermal loads, but thermal effects are not considered during \acrshort{SI}. 
The identified Young's modulus distribution for this scenario for the two thermal fields and their convergence plots are shown in Fig. \ref{f:footbridge_wo_temp}. 
It can be observed that for both the thermal field cases, no damage whatsoever has been identified. This is the most dangerous scenario because failing to account for thermal effects leads to complete failure in damage detection and localization. From the cost functions, it can be seen that the optimizer is unable to reduce the cost. 

Scenario 2, as described in Table \ref{tab:Table_1}, improves on Scenario 1 and accounts for thermal effects via a constant thermal field applied to the structure in the primal simulation during \acrshort{SI}. 
The identified Young's modulus distribution for this scenario for the two thermal fields and their convergence plots are shown in Fig. \ref{f:footbridge_w_const_temp}. An improvement in the results is observed, with some damage localization, but many false damage detections are also being identified. Overall, the identification is quite deficient. The cost functions are in the order $\mathcal{O}(10^{-2})$.

\begin{figure}[!t]
\centering
\includegraphics[trim=0 730 0 730, clip, width=0.99\linewidth]{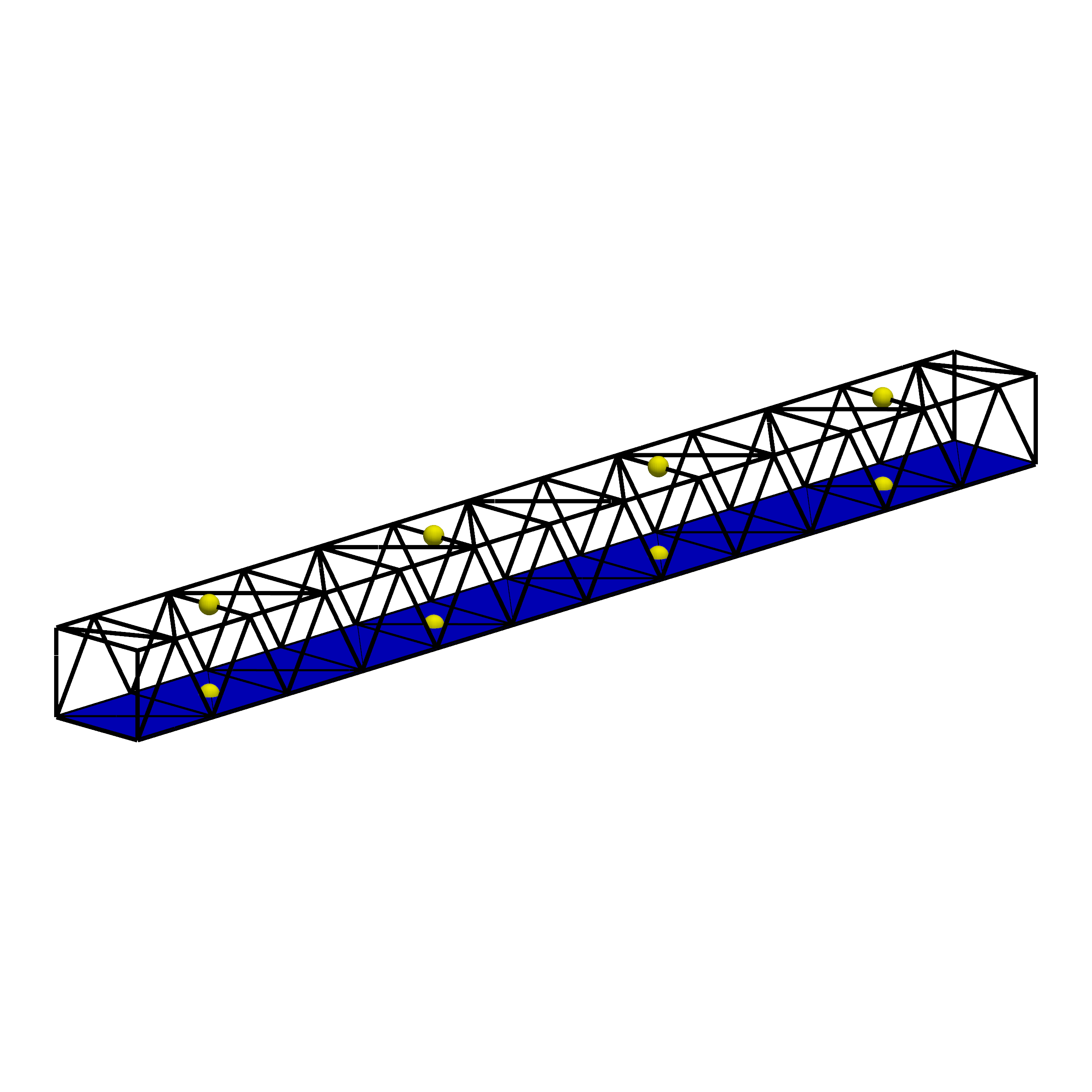}
\caption{\textbf{Footbridge.} Location of the 8 temperature sensors.}
\label{f:footbridge_temp_sensors}
\end{figure} 

Similar to the plate with hole example, adequate identification is not obtained in Scenarios 1 and 2, thus requiring a better estimation of the thermal field of the structure for meaningful damage localization. Hence, eight temperature sensors are placed on the structure. The temperature sensor locations are illustrated in Fig. \ref{f:footbridge_temp_sensors}. 

Scenario 3, as described in Table \ref{tab:Table_1}, is one in which the temperature distribution of the structure is approximated using interpolation of the temperature sensor measurements. In this example, interpolation was performed using \acrshort{kNN} with 3 nearest neighbors and inverse distance weighting. This interpolated temperature distribution is applied as a fixed thermal load on the structure in the primal simulation during \acrshort{SI}. The quantity being identified remains the Young's modulus distribution, as the \acrshort{SI} problem remains unchanged; only the manner in which thermal effects are accounted for has changed. 

Scenario 4, as described in Table~\ref{tab:Table_1}, is the proposed methodology in which the temperature distribution of the structure is inferred during \acrshort{SI}, alongside the Young's modulus distribution, using the monolithic and partitioned approaches. 
The monolithic and partitioned \acrshort{SI} cases were set up similarly to the plate with hole example. 
At the start of the optimization, the bridge had a uniform temperature of $\Delta \text{T}(\mathbf{x}) = 20$ \degree C. The bound for the nodal temperature control variable was set to the range $[-10, 40]$ \degree C. 
The `maximum measured value' sensor normalization was applied separately to each sensor type according to Eqs.~\eqref{eq:sensor_weight_a},~\eqref{eq:sensor_weight_b}.
\acrshort{VM} with a linear kernel and radius $r = 20$ was employed for the temperature field to help regularize the problem. All other settings related to the optimization algorithm, convergence criteria, etc, were kept the same.

For the partitioned approach, the `maximum measured value' sensor normalization was applied separately to each sub-optimization.
The relaxation factor was set to $\beta=1$, i.e., no relaxation was applied. 
Per coupling iteration, the two sub-optimizations were solved inexactly by defining loose convergence criteria of $30$\% reduction in the cost function or a maximum of 5 iterations. One solve of each sub-optimization constituted one coupling iteration. The convergence criteria for the outer coupling loop were set to a target composite cost function value $J = (J_D + J_T) \leq (1\cdot 10^{-7})$ or approximately $1000$ optimization iterations (summing iterations from the two sub-optimizations). The optimization algorithm and step were kept the same as all the previous scenarios.

\begin{figure}[!b]
    \centering
    \begin{minipage}[t]{\textwidth}
        \centering
        \begin{subfigure}[t]{\textwidth}
            \centering
            \includegraphics[trim=0 730 0 730, clip,width=0.9\textwidth]{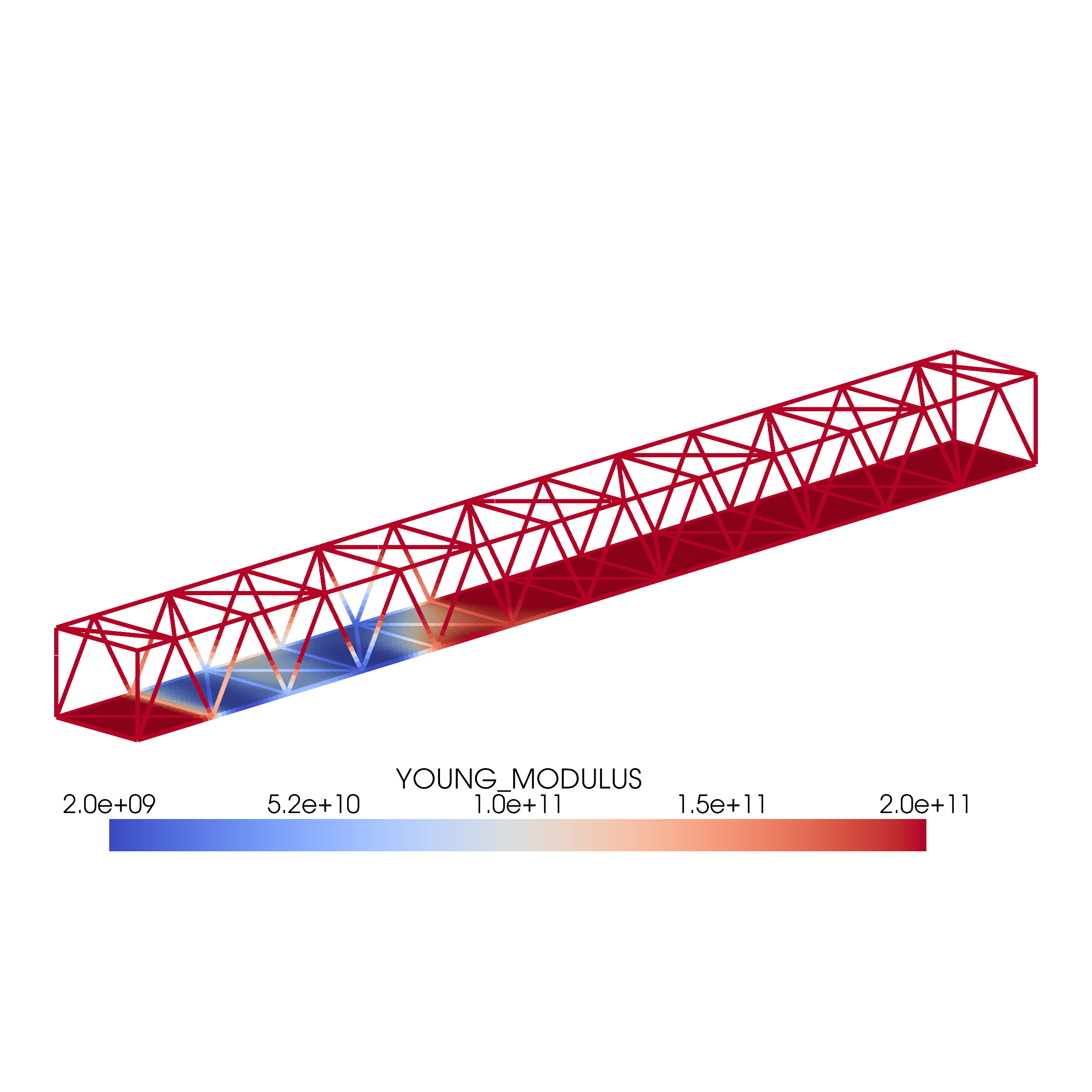}
            \caption{Identified Young's moduli when the thermal field is \textit{interpolated}.}
            \label{f:footbridge_temp_linear_1_a}
        \end{subfigure}
        \begin{subfigure}[t]{\textwidth}
            \centering
            \includegraphics[trim=0 730 0 730, clip,width=0.9\textwidth]{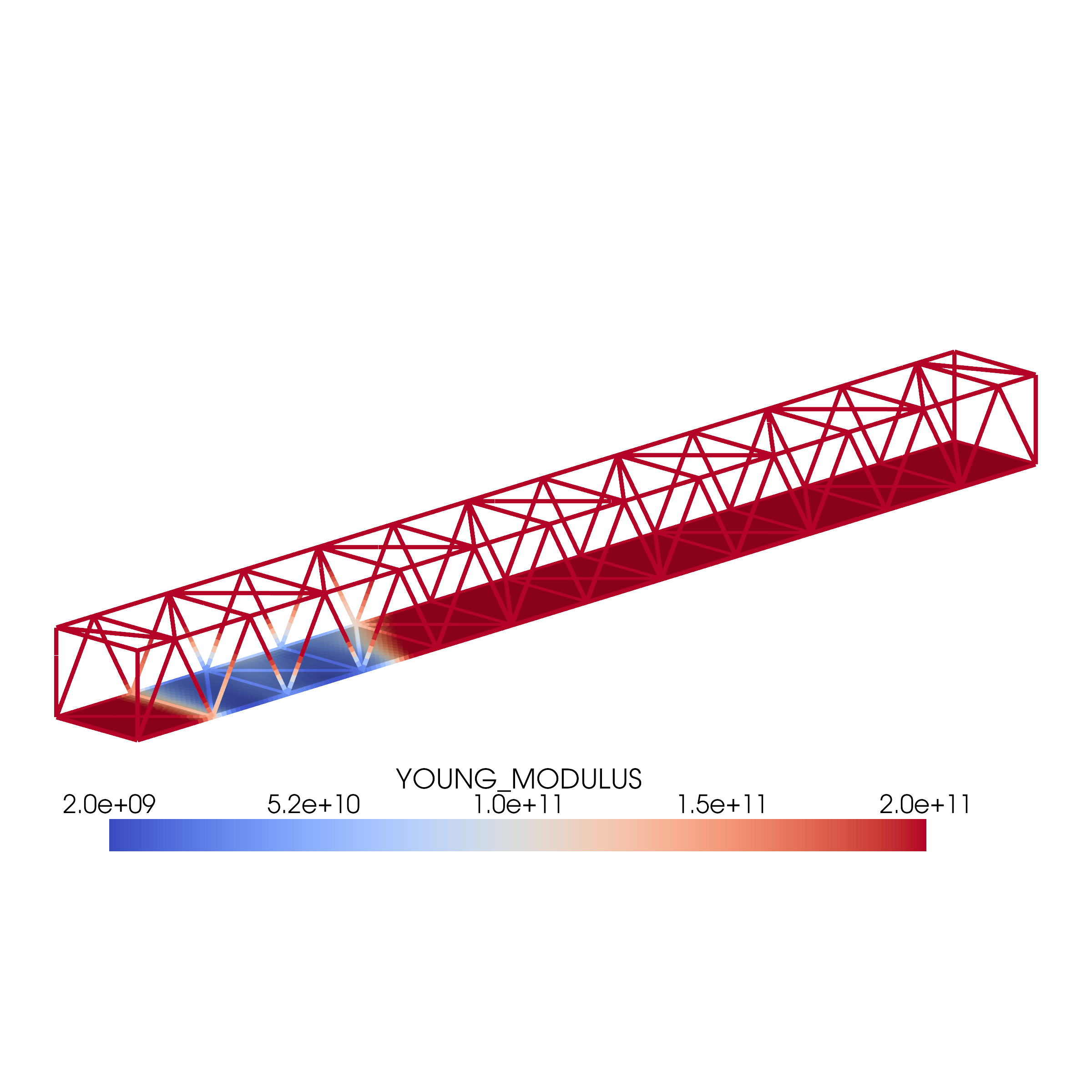}
            \caption{Identified Young's moduli when the thermal field is reconstructed during \acrshort{SI}: \textit{Monolithic} approach.}
            \label{f:footbridge_temp_linear_1_b}
        \end{subfigure}
        \begin{subfigure}[t]{\textwidth}
            \centering            
            \includegraphics[trim=0 730 0 730, clip,width=0.9\textwidth]{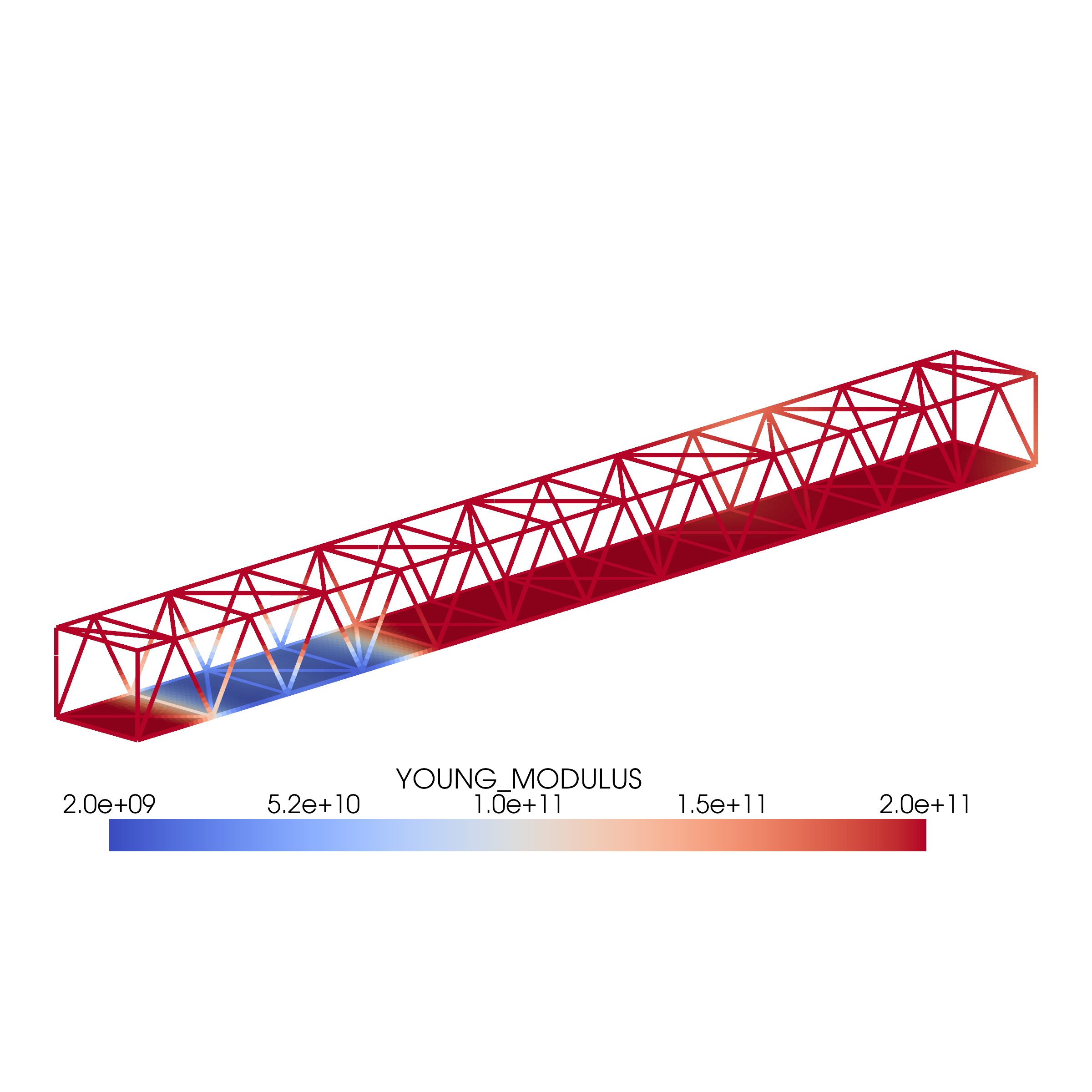}
           \caption{Identified Young's moduli when the thermal field is reconstructed during \acrshort{SI}: \textit{Partitioned} approach.}
            \label{f:footbridge_temp_linear_1_c}
        \end{subfigure}
        \vspace{-0.4em}
        \begin{subfigure}[t]{\textwidth}
            \centering            
            \includegraphics[trim=0 530 0 1631, clip,width=0.9\textwidth]{Figure_24c.png}
        \end{subfigure}
    \end{minipage}
    \hfill
    \caption{\textbf{Footbridge. Linearly varying thermal field.} Identified Young's moduli distributions when the thermal load is considered using different approaches during \acrshort{SI}.}
\label{f:footbridge_temp_linear_1}
\end{figure} 

\begin{figure}[!b]
    \centering
    \begin{minipage}[t]{\textwidth}
        \centering
        \begin{subfigure}[t]{\textwidth}
            \centering
            \includegraphics[trim=0 730 0 730, clip,width=0.9\textwidth]{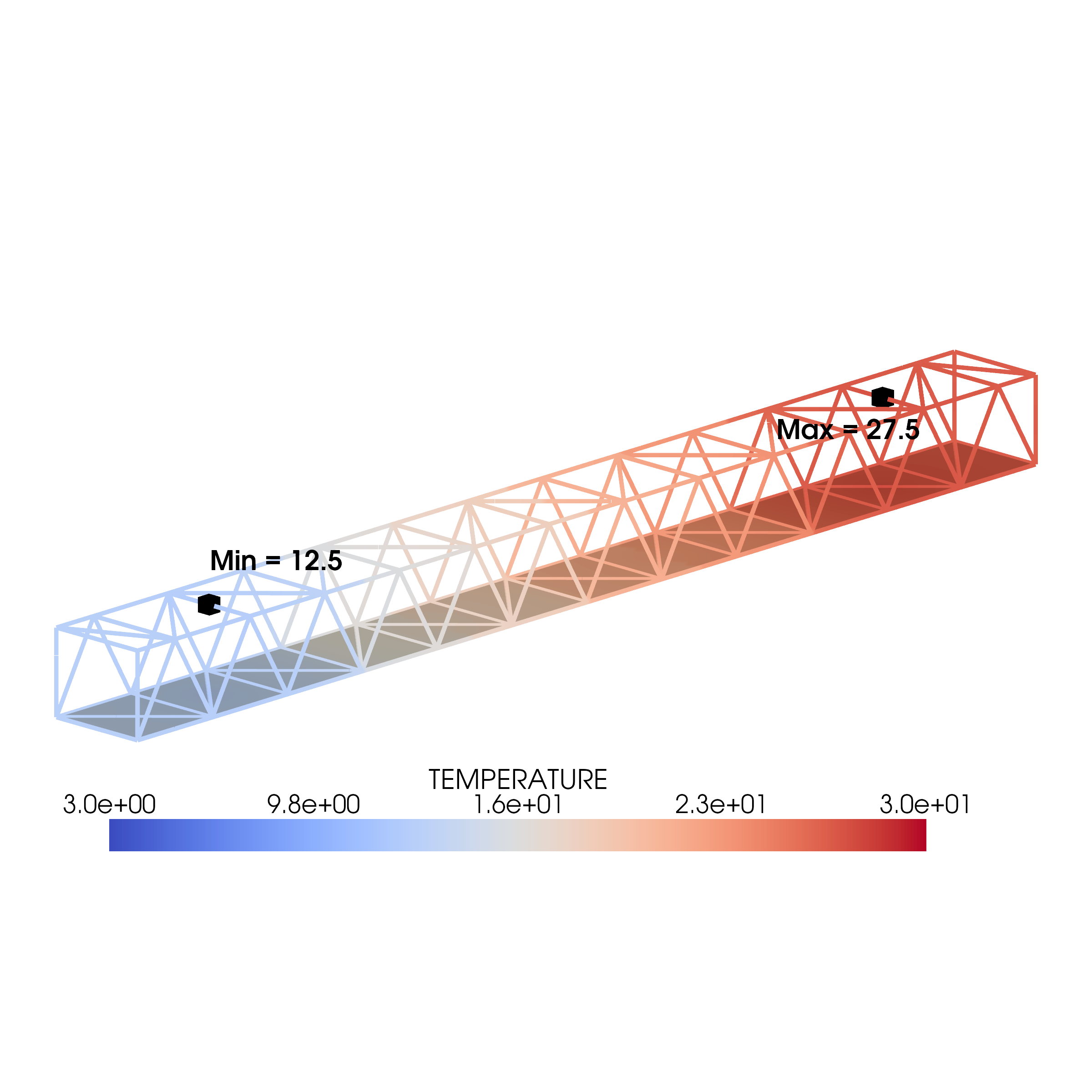}
            \caption{Interpolated temperature distribution when thermal field is \textit{interpolated} during \acrshort{SI}.}
            \label{f:footbridge_temp_linear_2_a}
        \end{subfigure}
        \begin{subfigure}[t]{\textwidth}
            \centering
            \includegraphics[trim=0 730 0 730, clip,width=0.9\textwidth]{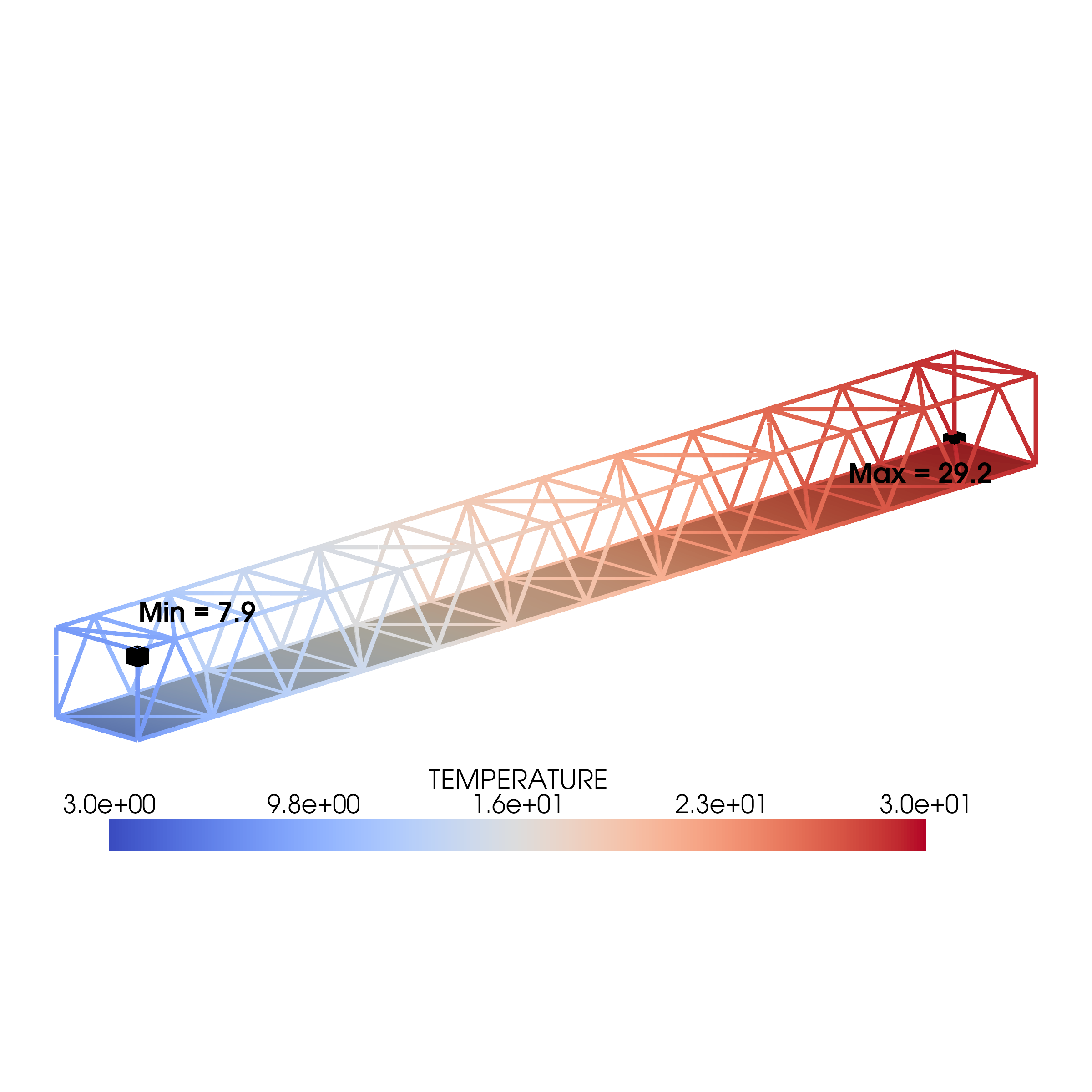}
            \caption{Identified temperature distribution when the thermal field is reconstructed during \acrshort{SI}: \textit{Monolithic} approach.}
            \label{f:footbridge_temp_linear_2_b}
        \end{subfigure}
        \begin{subfigure}[t]{\textwidth}
            \centering            
            \includegraphics[trim=0 730 0 730, clip,width=0.9\textwidth]{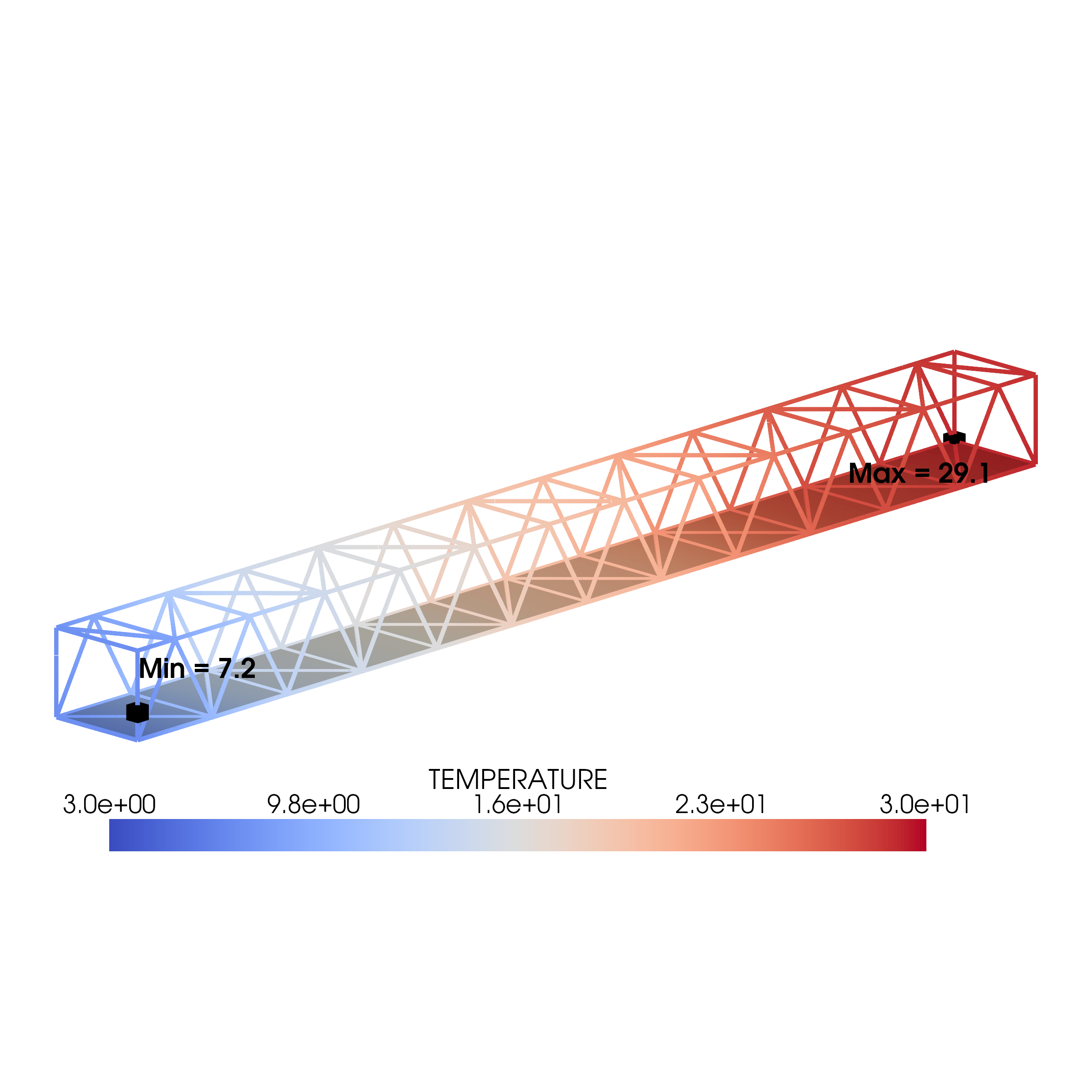}
           \caption{Identified temperature distribution when the thermal field is reconstructed during \acrshort{SI}: \textit{Partitioned} approach.}
           \label{f:footbridge_temp_linear_2_c}
        \end{subfigure}
      \vspace{-0.4em}
        \begin{subfigure}[t]{\textwidth}
            \centering            
            \includegraphics[trim=0 530 0 1631, clip,width=0.9\textwidth]{Figure_25c.png}
        \end{subfigure}
    \end{minipage}
    \hfill
    \caption{\textbf{Footbridge. Linearly varying thermal field.} Temperature field distributions when the thermal load is considered using different approaches during \acrshort{SI}. Peak temperatures are noted in the figures. }
\label{f:footbridge_temp_linear_2}
\end{figure} 

\begin{figure}[!t]
\centering
\begin{minipage}[t]{0.49\linewidth}
\centering
\includegraphics[width=\linewidth]{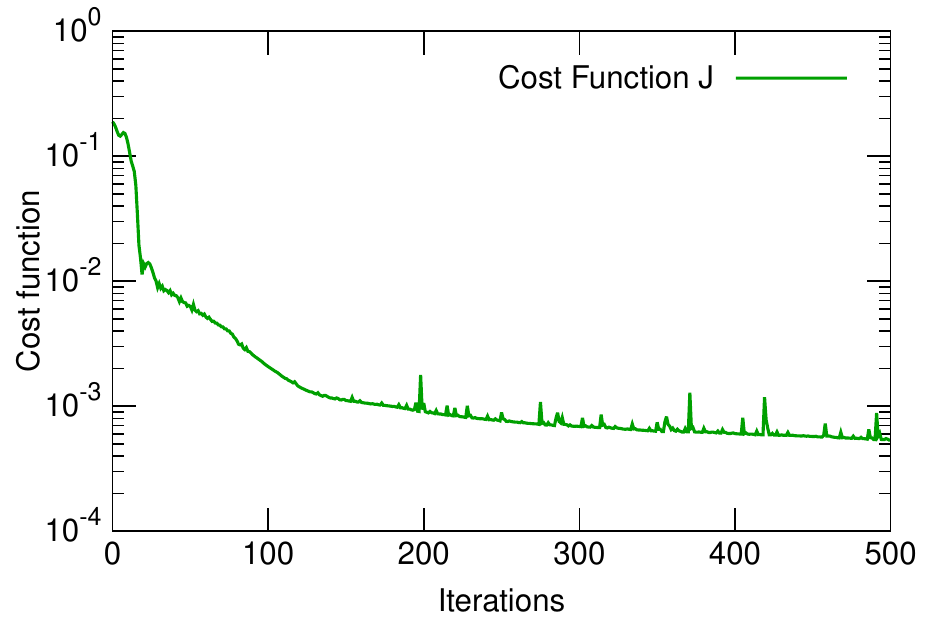}
\subcaption{When the temperature field is \textit{interpolated}.}
\label{f:footbridge_temp_linear_3_a}
\end{minipage}
\hfill
\begin{minipage}[t]{0.49\linewidth}
\centering
\includegraphics[width=\linewidth]{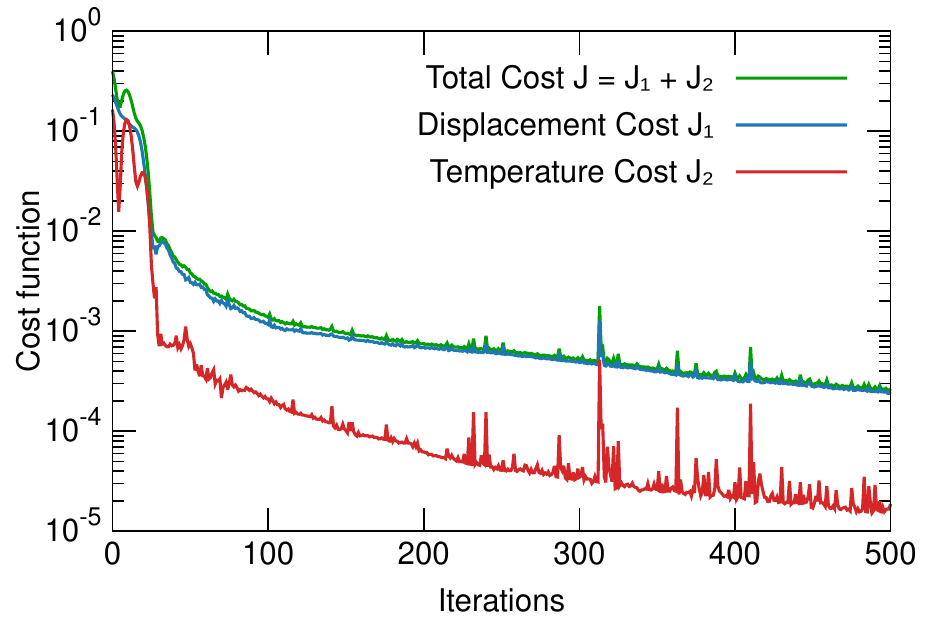}
\subcaption{When the temperature field is \textit{identified}: \textit{Monolithic} approach.}
\label{f:footbridge_temp_linear_3_b}
\end{minipage}
\begin{minipage}[t]{0.49\linewidth}
\centering
\includegraphics[width=\linewidth]{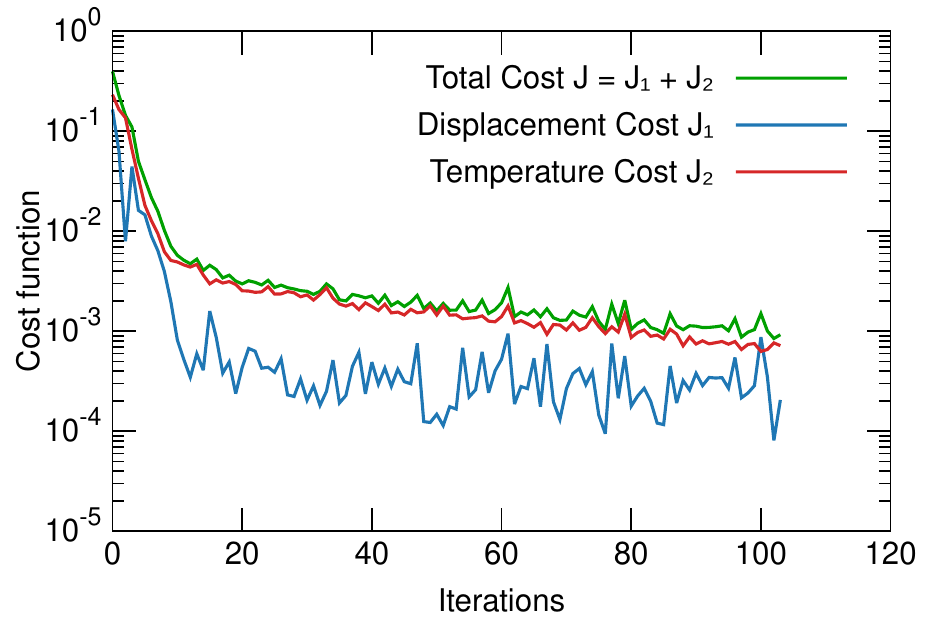}
\subcaption{When the temperature field is \textit{identified}: \textit{Partitioned} approach. Coupling iteration: 103, Overall iteration: 996.}
\label{f:footbridge_temp_linear_3_c}
\end{minipage}
\caption{\textbf{Footbridge. Linearly varying thermal field.}
Convergence plots when the thermal load is considered using different approaches during \acrshort{SI}.}
\label{f:footbridge_temp_linear_3}
\end{figure}

The results for Scenario 3 (temperature interpolation) and Scenario 4 (temperature identification) are presented together for easy comparison. They are grouped by the target temperature distribution.


\textbf{Case 1: Linearly varying thermal field.} 

Figures \ref{f:footbridge_temp_linear_1_a} and \ref{f:footbridge_temp_linear_2_a} show the identified Young's modulus distribution and temperature distribution for Scenario 3 (\acrshort{kNN} interpolated temperature), for the linearly varying thermal field case. It can be observed that a notable improvement in the Young's moduli identification is observed over Scenario 2 (constant temperature field; Fig. \ref{f:footbridge_w_const_temp_a}). Some of the false damage identifications have been cleared, leading to marginally better damage localization. In the temperature field, interpolation reconstructs a fairly adequate approximation of the thermal distribution in the range $[12.5, 27.5]$ \degree C. 
The corresponding convergence plot is depicted in Fig. \ref{f:footbridge_temp_linear_3_a}, where it can be observed that the cost function has a higher reduction than Scenario 2 (Fig. \ref{f:footbridge_w_const_temp_c}) and is in the order of $\mathcal{O}(10^{-3})$.

Figures \ref{f:footbridge_temp_linear_1_b} and \ref{f:footbridge_temp_linear_1_c} show the identified Young's modulus distribution for Scenario 4 (temperature field identified), for the linearly varying thermal field case. Results are shown for the monolithic and partitioned approaches, respectively.
Compared to Scenario 3, the Young's modulus distribution is even better localized for both the monolithic and partitioned approaches. However, some minimal false damage is observed in the partitioned approach identification.
The corresponding temperature distributions are illustrated in Figs. \ref{f:footbridge_temp_linear_2_b} and \ref{f:footbridge_temp_linear_2_c}
The identified temperature ranges are $[7.9, 29.2]$ \degree C and $[7.2, 29.1]$ \degree C for the monolithic and partitioned approaches, respectively. It can be noted that the temperature ranges are larger compared to the interpolation case.
The corresponding convergence plots for the monolithic and partitioned approaches are depicted in Figs. \ref{f:footbridge_temp_linear_3_b} and \ref{f:footbridge_temp_linear_3_c}, respectively. Similar to the plate with a hole example, the partitioned approach plot is plotted over the outer coupling iterations, while the monolithic approach plot is plotted over the optimization iterations. 
Also, a contrasting dominance of the cost function contributions is noticeable: the displacement error dominates in the monolithic approach, and the temperature error dominates in the partitioned approach, despite visually similar fields. The composite cost functions in both approaches are in the order of $\mathcal{O}(10^{-3})$.

\begin{figure}[!b]
    \centering
    \begin{minipage}[t]{\textwidth}
        \centering
        \begin{subfigure}[t]{\textwidth}
            \centering
            \includegraphics[trim=0 730 0 730, clip,width=0.9\textwidth]{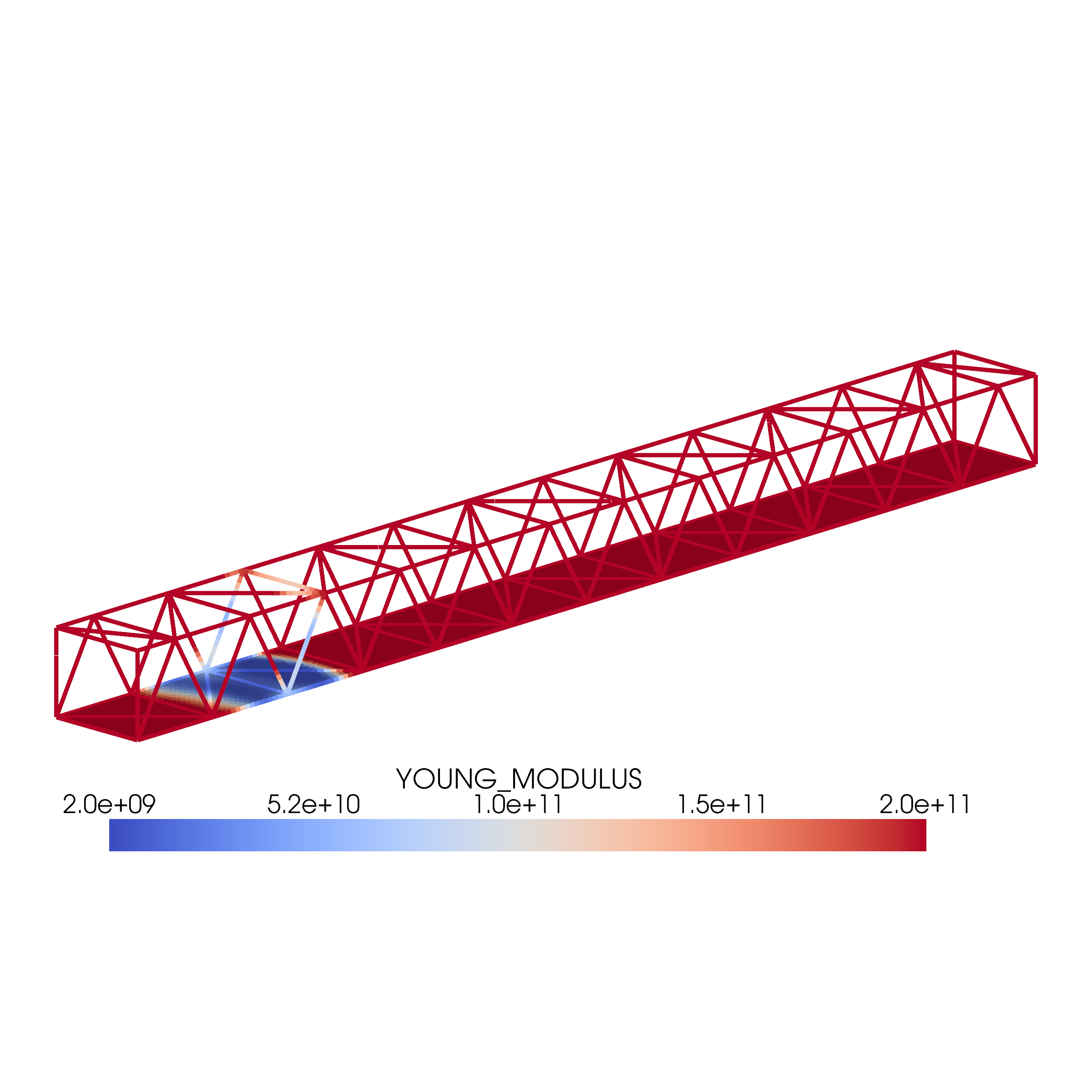}
            \caption{Identified Young's moduli when the thermal field is \textit{interpolated}.}
            \label{f:footbridge_temp_gaussian_1_a}
        \end{subfigure}
        \begin{subfigure}[t]{\textwidth}
            \centering
            \includegraphics[trim=0 730 0 730, clip,width=0.9\textwidth]{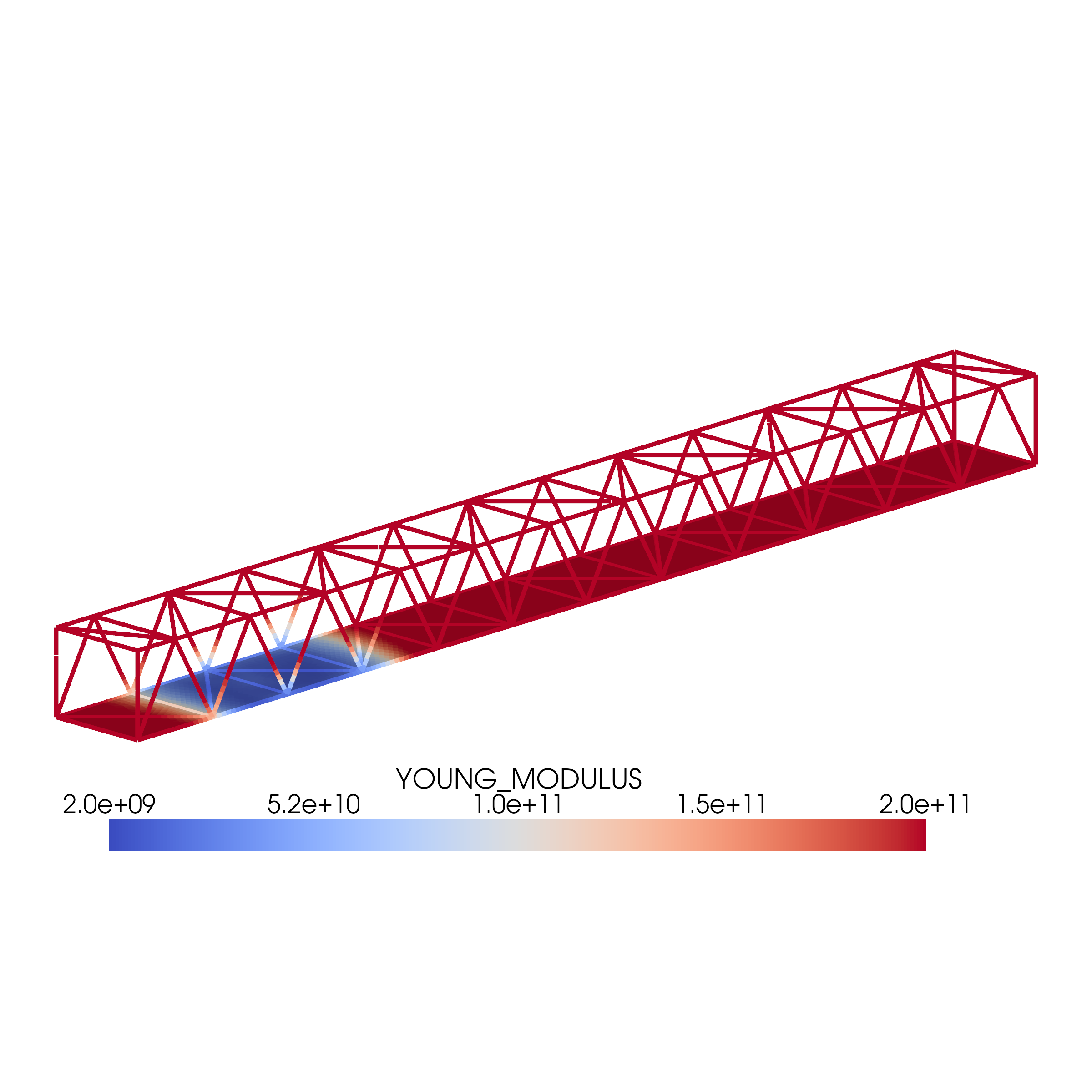}
            \caption{Identified Young's moduli when the thermal field is reconstructed during \acrshort{SI}: \textit{Monolithic} approach.}
            \label{f:footbridge_temp_gaussian_1_b}
        \end{subfigure}
        \begin{subfigure}[t]{\textwidth}
            \centering            
            \includegraphics[trim=0 730 0 730, clip,width=0.9\textwidth]{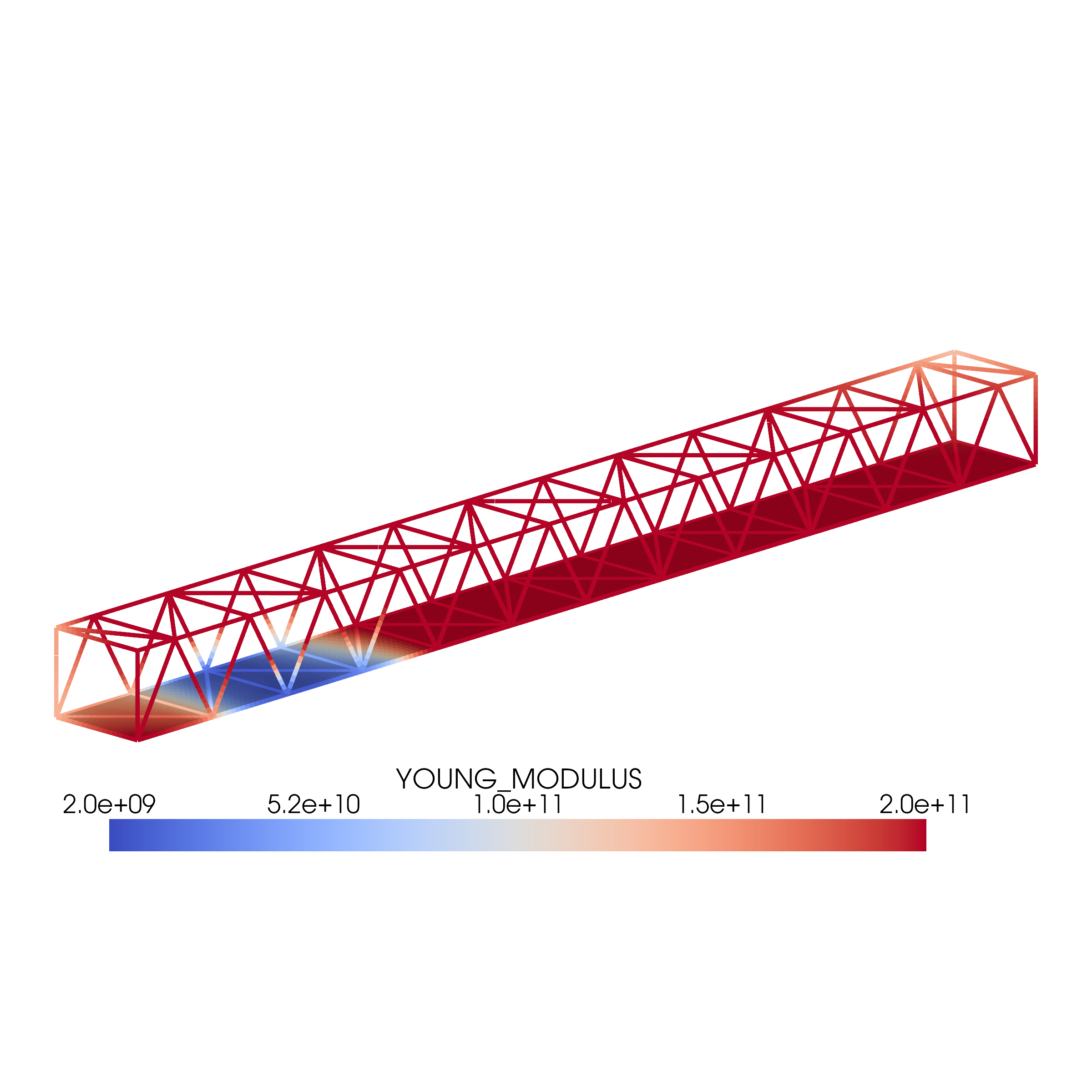}
           \caption{Identified Young's moduli when the thermal field is reconstructed during \acrshort{SI}: \textit{Partitioned} approach.}
           \label{f:footbridge_temp_gaussian_1_c}
        \end{subfigure}
        \vspace{-0.4em}
        \begin{subfigure}[t]{\textwidth}
            \centering            
            \includegraphics[trim=0 530 0 1631, clip,width=0.9\textwidth]{Figure_27c.png}
        \end{subfigure}
    \end{minipage}
    \hfill
    \caption{\textbf{Footbridge. Localized thermal field.} Identified Young's moduli distributions when the thermal load is considered using different approaches during \acrshort{SI}.}
\label{f:footbridge_temp_gaussian_1}
\end{figure} 

\begin{figure}[!b]
    \centering
    \begin{minipage}[t]{\textwidth}
        \centering
        \begin{subfigure}[t]{\textwidth}
            \centering
            \includegraphics[trim=0 730 0 730, clip,width=0.9\textwidth]{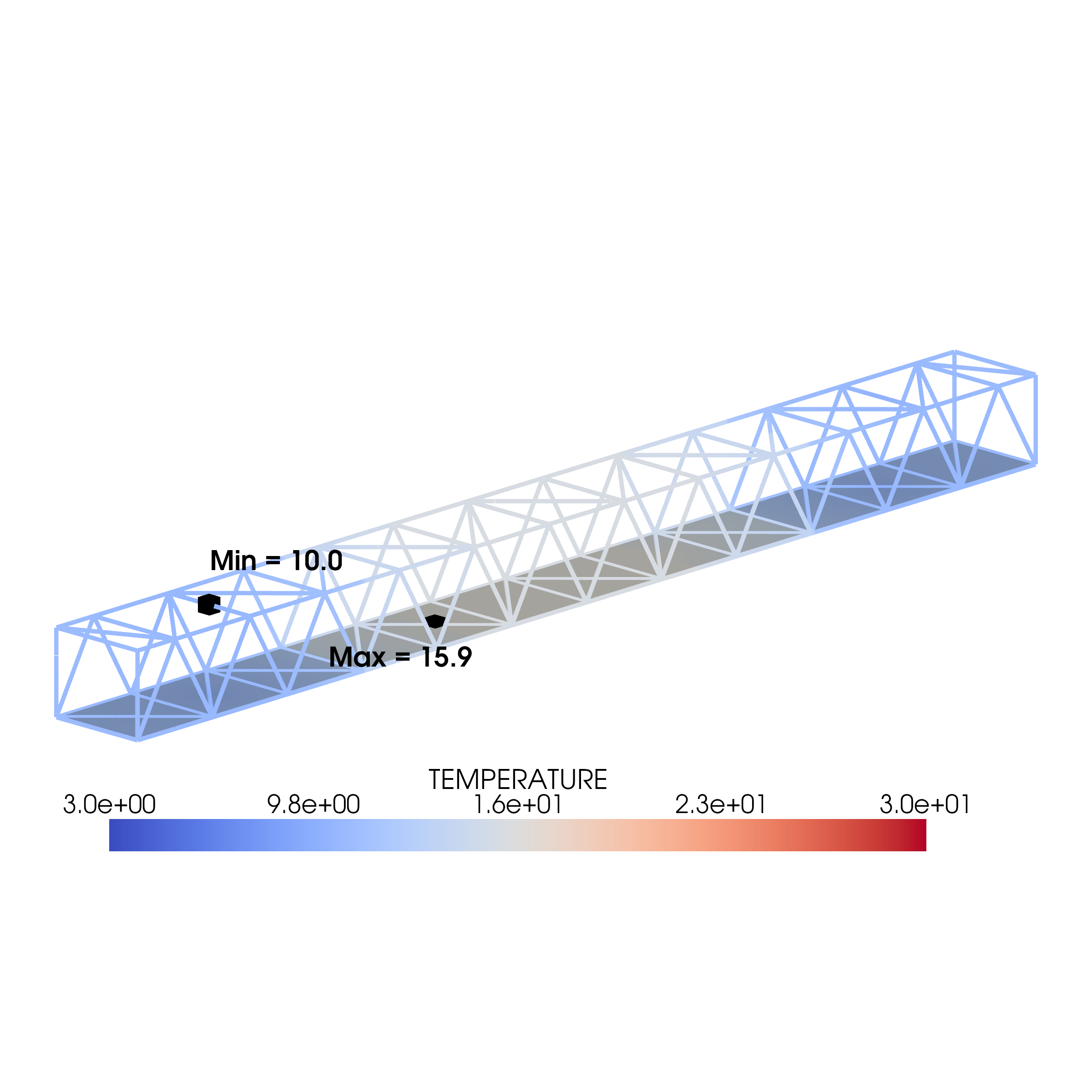}
            \caption{Interpolated temperature distribution when thermal field is \textit{interpolated} during \acrshort{SI}.}
            \label{f:footbridge_temp_gaussian_2_a}
        \end{subfigure}
        \begin{subfigure}[t]{\textwidth}
            \centering
            \includegraphics[trim=0 730 0 730, clip,width=0.9\textwidth]{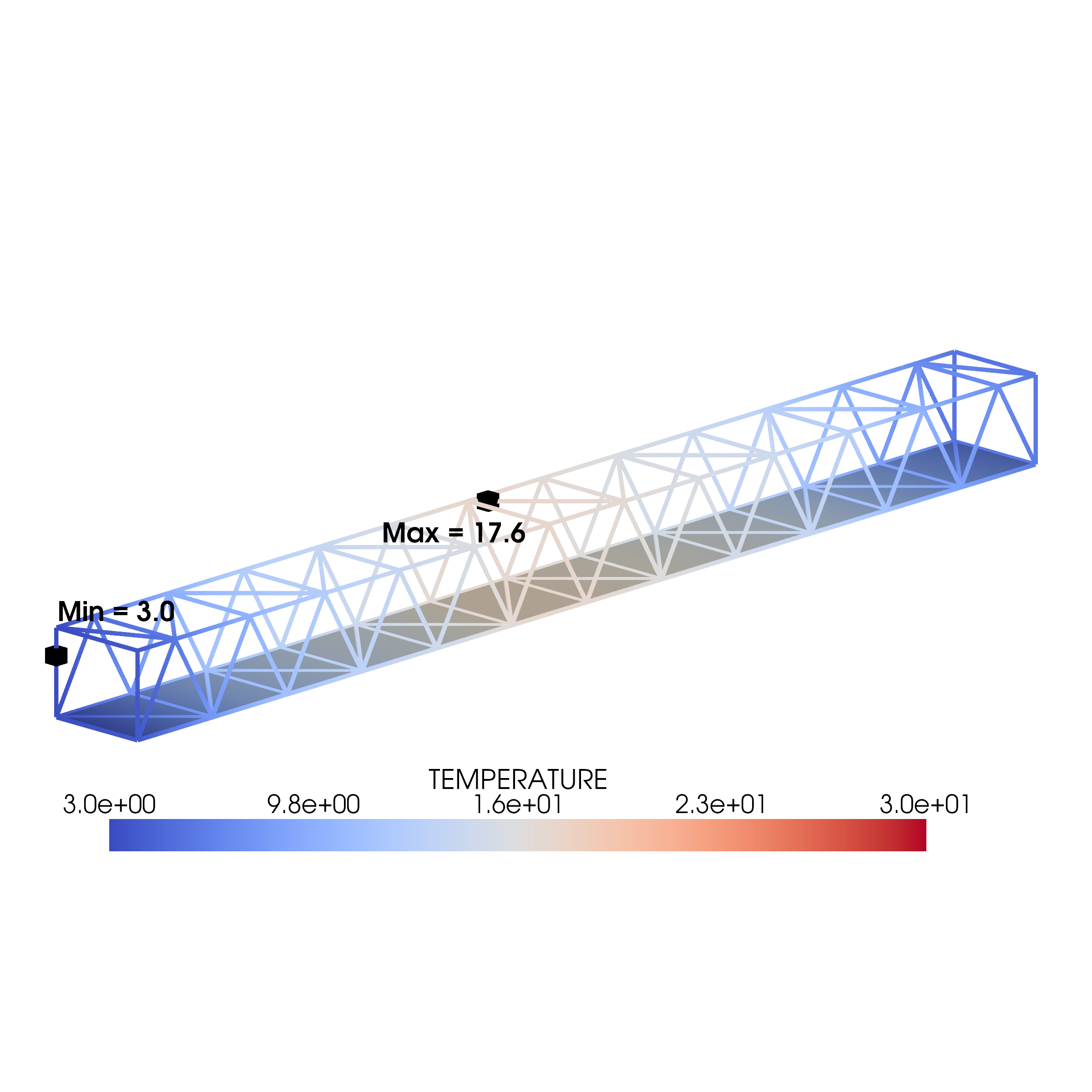}
            \caption{Identified temperature distribution  when the thermal field is reconstructed during \acrshort{SI}: \textit{Monolithic} approach.}
            \label{f:footbridge_temp_gaussian_2_b}
        \end{subfigure}
        \begin{subfigure}[t]{\textwidth}
            \centering            
            \includegraphics[trim=0 730 0 730, clip,width=0.9\textwidth]{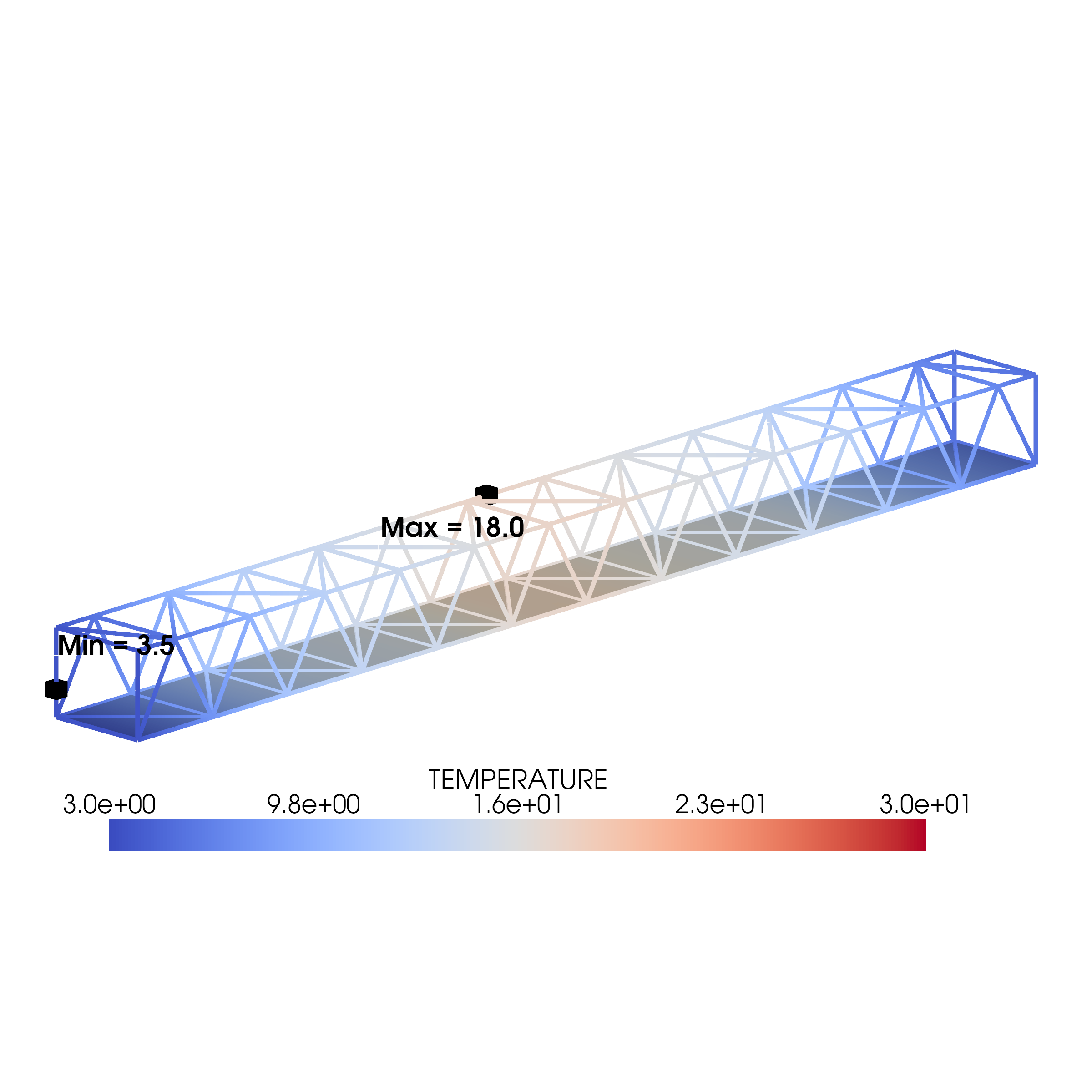}
           \caption{Identified temperature distribution  when the thermal field is reconstructed during \acrshort{SI}: \textit{Partitioned} approach.}
           \label{f:footbridge_temp_gaussian_2_c}
        \end{subfigure}
        \vspace{-0.4em}
        \begin{subfigure}[t]{\textwidth}
            \centering            
            \includegraphics[trim=0 530 0 1631, clip,width=0.9\textwidth]{Figure_28c.png}
        \end{subfigure}
    \end{minipage}
    \hfill
    \caption{\textbf{Footbridge. Localized thermal field.} Temperature field distributions when the thermal load is considered using different approaches during \acrshort{SI}. Peak temperatures are noted in the figures. }
\label{f:footbridge_temp_gaussian_2}
\end{figure} 

\begin{figure}[!h]
\centering
\begin{minipage}[t]{0.49\linewidth}
\centering
\includegraphics[width=\linewidth]{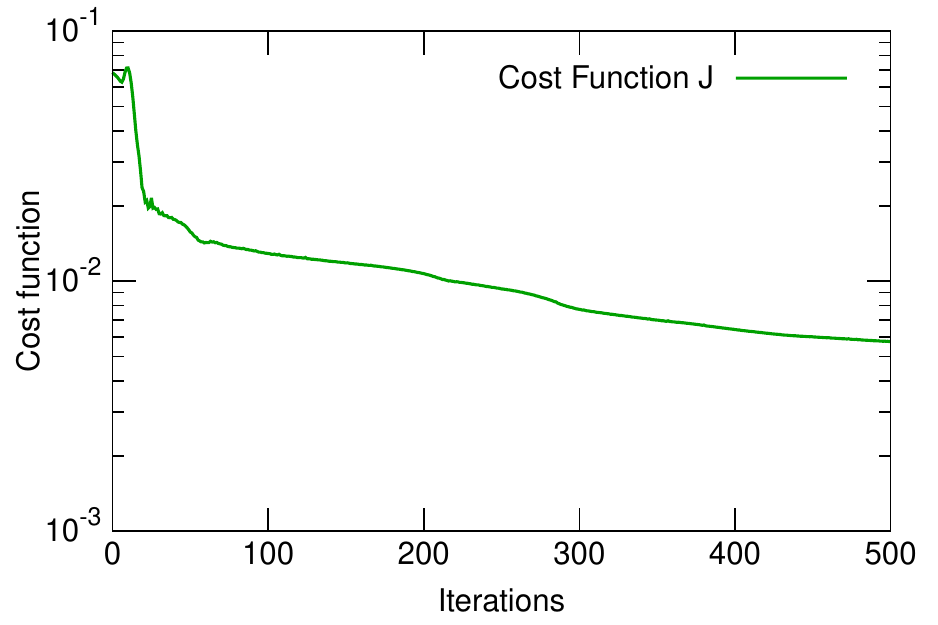}
\subcaption{When the temperature field is \textit{interpolated}.}
\label{f:footbridge_temp_gaussian_3_a}
\end{minipage}
\hfill
\begin{minipage}[t]{0.49\linewidth}
\centering
\includegraphics[width=\linewidth]{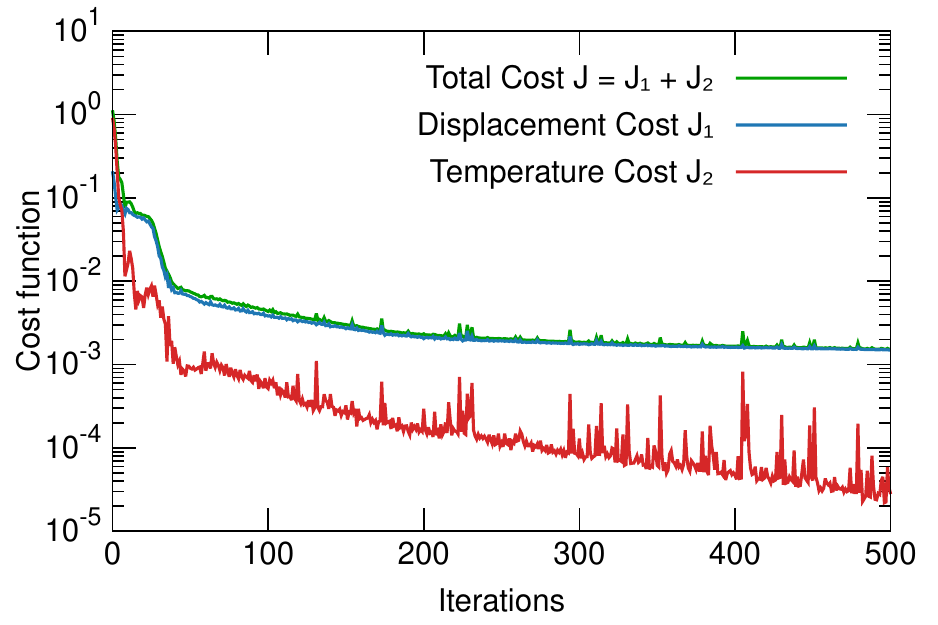}
\subcaption{When the temperature field is \textit{identified}: \textit{Monolithic} approach.}
\label{f:footbridge_temp_gaussian_3_b}
\end{minipage}
\begin{minipage}[t]{0.49\linewidth}
\centering
\includegraphics[width=\linewidth]{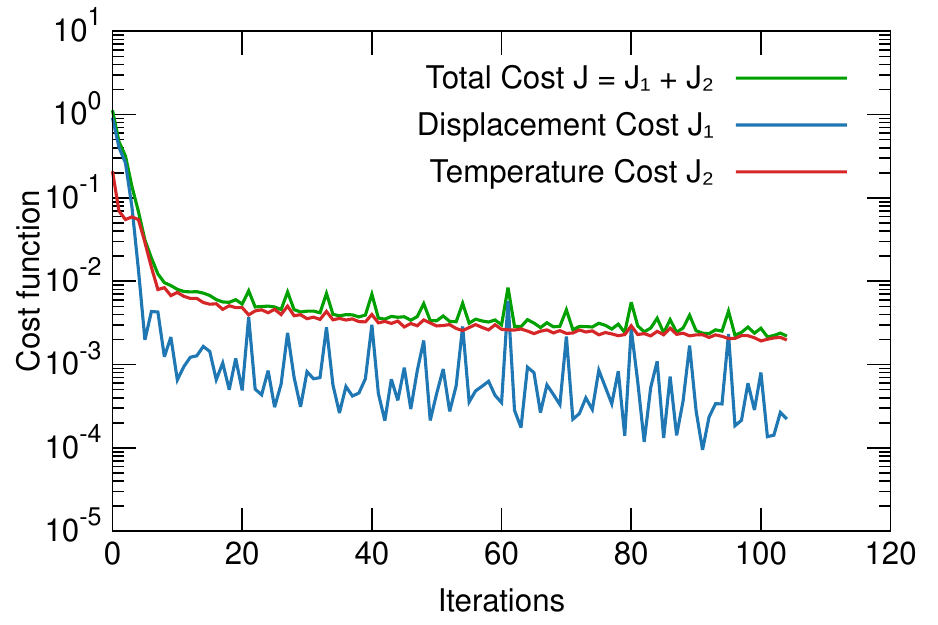}
\subcaption{When the temperature field is \textit{identified}: \textit{Partitioned} approach. Coupling iteration: 104, Overall iteration: 997.}
\label{f:footbridge_temp_gaussian_3_c}
\end{minipage}
\caption{\textbf{Footbridge. Localized thermal field.}
Convergence plots when the thermal load is considered using different approaches during \acrshort{SI}.}
\label{f:footbridge_temp_gaussian_3}
\end{figure}



\textbf{Case 2: Localized thermal field.} 

Figures \ref{f:footbridge_temp_gaussian_1_a} and \ref{f:footbridge_temp_gaussian_2_a} show the identified Young's modulus distribution and temperature distribution for Scenario 3 (\acrshort{kNN} interpolated temperature), for the localized thermal field case. 
Compared to the Scenario 2 (constant temperature field; Fig. \ref{f:footbridge_w_const_temp_a}) result, the interpolation results yield incredible enhancement in the damage localization quality. Apart from some false damage on a few side and top beam members, the damage localization is very accurate. 
Although not perfect, the interpolation does loosely detect the localized pattern in the thermal field with interpolated temperatures in the range $[10.0, 15.9]$ \degree C. 
The corresponding convergence plot is depicted in Fig. \ref{f:footbridge_temp_gaussian_3_a}, where it can be observed that even with just one order of magnitude reduction in the cost function, a very good approximation of the Young's modulus distribution is obtained.

Figures \ref{f:footbridge_temp_gaussian_1_b} and \ref{f:footbridge_temp_gaussian_1_c} show the identified Young's modulus for Scenario 4 (temperature field identified), for the localized thermal field case. Results are shown for the monolithic and partitioned approaches, respectively.
Compared to Scenario 3, the damage localization shows more smearing, but the false damage is eliminated in the monolithic case. In the partitioned case, most of the identified damage is concentrated at the correct location; however, some false damages are also observed.
The corresponding temperature distributions are illustrated in Figs. \ref{f:footbridge_temp_gaussian_2_b} and \ref{f:footbridge_temp_gaussian_2_c}.
The identified temperature ranges are $[3.0, 17.6]$ \degree C and $[3.5, 18.0]$ \degree C for the monolithic and partitioned approaches, respectively. Compared to the interpolation case, the distributions have higher maximum peak temperatures and lower minimum peak temperatures.
The corresponding convergence plots for the monolithic and partitioned approaches are depicted in Figs. \ref{f:footbridge_temp_gaussian_3_b} and \ref{f:footbridge_temp_gaussian_3_c}, respectively. 
In this case as well, the displacement error dominates the composite cost function in the monolithic case, and the temperature error dominates in the partitioned case. 
The composite cost functions in both approaches are in the order of $\mathcal{O}(10^{-2})$.




\vspace{5em}

\textbf{Quantitative Comparison}

\begin{table}[!b]
\caption{\textbf{Footbridge.}  Identification errors in Young’s modulus distributions for Scenarios \textbf{1 and 2} where thermal effects are ignored or accounted for using a constant temperature field. Percentage changes are computed relative to the errors at optimization start. }
\label{tab:Table_4}
\begin{tabularx}{\textwidth}{Xccc}
\hline \hline
\textbf{\begin{tabular}[c]{@{}c@{}} Approach for \\ system identification    \end{tabular}} 
& \textbf{\begin{tabular}[c]{@{}c@{}} Type of \\ thermal field \end{tabular}}  & \textbf{\begin{tabular}[c]{@{}c@{}} \boldsymbol{$\epsilon_{L_2}$}   {[}-{]} \end{tabular}}  &  \textbf{\begin{tabular}[c]{@{}c@{}} \boldsymbol{$\delta\epsilon_{L_2}$}   {[}\%{]} \end{tabular}} \\ 
\hline
\multicolumn{4}{c}{\textbf{\begin{tabular}[c]{@{}c@{}}Scenario 1: Thermal load not accounted for during \acrshort{SI}\end{tabular}}}      \\ \hline
\multirow{2}{=}{\acrshort{SI} without considering temperature}          & Linearly varying      & $3.623\cdot10^{-1}$  & $1.084$   \\ \cline{2-4}    & Localized     &$ 3.623\cdot10^{-1}$  & $1.084$     \\ \hline        &     &    &    \\ \hline
\multicolumn{4}{c}{\textbf{\begin{tabular}[c]{@{}c@{}}Scenario 2: Thermal load accounted as a constant \\ temperature field during \acrshort{SI}\end{tabular}}}     \\ \hline
\multirow{2}{=}{\acrshort{SI} considering constant $20$\degree C temperature field} & Linearly varying   & $4.466\cdot10^{-1}$	& $24.599$
      \\ \cline{2-4} & Localized   & $4.892\cdot10^{-1}	$& $36.479$
     \\ \hline \hline
\end{tabularx}
\end{table}

A quantitative comparison based on the relative discrete \acrshort{$L_2$} error ($\epsilon_{L_2}$) and its percentage change ($\delta \epsilon_{L_2}$) is presented in this example as well. 
At the start of system identification, the relative discrete \acrshort{$L_2$} error in Young's modulus distribution was $3.584\cdot10^{-1}$, in the temperature field for the linearly varying thermal field was $2.829\cdot10^{-1}$, while for the localized thermal field was $5.732\cdot10^{-1}$. These errors are considered as the reference errors to evaluate the performance of different approaches and also serve as the baseline to compute percentage change in the relative discrete \acrshort{$L_2$} errors ($\delta \epsilon_{L_2}$) for each case.

Table \ref{tab:Table_4} tabulates the $\epsilon_{L_2}$ and $\delta\epsilon_{L_2}$ for the Young’s modulus distributions for Scenarios 1 and 2, where thermal effects were unaccounted, and accounted for using a constant $20$ \degree C temperature field. 
The superscripts `$\mathbf{E}$' and `$\boldsymbol{\Delta\mathbf{T}}$' in $\epsilon$ and $\delta\epsilon$ refer to errors in the Young's modulus and temperature fields.

From Scenario 1, it can be observed that failing to account for thermal effects during \acrshort{SI} results in higher errors than at the start of the optimization. The increase in $\epsilon_{L_2}$ is not that high numerically, because during \acrshort{SI}, the Young's modulus throughout the structure increases from the initial value of $1.98\cdot10^{11}$ and stagnates at the bound of $2\cdot10^{11}$. This small change results in a small increase in the $\epsilon_{L_2}$. However, the overarching consequence of this is that it leads to complete failure to detect any damage on the structure as observed in the visual examination in Fig. \ref{f:footbridge_wo_temp}.

In Scenario 2, where a simplistic (constant $20$ \degree C thermal field) approach to account for thermal effects was applied, the errors increased and were much higher than in Scenario 1. From the convergence plots for this case in Fig. \ref{f:footbridge_w_const_temp}, it was observed that the cost function reduced a few magnitudes, which happens because the algorithm, unaware of the correct thermal field, attempts to reduce the displacement sensor errors by changing the Young's modulus distribution. Therefore, even though the optimization cost function reduces, the $\epsilon_{L_2}$, which compares the identified field to the target field, increases. 

Table \ref{tab:Table_5} tabulates the $\epsilon_{L_2}$ and $\delta\epsilon_{L_2}$ for the Young’s modulus and temperature distributions, where the temperature sensor information was incorporated to account for thermal effects, i.e., Scenario 3 (via temperature interpolation) and Scenario 4 (via temperature identification). The table is grouped by thermal field type. 

For the linearly varying thermal field case, the interpolation case shows some improvement in the Young's modulus identification ($\delta \epsilon_{L_2}^{E}=-15.907$\%), which can be attributed to the better approximation of the thermal field ($\delta \epsilon_{L_2}^{T}=-80.212$\%). The monolithic and partitioned approaches in the identification case demonstrated much better accuracy than the interpolation case for Young's modulus identification with reductions in $\delta \epsilon_{L_2}^{E}$ of over $40$\% and reductions in $\delta \epsilon_{L_2}^{T}$ of over $86$\% in the temperature distribution errors. 
These values are consistent with the visual examination in Figs. \ref{f:footbridge_temp_linear_1} and \ref{f:footbridge_temp_linear_2}.

\begin{table}[!t]
\caption{\textbf{Footbridge.} Identification Errors in Young's modulus (\textbf{E}) and temperature ($\boldsymbol{\Delta \mathbf{T}}$) fields for Scenarios (\textbf{3 and 4}) \textbf{incorporating temperature sensor information} to account for thermal effects. Percentage changes are computed relative to the errors at optimization start. }
\label{tab:Table_5}
\begin{tabularx}{\textwidth}{Xcccc}
\hline \hline
\multirow{2}{=}{\parbox{\linewidth}{\centering
\textbf{Approach accounting}\\
\textbf{for thermal effects}}}                                        & \multicolumn{2}{c}{\textbf{\begin{tabular}[c]{@{}c@{}} \boldsymbol{$\epsilon_{L_2}$}   {[}-{]}\end{tabular}}} & \multicolumn{2}{c}{\textbf{\begin{tabular}[c]{@{}c@{}} \boldsymbol{$\delta\epsilon_{L_2}$}   {[}\%{]}\end{tabular}}} \\ \cline{2-5}   & \textbf{E}  & $\boldsymbol{\Delta \mathbf{T}}$   & \textbf{E}  & $\boldsymbol{\Delta \mathbf{T}}$   \\ \hline
\multicolumn{5}{c}{\textbf{Case 1: Linearly varying thermal field}}       \\ \hline
temperature interpolation   & $3.014\cdot10^{-1}$   & $5.598\cdot10^{-2}$    & $-15.907 $   & $-80.212$    \\ \hline
temperature identification   &    &     &     &   \\
\multicolumn{1}{c}{Monolithic approach}   & $1.987\cdot10^{-1}$  & $2.501\cdot10^{-2}$  & $-44.557$  & $-91.161$    \\ 
\multicolumn{1}{c}{Partitioned approach}    & $2.047\cdot10^{-1} $        & $3.868\cdot10^{-2}$  & $-42.881 $    & $-86.329$    \\ \hline & \multicolumn{1}{c}{}  & \multicolumn{1}{c}{}      & \multicolumn{1}{c}{}  & \multicolumn{1}{c}{}   \\ \hline
\multicolumn{5}{c}{\textbf{Case 2: Localized thermal field}}  \\ \hline
temperature interpolation    & $1.839\cdot10^{-1}$    & $3.237\cdot10^{-1}$    & $-48.690$   & $-43.524$   \\ \hline
temperature identification     &    &   &   &    \\
\multicolumn{1}{c}{Monolithic approach}    & $1.698\cdot10^{-1}$   & $3.209\cdot10^{-1}$   & $-52.620$  & $-44.011 $  \\ 
\multicolumn{1}{c}{Partitioned approach}  & $1.664\cdot10^{-1} $   & $3.220\cdot10^{-1}$   & $-53.565$  & $-43.834$  \\ \hline \hline
\end{tabularx}
\end{table}

For the localized thermal field case, the interpolation case yielded good approximations for the Young's modulus and temperature distributions with over $40$\% reductions in $\delta \epsilon_{L_2}^{E}$ and $\delta \epsilon_{L_2}^{T}$. 
The monolithic and partitioned approaches in the identification case outperformed the interpolation case with marginally better identifications of Young's modulus and temperature fields. These values corroborate the observations in Figs. \ref{f:footbridge_temp_gaussian_1} and \ref{f:footbridge_temp_gaussian_2}.

Between the two thermal field cases, the linearly varying thermal field case saw almost double the reduction in $\delta \epsilon_{L_2}^{T}$ than the localized thermal field case. This indicates that the current 8 temperature sensors and their positioning are better suited to capture the linearly varying thermal field than the localized thermal field used here. This can also be inferred by looking at the sensor locations in Fig. \ref{f:footbridge_temp_sensors}.

\FloatBarrier

\section{Conclusion}
\label{s:conclusion}

In this work, an optimization-driven adjoint-based high-fidelity system identification procedure for localizing weakness and identifying the temperature field in one-way thermo-mechanical coupled structures using displacement and temperature measurements is proposed. This work highlights the need to accurately account for thermal effects on structures during \acrshort{SI}. The underlying framework formulates the identification problem as an optimization problem that minimizes the aggregated weighted errors between the measured and simulated responses. 
Results show that failing to account for thermal effects during \acrshort{SI} leads to catastrophic damage localization with many false positives and, even more severely, to complete localization failure (i.e., no damage is detected). 

Therefore, building upon the standard \acrshort{SI} method, this work proposed two alternative approaches, namely, monolithic and partitioned approaches, to account for thermal effects within the \acrshort{SI} framework for two-field (Young’s modulus and temperature fields) identification in one-way thermo-mechanical coupled structures, where temperature affects the mechanical response. 
The \textit{monolithic} approach formulates the \acrshort{SI} problem as a single global optimization problem that minimizes the combined errors from displacement and temperature measurements and the model response, and simultaneously identifies the Young's modulus and temperature distributions. 
On the other hand, the \textit{partitioned} approach divides the \acrshort{SI} problem into two sub-optimizations, each responsible for updating one control field, and uses a Gauss-Seidel-type fixed-point iteration loop on top to couple these two sub-optimizations. The sub-optimizations are solved inexactly with loose convergence, accumulating incremental refinements to both fields and converging to the overall solutions.

The proposed methodology was tested using two numerical examples of a Plate With a Hole and a Footbridge \acrshort{FE} model. The Plate With a Hole example was designed to study the effect of the number and locations of temperature sensors on \acrshort{SI} and was investigated for 6- and 16-temperature sensor configurations, with a fixed 14-displacement sensor configuration. The Footbridge example was intended to examine the methodology's applicability to a realistic, fine-meshed, mixed (beams and shells) \acrshort{FE} model, and was investigated for an 8-temperature and 8-displacement sensor configuration. Both examples were tested on linearly varying and localized (Gaussian-type) thermal load cases. 

For all cases, `target' Young's modulus and temperature distributions were prescribed, together with external mechanical loads. The `measurement data' were obtained by extracting the displacements and temperatures at the measurement locations by simulating this target model. The structure was assumed to be in a near-pristine (undamaged) Young's modulus condition, and an initial temperature distribution was assumed at the start of the optimization. The optimization algorithm, then, depending on the \acrshort{SI} approach (monolithic or partitioned), minimized sensor errors by tuning the control fields (Young's modulus and temperature distributions) to recover the prescribed target distributions. 
In this work, different optimization algorithms with different step size methods were used to demonstrate independence from the algorithm type. 
It is well established in the literature that these high-fidelity inverse problems are highly ill-conditioned due to a much smaller number of measurement points and a much larger number of design variables. This issue is further compounded in these two fields, Young's modulus and temperature fields identification, because not only can different combinations of variables within one field produce the same model response, but also different combinations of variables between the two fields also produce the same response. This can lead the optimization algorithm to a completely wrong local optimum, even when the cost function is similar to that of the correct target solution. To help mitigate this issue, Vertex Morphing filtering was used in this work, which can be understood as a convolution operation on variables and gradients to smooth the discrete, noisy solutions. 

In all cases, both the proposed monolithic and partitioned \acrshort{SI} approaches accurately identified the Young's modulus distribution (i.e., localized damage) and the temperature distribution (i.e., thermal load) in the structure.

In general, it was observed that the linearly varying thermal field was identified with better accuracy ($\epsilon_{L_2}^{T}=\mathcal{O}(10^{-3}-10^{-2})$) than the localized thermal field ($\epsilon_{L_2}^{T}=\mathcal{O}(10^{-1})$). This can be attributed to the linear trend being more easily captured than the localized trend, especially when the heated zone is devoid of temperature sensors. Additionally, \acrshort{VM} filtering smooths and smears the fields, thus losing the capability to capture sharp gradients in the solution fields.

For the Plate With Hole example, across the 6- and 16-temperature sensor configurations, the 16-sensor configuration yielded better results than the 6-sensor configuration for the linearly varying thermal field. However, the 16-sensor configuration underperformed compared to the 6-sensor configuration for the localized (Gaussian-like) thermal field. Thus, having more temperature sensors yields better results when they collectively capture the thermal field more accurately. If the sensors are placed such that they do not accurately capture the features of the thermal field, then even more sensors do not yield better results, as observed in the 16-sensor configuration with a localized thermal field where the middle heated zone is devoid of temperature sensors. 

Across both examples, the performance of the monolithic and partitioned approaches was very similar, with minor differences due to ill-conditioning. However, in most cases, Young's modulus identification using the monolithic approach was marginally better than that using the partitioned approach, especially in reducing false damage identifications. 

It was interesting to observe that even though the solution fields from both approaches were very similar, visually and quantitatively ($\epsilon_{L_2}$), the cost function convergence showed a stark contrast in trend. In the monolithic cases, the cost function was dominated by the displacement error term, which was 1-3 orders of magnitude higher than the temperature error term. On the contrary, in the partitioned cases, the cost function was dominated by the temperature error term, which was 1-3 orders of magnitude higher than the temperature error term. This is the clearest evidence of an ill-conditioned \acrshort{SI}, where both approaches, despite very similar costs and solution fields, have landed in different local optima.

The results from the proposed methodology were also compared against other methods of accounting for thermal effects. In the simplest case, a constant-temperature distribution was assumed throughout the structure. In general, some improvement was observed compared to the case in which thermal effects were entirely unaccounted for, but major deficiencies remained, with large regions of false damage identifications. 
To obtain a better estimate of the thermal field, temperature sensor measurements were used to interpolate it throughout the structure. A major improvement was observed in the interpolation cases; however, there were still major deficiencies in the localized thermal field cases, where the damage localization was poor. 

The proposed methodology proved most beneficial in cases where the temperature sensors, due to their number and/or location, were unable to fully capture the underlying trend in the temperature distribution through interpolation. 
Despite the proposed monolithic and partitioned approaches successfully identifying the Young's modulus and temperature distributions in this one-way thermo-mechanical coupled \acrshort{SI}, further research is required to address several open questions:
\begin{itemize}
    \item Extending the methodology to strongly coupled multi-physics \acrshort{SI} problems, such as fluid-structure-interaction, and fluid-structure-thermal-interaction,
    \item Effect of relative weighting of the displacement cost and temperature cost terms in the composite cost function on the \acrshort{SI} solution,
    \item Effect of relaxation on the Gauss-Seidel-based partitioned \acrshort{SI} approach,
    \item Exploring strategies to promote \textit{balanced convergence} of multiple identified fields, preventing dominance of individual cost terms during \acrshort{SI},
    \item Extending the methodology to account for other relevant environmental effects such as wind, rain, cloud cover, etc.    
\end{itemize}


\section*{Declaration of competing interest}

The authors declare that they have no known competing financial interests or personal relationships that could have appeared to influence the work reported in this paper.

\section*{DATA AVAILABILITY STATEMENT}

All the examples shown in the article are clearly defined and do not require any further data for reproducibility.

\section*{CRediT authorship contribution statement}

Conceptualization, T.A., S.W., R.L. and R.W.; methodology, T.A., S.W., and I.A.; software, T.A. and S.W.; validation, T.A. and I.A.; formal analysis and investigation, T.A.; data curation, T.A. and S.W.; writing---original draft preparation, T.A.; writing---review and editing, T.A., S.W., I.A., H.A., R.L., and R.W.; visualization, T.A.; supervision, H.A., R.L., and R.W.; project administration, R.L. and R.W.; funding acquisition, R.L. and R.W.

\section*{Acknowledgments}
This research was funded by SIEMENS AG and the Technical University of Munich, Institute for Advanced Study, Lichtenbergstrasse 2a, D-85748, Garching, Germany. The GMU team is also partially supported by the Office of Naval Research (Award No: N00014-24-1-2147), NSF grant DMS-2408877, and the Air Force Office of Scientific Research (Award No: FA9550-25-1-0231). We would like to acknowledge support from the Leonhard Obermeyer Center, Technical University of Munich, and funding by the Deutsche Forschungsgemeinschaft (DFG, German Research Foundation) under Germany´s Excellence Strategy – EXC 2163/1 - Sustainable and Energy Efficient Aviation – Project-ID 390881007. The authors acknowledge Dr. Facundo Airaudo for valuable discussions related to this work.

\bibliographystyle{elsarticle-harv} 
\bibliography{references,ref2,ref3}

\end{document}